\pdfoutput=1
\RequirePackage{ifpdf}
\ifpdf 
\documentclass[pdftex]{sigma}
\else
\documentclass{sigma}
\fi

\usepackage{amsfonts}
\usepackage[all]{xy}

\usepackage{bbm}

\numberwithin{equation}{section}

\newtheorem{Theorem}{Theorem}[section]
\newtheorem*{Theorem*}{Theorem}
\newtheorem{Corollary}[Theorem]{Corollary}
\newtheorem{Lemma}[Theorem]{Lemma}

 { \theoremstyle{definition}

\newtheorem{Remark}[Theorem]{Remark} }

\newcommand{\rsdraw}[3]{\raisebox{-#1\height}{\scalebox{#2}{\includegraphics{#3}}}}
\newcommand{\co}{\colon}
\newcommand{\vareps}{\varepsilon}
\newcommand{\tens}{\otimes}
\newcommand{\id}{\mathrm{id}}
\newcommand{\idrm}{\mathrm{id}}
\newcommand{\fs}[1]{[#1]}
\newcommand{\Zen}[1]{\mc {Z}(#1)}
\newcommand{\ltr}[1]{\text{tr}_l(#1)}
\newcommand{\rtr}[1]{\text{tr}_r(#1)}
\newcommand{\mc}{\mathcal}
\newcommand{\opp}{\mathrm{op}}
\newcommand{\N}{\mathbb{N}}
\newcommand{\R}{\mathbb{R}}

\newcommand{\Hom}{\mathrm{Hom}}
\newcommand{\lra}{\leftrightarrow}
\newcommand{\End}{\text{End}}
\newcommand{\Cyl}{\mathrm{Cyl}}
\newcommand{\Mod}{\mathrm{Mod}}

\newcommand{\Ob}[1]{\mathrm{Ob}(#1)}
\newcommand{\Cob}{\textbf{Cob}}
\newcommand{\conncob}{\textbf{3Cob}_1}
\newcommand{\kk}{\Bbbk}
\newcommand{\uu}{\mathbbm{1}}
\newcommand{\evl}{\text{ev}}
\newcommand{\evr}{\widetilde{\text{ev}}}
\newcommand{\coevl}{\text{coev}}
\newcommand{\coevr}{\widetilde{\text{coev}}}
\newcommand{\RT}{\text{RT}}
\newcommand{\TV}{\text{TV}}

\newcommand{\aend}{\mathbb{A}}
\newcommand{\coend}{\mathbb{F}}
\newcommand{\coendzc}{\mathbb{G}}
\newcommand{\Reid}[1]{\textbf{R}#1}
\newcommand{\rep}[1]{\text{rep}_{#1}}
\newcommand{\ltrn}{\triangleleft}
\newcommand{\rtrn}{\triangleright}
\newcommand{\RomanNumeralCaps}[1]
{\MakeUppercase{\romannumeral #1}}
\newlength{\gnat}

\begin{document}

\allowdisplaybreaks

\newcommand{\arXivNumber}{2306.07216}

\renewcommand{\PaperNumber}{074}

\FirstPageHeading

\ShortArticleName{Cyclic Objects from Surfaces}

\ArticleName{Cyclic Objects from Surfaces}

\Author{Ivan BARTULOVI\'{C}}

\AuthorNameForHeading{I.~Bartulovi\'{c}}

\Address{Institut f\"{u}r Geometrie, Zellescher Weg 12-14, 01062 Dresden, Germany}
\Email{\href{mailto:ivan.bartulovic@tu-dresden.de}{ivan.bartulovic@tu-dresden.de}}
\URLaddress{\url{https://sites.google.com/view/ivan-bartulovic}}

\ArticleDates{Received May 27, 2024, in final form August 20, 2025; Published online September 06, 2025}

\Abstract{In this paper, we endow the family of all closed genus $g \ge 1$ surfaces with a~structure of a (co)cyclic object in the category of 3-dimensional cobordisms. As a corollary, any $3$-dimensional TQFT induces a (co)cyclic module, which we compute algebraically for the Reshetikhin--Turaev TQFT.}


\Classification{18N50; 57K16}

\section{Introduction}

\subsection{Background motivation and the main result} In this paper, we study cyclic objects and their interplay with topological quantum field theories.
A (co)cyclic object in a category is, roughly speaking, a (co)simplicial object with
compatible actions of the cyclic groups. Cyclic homology of algebras was independently introduced by Connes~\cite{connes_non-commutative_1985} and Tsygan~\cite{tsygan}. To any algebra over a commutative ring $\kk$ is associated a certain (co)cyclic $\kk$-module, that is, a (co)cyclic object in the category of $\kk$-modules. The (co)faces and the (co)degeneracies of the associated (co)cyclic module are induced respectively by the multiplication and the unit of the algebra and the (co)cyclic operators are given by cyclic permutations on tensor products.
This construction was generalized to the braided setting by Akrami and Majid~\cite{cycliccocycles}, who associate a cocyclic $\kk$-module to any ribbon algebra in a braided monoidal category.

Motivated by the Hopf-algebraic study of ribbon string links (which are like framed pure braids, but they can double back on itself) from~\cite{bv}, the author of the present paper equipped in \cite{bartulovic2022cyclicsetsribbonstring} the set of isotopy classes of ribbon string links with a structure of a cocyclic set by setting\looseness=1
\begingroup
\allowdisplaybreaks
\begin{gather*}
\delta_0^n(T)= \,
\put(10,-30){\small $1$}
\put(40,-30){\small $n$}
\put(25,-2.0){\small $T$}
\put(22,17){\small $\cdots$}
\put(22,-22){\small $\cdots$}
\raisebox{-11.5mm}{\includegraphics{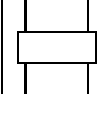}}\,,
\qquad
\delta_i^n(T)=\,
\put(12.5,-2.0){\small $T$}
\put(2.6,-30){\small $1$}
\put(25,-30){\small $i$}
\put(40.3,-30){\small $i+1$}
\put(69.4,-30){\small $n$}
\put(10.5,17){\small $\cdots$}
\put(10.5,-22){\small $\cdots$}
\put(55.5,17){\small $\cdots$}
\put(55.5,-22){\small $\cdots$}
\raisebox{-11.5mm}{\includegraphics{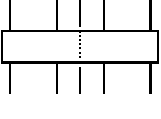}}\,,
\qquad
\delta_n^n(T)=\,
\put(17,-2.0){\small $T$}
\put(2.6,-30){\small $1$}
\put(31.7,-30){\small $n$}
\put(14.5,17){\small $\cdots$}
\put(14.5,-22){\small $\cdots$}
\raisebox{-11.5mm}{\includegraphics{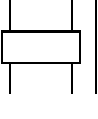}}\,, \\
\sigma_j^n(T)=\,
\put(22.5,-1.3){\small $T$}
\put(2.4,-28.4){\small $0$}
\put(30,17){\small $j$}
\put(36,-28.5){\small $j+1$}
\put(63,-28.5){\small $n+1$}
\put(10,19){\small $\cdots$}
\put(10,-21){\small $\cdots$}
\put(55.5,19){\small $\cdots$}
\put(55.5,-21){\small $\cdots$}
\raisebox{-12.2mm}{\includegraphics{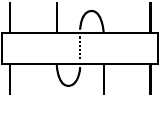}}\,,
\qquad
\tau_n(T)=\,
\put(26,-1.3){\small $T$}
\put(2.4,-28.4){\small $0$}
\put(21,-28.5){\small $n-1$}
\put(48,-28.5){\small $n$}
\put(14,20){\small $\cdots$}
\put(14,-21){\small $\cdots$}
\raisebox{-11.5mm}{\includegraphics{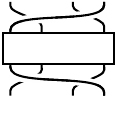}}\,,
\end{gather*}
\endgroup
and of a cyclic set by defining
\begin{gather*}
d_i^n(T)=\,
\put(22,-1.3){\small $T$}
\put(2.35,-28.6){\small $0$}
\put(36,-28.6){\small $i$}
\put(69,-28.6){\small $n$}
\put(10.5,20){\small $\cdots$}
\put(10.5,-21){\small $\cdots$}
\put(55.5,20){\small $\cdots$}
\put(55.5,-21){\small $\cdots$}
\raisebox{-11.5mm}{\includegraphics{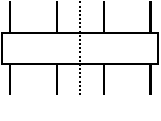}}\,,\qquad
s_j^n(T)=\,
\put(37,-1.3){\small $T$}
\put(2.35,-28.6){\small $0$}
\put(37,-28.6){\small $j$}
\put(73,-28.6){\small $n$}
\put(10.5,20){\small $\cdots$}
\put(10.5,-21){\small $\cdots$}
\put(59,20){\small $\cdots$}
\put(59,-21){\small $\cdots$}
\raisebox{-12.3mm}{\includegraphics{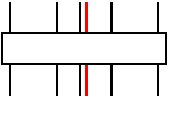}}\,,\qquad
t_n(T)=\,
\put(25,-1.3){\small $T$}
\put(2.5,-28.6){\small $0$}
\put(18,-28.6){\small $1$}
\put(47,-28.6){\small $n$}
\put(29,20){\small $\cdots$}
\put(29,-21){\small $\cdots$}
\raisebox{-11.2mm}{\includegraphics{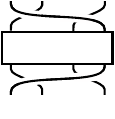}}\,.
\end{gather*}
In the latter, the face operator $d_i^n$ is given by deleting the component labeled by $i$ and the degeneracy operator~$s_j^n$ is given by duplicating (along the framing) the component labeled by~$j$ of a given string link.
By mimicking the construction of the cocyclic module associated to so called ribbon algebras given by Akrami and Majid in~\cite{cycliccocycles}, one can write the (co)cyclic module for (co)algebras in a balanced category. In parti\-cular, this applies to the Lyubashenko coend $\coend$ of a ribbon category $\mc B$, which is a Hopf algebra object in $\mc B$ and which features in constructions of quantum invariants of links and $3$-manifolds~\cite{lyubashenko1995invariants}.
 It is also shown in \cite{bartulovic2022cyclicsetsribbonstring} that the quantum invariants à la Reshetikhin--Turaev from \cite{bv} form a morphism from the (co)cyclic set from ribbon string links to the (co)cyclic set \smash{$\bigl\{\Hom_{\mc B}\bigl(\coend^{\tens n+1}, \uu\bigr)\bigr\}_{n\in \N}$} associated to the coend $\coend$ of $\mc B$ (using its underlying (co)algebra structure). In this way, the (co)cyclic sets of geometric inspiration are initial among (co)cyclic sets derived from ribbon categories.

On the other hand, compact surfaces are relatively well-understood $2$-manifolds which appear in many areas of mathematics. In particular, closed oriented surfaces are objects of a symmetric monoidal category of $3$-dimensional cobordisms $\textbf{3Cob}_0$. A morphism between two surfaces is given by a homeomorphism class of $3$-cobordisms between the given surfaces. Introduced by Atiyah~\cite{At}, a $3$-dimensional topological quantum field theory (or shortly, TQFT) is a strong symmetric monoidal functor from $\textbf{3Cob}_0$ to the category of modules over a commutative ring~$\kk$.
 A fundamental construction of a~$3$-dimensional TQFT in this sense is the Reshetikhin--Turaev TQFT~\cite{reshetikhin_invariants_1991,turaevqinvariants}.
Its main algebraic ingredient is a modular category $\mc B$ (see Section~\ref{modular}), which is in particular $\kk$-linear, braided (not necessarily symmetric) and semisimple.
The $\kk$-module associated to a surface of genus $n$ (often called the state space) is isomorphic to~$\Hom_{\mc B}\bigl(\coend^{\tens n}, \uu\bigr)$.
By $\kk$-linearity of $\mc B$, the above-mentioned (co)cyclic set~\smash{$\bigl\{\Hom_{\mc B}\bigl(\coend^{\tens n+1}, \uu\bigr)\bigr\}_{n\in \N}$} forms in fact a~(co)cyclic $\kk$-module.
The motivating question of this paper was whether there is a (co)cyclic object in the category of~$3$-dimensional cobordisms, which is sent by the Reshetikhin--Turaev functor to the (co)cyclic $\kk$-module~\smash{$\bigl\{\Hom_{\mc B}\bigl(\coend^{\tens n+1}, \uu\bigr)\bigr\}_{n\in \N}$}.
The main results of this paper (see Theorems~\ref{CYCINCOB} and \ref{CYCINCOBRT}) answer this question positively.
As a main corollary, any $3$-dimensional TQFT induces a (co)cyclic $\kk$-module.
Also, we discuss some potentially related work in the setting of the category of connected cobordisms $\textbf{3Cob}_1$, which first appeared in~\cite{kerler} and is different from the cobordism category used throughout the paper. For instance, it is a non-symmetric braided category.
In this setting, we outline a construction of the so called para(co)cyclic objects associated to the one-holed torus (see Sections \ref{parastuff} and \ref{cy para}). By composition, the braided monoidal functor $J_3$ from \cite{beliakova-derenzi} induces para(co)cyclic objects associated to the end of unimodular ribbon factorizable categories.

\subsection{Organization of the paper}
The paper is organized as follows.
In Section~\ref{cyclicprelimini}, we recall the notion of a (co)cyclic object in a~category.
In Section~\ref{seccobs}, we recall some facts about $3$-cobordisms and their presentation via special ribbon graphs. In Section~\ref{cobcyclic}, we construct (co)cyclic objects in the category of $3$-dimensional cobordisms.
In Section~\ref{catprelimini}, we review ribbon categories and their graphical calculus, braided Hopf algebras, and related concepts.
Section~\ref{algebraiccyclic} is dedi\-cated to (co)cyclic modules from categorical (co)algebras.
In Section~\ref{tftsrel}, we relate, via the Reshetikhin--Turaev TQFT, the (co)cyclic objects from surfaces with (co)cyclic modules associated to the coend of an anomaly free modular category. In Section~\ref{related}, we discuss paracyclic objects in the category of connected cobordisms.

\subsection{Notation}
Unless otherwise stated, by $\kk$ we denote any commutative ring. The class of objects in a~category~$\mc C$ is denoted by~$\Ob{\mc C}$.
By~$\mathbb{N}$ we denote the set of natural numbers including zero and we put~$\mathbb{N}^* = \mathbb{N} \setminus \{0\}$.

\section{Cyclic objects} \label{cyclicprelimini}

In this section, we recall the notions of (co)simp\-li\-cial and (co)cyclic objects in a category.

\subsection{The simplicial category}
The \emph{simplicial category}~$\Delta$ is defined as follows.
The objects of~$\Delta$ are the non-negative integers. For $n\in \N$, denote $\fs{n}=\{0, \dots, n\}$.
A morphism~$n\to m$ in~$\Delta$ is an increasing map between sets~$\fs{n}$ and~$\fs{m}$.
For~$n\in \N^*$ and~$0\le i \le n$, the~$i$-th \emph{coface}~$\delta_i^n \co n-1 \to n$ is the unique increasing injection from~$\fs{n-1}$ into~$\fs{n}$ which misses~$i$.
For~$n\in \N$ and~$0\le j \le n$, the~$j$-th~\emph{codegene\-racy}~$\sigma_j^n \co n+1 \to n$ is the unique increasing surjection from~$\fs{n+1}$ onto~$\fs{n}$ which sends both~$j$ and~$j+1$ to~$j$.

It is well known (see~\cite[Lemma 5.1]{MCLHomology}) that morphisms in~$\Delta$ are generated by cofaces $\{\delta_i^n\}_{n\in \N^*, 0\leq i \leq n}$ and codegeneracies~\smash{$\bigl\{\sigma_j^n\bigr\}_{n\in \N, 0\leq j \leq n}$} subject to the \textit{simplicial relations}:
\begingroup
\allowdisplaybreaks
\begin{gather}
\label{cofaces} \delta_j^{n+1} \delta_i^n= \delta_i^{n+1} \delta_{j-1}^n \hspace{7.9mm} \text{for all }0 \le i<j \le n+1, \\
\label{codegeneracies} \sigma_j^{n}\sigma_i^{n+1}=\sigma_i^n\sigma_{j+1}^{n+1} \hspace{9.2mm} \text{for all } 0\le i \leq j \le n,
\\
\label{compcofcod}\sigma_j^n\delta_i^{n+1}= \begin{cases}
\delta_i^n \sigma_{j-1}^{n-1} & \text{for all } 0\le i<j \le n, \\
\id_{n} & \text{for all } 0 \le i=j \le n \text{ or } 1\le i=j+1 \le n+1,\\
\delta_{i-1}^n\sigma_j^{n-1} & \text{for all } 1\le j+1 < i \le n+1.
\end{cases}
\end{gather}
\endgroup
In the op\-po\-site ca\-te\-go\-ry~$\Delta^{\opp}$, every co\-face~$\delta_i^n$ and every co\-de\-ge\-ne\-ra\-cy~$\sigma_j^n$ are re\-spec\-tively de\-noted by
$
d_i^n \co n \to n-1$ and $s_j^n \co n \to n+1$.
The morphisms $\{d_i^n\}_{n\in \N^*, 0\le i \le n}$ are called \textit{faces} and the morphisms \smash{$\bigl\{s_j^n\bigr\}_{n\in \N, 0 \le j \le n}$} are called \textit{degeneracies}.

\subsection{The cyclic category}
The \emph{cyclic category} $\Delta C$ is introduced by Connes in \cite{Connesext}. We will use a more combinatorial definition from \cite[Section~6.1]{loday98}, which is as follows. The objects of~$\Delta C$ are the non-negative integers.
The morphisms in this category are generated by morphisms~$\{\delta_i^n\}_{n\in \N^*, 0\le i \le n}$, called \textit{cofaces}, morphisms~\smash{$\bigl\{\sigma_j^n\bigr\}_{n\in \N, 0\le j \le n}$}, called \textit{codegeneracies}, and isomorphisms $\{\tau_n\co n \to n\}_{n\in \N}$, called~\emph{cocyclic operators}, which satisfy the simplicial relations and additionally:
\begin{alignat}{3}
\label{compcoccof}
&\tau_n \delta_i^n = \delta_{i-1}^n \tau_{n-1} \qquad&& \text{for all } 1 \leq i \leq n,& \\
\label{tndelta0}
&\tau_n\delta_0^n=\delta_n^n \qquad &&\text{for all } n\ge 1,& \\
\label{compcoccod}
&\tau_n\sigma_i^n=\sigma_{i-1}^n\tau_{n+1} \qquad &&\text{for all } 1 \leq i \leq n,&\\
\label{tnsigma0}
&\tau_n\sigma_0^n=\sigma_n^n\tau_{n+1}^2 \qquad&& \text{for all } n\ge 0,&\\
\label{cocyclicity}
&\tau_n ^{n+1} = \id_{n} \qquad&& \text{for all } n \in \N.&
\end{alignat}
Note that $\tau_0=\id_{0}$.
In the op\-po\-site ca\-te\-go\-ry~$\Delta C^{\opp}$, every co\-face~$\delta_i^n$, every co\-de\-ge\-ne\-ra\-cy~$\sigma_j^n$, and every cocyclic operator $\tau_n$ are re\-spec\-tively de\-noted by $d_i^n \co n \to n-1$, $s_j^n \co n \to n+1$, and $t_n \co n \to n$.
The morphisms $\{d_i^n\}_{n\in \N^*, 0\le i \le n}$ are called \textit{faces}, the morphisms \smash{$\big\{s_j^n\big\}_{n\in \N, 0 \le j \le n}$} are called \textit{degeneracies}, and the morphisms $\{t_n\}_{n\in \N}$ are called \textit{cyclic operators}.

\subsection{(Co)simplicial and (co)cyclic objects in a category} \label{simpcyc}
Let~$\mc C$ be any category.
A~\textit{simplicial object} in~$\mc C$ is a functor~$\Delta^\opp \to \mc C$ and a~\emph{cyclic object} in~$\mc C$ is a functor~$\Delta C^{\opp} \to \mc C$.
Dually, a~\emph{cosimplicial object} in~$\mc{C}$ is a functor~$\Delta \to \mc C$ and a~\emph{cocyclic object} in~$\mc C$ is a functor~$\Delta C \to \mc C$. A (co)simplicial/(co)cyclic object in the category of~$\kk$-modules is called a \textit{(co)simplicial/(co)cyclic~$\kk$-module}.
A \textit{morphism} between two (co)simplicial/(co)cyclic objects is a natu\-ral transformation between them.
For shortness, one often denotes the image of a morphism~$f$ under a (co)simplicial/(co)cyclic object in $\mc C$ by the same letter $f$.

Since the categories~$\Delta$ and~$\Delta C$ are defined by generators and relations, a (co)simpli\-cial/(co)\-cy\-clic object in a category is entirely determined by the images of the generators satisfying the corresponding relations. For example, a cocyclic object~$X$ in~$\mc C$ may be explicitly described as a~family~\smash{$X^\bullet = \{X^n\}_{n\in \N}$} of objects in~$\mc C$, equipped with morphisms~\smash{$\bigl\{\delta_i^n\co X^{n-1} \to X^{n}\bigr\}{}_{n\in \N^*,0\leq i \leq n}$}, called~\textit{cofaces}, morphisms~\smash{$\bigl\{\sigma_j^n\co X^{n+1} \to X^n\bigr\}{}_{n\in \N, 0 \leq j \leq n}$}, called~\textit{codegeneracies}, and isomorphisms \smash{$\{\tau_n\co X^n \to X^n\}_{n\in \N}$}, called~\textit{cocyclic operators}, which satisfy~\eqref{cofaces}--\eqref{cocyclicity}. Note that~$\tau_0$ is the identity.
A morphism~$\alpha^\bullet \co X^\bullet \to Y^\bullet$ between cocyclic objects~$X^\bullet$ and~$Y^\bullet$ in~$\mc C$ is then described by a fami\-ly~$\alpha^\bullet=\{\alpha^n\co X^n \to Y^n\}_{n\in \N}$ of morphisms in~$\mc C$ such that
\begin{alignat}{3}
\label{commcofac}&\delta_i^n\alpha^{n-1}= \alpha^n\delta_i^n \qquad&& \text{for all } n\in \N^* \text{ and } 0\leq i \leq n,& \\
\label{commcodegen} &\sigma_j^n \alpha^{n+1}= \alpha^n\sigma_j^n \qquad&& \text{for all } n\in \N \text{ and } 0\leq j \leq n,& \\
\label{eq-commute-cocyc}
&\alpha^n\tau_n=\tau_n\alpha^n \qquad &&\text{for all } n \in \N.&
\end{alignat}
Similarly, a cyclic object~$X$ in~$\mc C$ may be seen as a family~$X_\bullet=\{X_n\}_{n\in \N}$ of objects in~$\mc C$ equipped with morphisms~\smash{$\{d_i^n\co X_n \to X_{n-1}\}_{n\in \N^*,0\leq i \leq n}$}, called~\textit{faces}, morphisms {$\bigl\{s_j^n\co X_n \to \smash{X_{n+1}\bigr\}{}_{n\in \N, 0\leq j \leq n}}$}, called~\textit{degeneracies}, and isomorphisms~$\{t_n\co X_n \to X_n\}_{n\in \N}$, called~\textit{cyclic opera\-tors}, which satisfy the corresponding relations in~$\Delta C^\opp$.
Also, a morphism~$\alpha_\bullet \co X_\bullet \to Y_\bullet$ between two cyclic objects~$X_\bullet$ and~$Y_\bullet$ in~$\mc C$ is described by a fa\-mi\-ly~$\alpha_\bullet=\{\alpha_n\co X_n \to Y_n\}_{n\in \N}$ of morphisms in~$\mc C$ commuting with faces, degeneracies and cyclic operators of~$X_\bullet$ and~$Y_\bullet$.

\subsection{Cyclic duality and reindexing involution automorphism} \label{connes loday dual}
It is well known that the cyclic category is isomorphic to its opposite category. The isomorphism established by Connes in~\cite{Connesext} is called \textit{cyclic duality}.
In its version due to Loday \cite[Proposition~6.1.11]{loday98}, the cyclic duality~$L \co \Delta C^\opp \to \Delta C$ is the identity on objects and it is defined on morphisms as follows.
For~$n\in \N^*$ and~$0\leq i \leq n$,
\[L(d_i^n)=
\begin{cases}
\sigma_i^{n-1} & \text{if } 0\leq i \leq n-1, \\
\sigma_0^{n-1}\tau_{n}^{-1} & \text{if } i=n,
\end{cases}\]
and for $n\in \N$ and $0\leq j \leq n$,
\[L(s_j^n)=\delta_{j+1}^{n+1} \qquad \text{and} \qquad L(t_n)=\tau_n^{-1}.\]

Given a category $\mc C$, the cyclic duality transforms any cocyclic object $X \co \Delta C \to {\mc C}$ in $\mc C$ into the cyclic object $XL\co \Delta C^\opp \to \mc C$. Similarly, the opposite functor~$L^\opp$ turns any cyclic object~$Y \co \Delta C^\opp \to \mc C$ in $\mc C$ into the cocyclic object~$YL^\opp \co \Delta C \to {\mc C}$.
Following Loday~\cite[Sec\-tion~6.1.14]{loday98}, we also recall the \textit{reindexing involution automorphism}~$\Phi$ of the cyclic category. It is identity on objects and it is defined on morphisms by formulas
\[\Phi(\delta_i^n)=\delta_{n-i}^n, \qquad \Phi(\sigma_j^n)=\sigma_{n-j}^n, \qquad \Phi(\tau_n)= \tau_n^{-1}.
\]

\section{3-cobordisms} \label{seccobs}
In this section we recall some facts about the category~$\textbf{3}\Cob_0$ of~$3$-dimensional cobordisms (or shortly, $3$-cobordisms) and their surgery presentations via ribbon graphs.
We denote by~$D^n$ the closed unit ball in~$\R^n$. The~$n$-dimensional sphere is denoted by~$S^n$.
All knots considered in this paper are smoothly embedded. For more details, we suggest some of the standard references on knot theory, such as \cite{burde2002knots, lickorish97}.

\subsection{3-cobordisms}
A~\emph{$3$-cobordism} is a quadruple~$(M,h,\Sigma_1,\Sigma_2)$, where~$M$ is a compact oriented~$3$-manifold,~$\Sigma_1$ and~$\Sigma_2$ are two closed oriented surfaces, and~$h$ is an orientation preserving homeomorphism $h\co (-\Sigma_1) \sqcup \Sigma_2 \to \partial M$. The surface~$\Sigma_1$ is called the~\textit{bottom base} and the surface~$\Sigma_2$ is called the~\textit{top base} of the cobordism~$M$.
Two cobordisms~$(M,h,\Sigma_1,\Sigma_2)$ and~$(M',h',\Sigma_1,\Sigma_2)$ are \textit{homeomorphic}, if there is an orientation preserving homeomorphism~$g\co M\to M'$ such that~$h'=g_{\vert\partial M}h$.
When clear, we will denote a cobordism~$(M,h,\Sigma_1,\Sigma_2)$ only by $M$.

The composition of two cobordisms $(M_1,h_1,\Sigma_1,\Sigma_2)$ and $(M_2,h_2,\Sigma_2,\Sigma_3)$ is the cobordism $(M,h,\Sigma_1,\Sigma_3)$, where $M$ is obtained by gluing $M_1$ to $M_2$ along $h_2h_1^{-1}\co h_1(\Sigma_2) \to h_2(\Sigma_2)$ and the homeomorphism~$h$ is given by
\[
h=h_1\vert_{\Sigma_1} \sqcup h_2\vert_{\Sigma_3} \co\ (-\Sigma_1) \sqcup \Sigma_3 \to \partial M.
\]
We say that cobordism~$M$ is obtained by gluing cobordisms~$M_1$ and~$M_2$ along~$\Sigma_2$.

\subsection{The category of 3-cobordisms}\label{sect-cob3}
The category~$\textbf{3}\Cob_0$ of~$3$-cobordisms is defined as follows.
The objects are closed oriented surfaces.
A morphism~$f\co \Sigma_1 \to \Sigma_2$ in~$\textbf{3}\Cob_0$ is a homeomorphism class of cobordisms between~$\Sigma_1$ and~$\Sigma_2$.
In $\textbf{3}\Cob_0$, the identity of a closed oriented surface~$\Sigma$ is re\-presented by \emph{identity cobordism}~$(C_\Sigma,e,\Sigma,\Sigma)$, where~$C_\Sigma=\Sigma\times [0,1]$ is a cylinder over~$\Sigma$ together with the product orientation, and~$e\co (-\Sigma)\sqcup \Sigma \to \partial C_{\Sigma}$ is the homeomorphism with~$e\vert_{-\Sigma}(x,0)=(x,0)$ and~$e\vert_{\Sigma}(x,1)=(x,1)$.
The composition of morphisms~$\Sigma_1\to \Sigma_2$ and~$\Sigma_2 \to \Sigma_3$ in~$\textbf{3}\Cob_0$, re\-presented respectively by cobordisms $M$ and $N$, is re\-presented by the cobordism obtained by gluing cobordisms~$M$ and~$N$ along $\Sigma_2$.
The category~$\textbf{3}\Cob_0$ is symmetric monoidal (see Section~\ref{braidedcats}). The monoidal product is given by disjoint
union and the monoidal unit is the empty surface. For more details, see \cite{moncatstft}.

\subsection{Surgery presentation of closed 3-manifolds} \label{dehn}
Let~$L$ be an~$n$-component framed link in the~$3$-sphere~$S^3$. Pick a closed tubular neighborhood~$N_L$ of $L$.
Since~$N_L$ is ho\-meo\-morphic to~\smash{$\bigsqcup_{i=1}^n S^1 \times D^2$}, the boundary of the~$3$-manifold~$S^3 \setminus \mathrm{Int}(N_L)$ is homeomorphic to the disjoint union of~$n$-tori~$S^1 \times S^1$.
The \textit{Dehn surgery on~$S^3$ along}~$L$ is the closed manifold
\[
S^3_L=\bigl(S^3\setminus N_L\bigr)\bigcup_{\phi} \biggr(\bigsqcup_{i=1}^n D^2\times S^1\biggr),
\]
where~\smash{$\phi\co \partial\bigl(S^3\setminus \mathrm{Int}(N_L)\bigr)\to \bigsqcup_{i=1}^n S^1\times S^1$} is a ho\-me\-o\-mor\-phism ex\-changing me\-ri\-dians and pa\-ral\-lels.
Any connected, oriented, closed~$3$-manifold~$M$ is, according to the Lickorish's theorem~\cite[Section~12]{lickorish97}, homeomorphic to~$S^3_L$ for some framed link~$L\subset S^3$.

\subsection{Ribbon graphs} \label{ribbongrph}
A~\emph{circle} is a~$1$-manifold homeomorphic to~$S^1$. An \emph{arc} is a~$1$-manifold homeomorphic to the closed interval~$[0,1]$.
The boundary points of an arc are called its \emph{endpoints}.
A \emph{rectangle} is a~$2$-manifold with corners homeomorphic to~$[0,1] \times [0,1]$. The four corner points of a rectangle split its boundary into four arcs called the~\emph{sides}. A~\emph{coupon} is an oriented rectangle with a~distinguished side called the~\emph{bottom base}, the opposite side being the~\emph{top base}.

A~\emph{plexus} is a topological space obtained from a disjoint
union of a finite number of oriented circles, oriented arcs, and coupons
by gluing some endpoints of the arcs to the bases of the coupons. It is required that different endpoints of the arcs are never glued to the same point of a~(base of a) coupon. The endpoints of the arcs that are not glued to coupons are called~\emph{free ends}. The set of free ends of a plexus~$\gamma$ is denoted by~$\partial \gamma$.

Given non-negative integers~$g$ and~$h$, a \textit{ribbon}~$(g,h)$-\textit{graph}~$\Gamma$ is a plexus~$\Gamma$ embedded in ${\R^2 \times [0,1]}$ and equipped with a framing such that
\[
\partial \Gamma=\Gamma\cap \partial \bigl(\R^2 \times [0,1]\bigr) =\{(1,0,0), \dots, (g,0,0)\} \cup \{(1,0,1), \dots, (h,0,1)\}
\]
and such that the arcs of~$\Gamma$ are transverse to~$\partial \bigl(\R^2 \times [0,1]\bigr)$ at all points of~$ \partial \Gamma$.
The free end~$(i,0,0)$ is called the~\emph{$i$-th input} and the free end~$(j,0,1)$ is called the~\emph{$j$-th output} of~$\Gamma$.
For example, ribbon graphs without free ends and without coupons are nothing but framed oriented links in~$\R^2 \times (0,1) \cong \R^3$.

We represent ribbon graphs by plane diagrams with blackboard framing. We require that for each coupon, its orientation is that of the plane, its bases are horizontal, and its bottom base is below its top base. By Reidemeister theorem (see \cite{reidemeister}), two diagrams represent isotopic ribbon graphs if and only if they are related by a finite sequence of plane isotopies, ribbon Reidemeister moves~\Reid{1}--\Reid{3}
	\[
 \put(29,-27){\small $\textbf{R}1$}
 \raisebox{-6mm}{\includegraphics[scale=0.85]{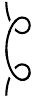}}\, \leftrightarrow\,
 \raisebox{-6mm}{\includegraphics[scale=0.85]{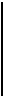}}\, \leftrightarrow\,
 \raisebox{-6mm}{\includegraphics[scale=0.85]{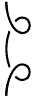}}\,,
 \put(57.5,-27){\small $\textbf{R}2$}
 \qquad
 \raisebox{-6mm}{\includegraphics[scale=0.85]{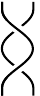}}\, \leftrightarrow\,
 \raisebox{-6mm}{\includegraphics[scale=0.85]{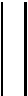}}\, \leftrightarrow\,
 \raisebox{-6mm}{\includegraphics[scale=0.85]{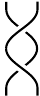}}\,,
 \put(67.5,-27){\small $\textbf{R}3$}
 \qquad
 \raisebox{-6mm}{\includegraphics[scale=0.85]{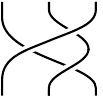}}\, \leftrightarrow\,
 \raisebox{-6mm}{\includegraphics[scale=0.85]{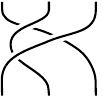}}\,,
 \]
and the following move:
\[
\raisebox{-9.3mm}{\includegraphics[scale=0.75]{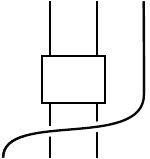}}
\put(-33.5,25){\small $\cdots$}
\put(-33.5,-25){\small $\cdots$}
\put(-31.5,-1){\small $C$}
\,\leftrightarrow\,
\raisebox{-9.4mm}{\includegraphics[scale=0.75]{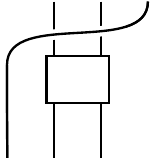}}\,.
\put(-36.5,25.2){\small $\cdots$}
\put(-36.5,-26){\small $\cdots$}
\put(-35.7,-1.5){\small $C$}
\]

\subsection{Standard ribbon graphs and surfaces}\label{cangrph}
For $g \geq 0$, consider the ribbon graph
\[
G_g^+ =\, \raisebox{-7.8mm}{\includegraphics[scale=0.7]{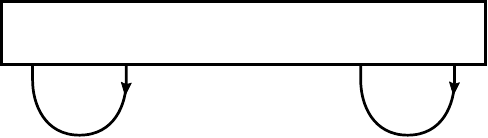}}
\put(-88,-5){\small $\cdots$}
\]
consisting of one coupon and $g$ unknotted untwisted cups oriented from right to left successively attached to the bottom base of the coupon. We fix a closed regular neighborhood~$H_g^+\subset \R^2\times[0,1]$ of~$G_g^+$. This is a handlebody of genus~$g$ and is provided with the right-handed orientation.
Consider also the ribbon graph
\[
G_g^- =\, \raisebox{-7.2mm}{\includegraphics[scale=0.7]{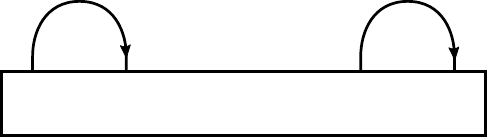}}
\put(-88,3){\small $\cdots$}
\]
consisting of one coupon and $g$ unknotted untwisted caps oriented from left to right successively attached to the top base of the coupon.
We fix a closed regular neighborhood~$H_g^-\subset \R^2\times[0,1]$ of~$G_g^-$. This is a handlebody of genus~$g$ and is provided with the right-handed orientation.
We choose the neighborhoods~$H_g^+$ and~$H_g^-$ so that the mirror reflection of~$\mathbb{R}^3$ with respect to the line~\smash{$\mathbb{R}^2\times\bigl\{\frac{1}{2}\bigr\}$} induces an orientation preserving homeomorphism between $-\bigl(H_g^-\bigr)$ and~$H_g^+$, where~$-\bigl(H_g^-\bigr)$ is $H_g^-$ with opposite orientation.

The boundary~$S_g$ of~$H_g^+$ is a closed connected oriented surface of genus~$g$ called the~\textit{standard surface of genus~$g$}. The orientation of~$S_g$ is induced by that of $H_g^+$. (We use the ``outward vector first'' convention for the induced orientation of the boundary.) The above mirror reflection induces a canonical orientation preserving homeomorphism between \smash{$\partial\bigl(H_g^-\bigr)$} and $-S_g$.

\subsection{Special ribbon graphs and 3-cobordisms} \label{specialribbon} Let~$g$ and~$h$ be non-negative integers.
A \emph{special ribbon~$(g,h)$-graph} is a ribbon~$(2g,2h)$-graph~$\Gamma$ with no coupons such that
\begin{itemize}\itemsep=0pt
\item for all $1 \le i \le g$, the $(2i-1)$-th and $2i$-th inputs of $\Gamma$ are connected by an arc oriented from the $(2i-1)$-th input to the $2i$-th input,
\item for all $1 \le j \le h$, the $(2j-1)$-th and $2j$-th outputs of $\Gamma$ are connected by an arc oriented from the $2j$-th output to the $(2j-1)$-th output.
\end{itemize}

Any special ribbon~$(g,h)$-graph $\Gamma$ gives rise to a connected $3$-cobordism $M_\Gamma$ between the standard surfaces $S_g$ and $S_{h}$ (see Section~\ref{cangrph}), which is defined as follows.
Attach coupons
\[
Q^-=[0,2g+1]\times \{0\}\times [-1,0] \qquad \text{and} \qquad Q^+=[0,2h+1]\times \{0\}\times [1,2]
\]
to~the bottom and the top of~$\Gamma$,~respectively. The result is a ribbon graph without free ends~$\tilde{\Gamma}$ in~$\R^3\cong \R^2\times (0,1)$, with two coupons~$Q^\pm$, and with finitely many circle components which form a framed link~$L$ in~$S^3$.
The arcs connecting the inputs of~$\Gamma$ become caps attached on the top base of~$Q^-$, so there is an embedding~$f^-\co G_g^- \to \tilde{\Gamma}$ mapping the coupon of~$G_g^-$ to $Q^-$ and mapping the caps attached to~$G_g^-$ to those attached to~$Q^-$.
Similarly, the arcs connecting the outputs of~$\Gamma$ become cups attached on the bottom base of~$Q^+$, so there is an embedding~$f^+\co G_{h}^+ \to \tilde{\Gamma}$ mapping the coupon of~$G_{h}^+$ to $Q^+$ and mapping the cups attached to~$G_{h}^+$ to those attached to~$Q^+$. Consider a tubular neighborhood~$N_L$ of~$L$ and embeddings~$\tilde{f}^-\colon H_g^-\to S^3 \setminus N_L$ and~$\tilde{f}^+\colon H_{h}^+\to S^3\setminus N_L$ respectively extending $f^-$ and~$f^+$.
Let $S^3_L$ be the Dehn surgery of~$S^3$ along~$L$ (see Section~\ref{dehn}).
The ma\-ni\-fold
\[
S^3_\Gamma=S^3_L \setminus \bigl(\tilde{f}^-\bigl(\mathrm{Int}\bigl(H_g^-\bigr)\bigr) \cup \tilde{f}^+\bigl(\mathrm{Int}\bigl(H_{h}^+\bigr)\bigr)\bigr)
\]
is a connected oriented compact $3$-manifold, whose boundary \smash{$\tilde{f}^-\bigl(\partial\bigl(H_g^-\bigr)\bigr)\sqcup \tilde{f}^+\bigl(\partial\bigl(H_{h}^+\bigr)\bigr)$}, following Section~\ref{cangrph}, is canonically homeomorphic to $(-S_g) \sqcup S_{h}$. This gives rise to a connected $3$-cobordism $M_\Gamma \co S_g \to S_{h}$.

For example, by \cite[Section~IV, Lemma 2.6]{turaevqinvariants}, the identity cobordism of the standard surface~$S_g$ (see Section~\ref{sect-cob3}) is represented by the following special ribbon $(g,g)$-graph:
\begin{equation} \label{eq-id-graph}
I_g=\, \raisebox{-15mm}{\includegraphics[scale=0.7]{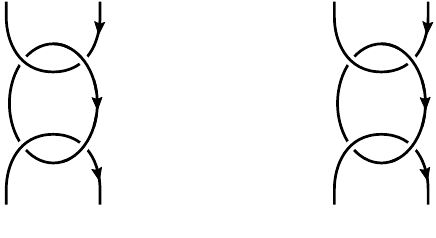}}\, .
\put(-82,0){\small $\cdots$}
\put(-136,-40){\small $1$}
\put(-25,-40){\small $g$}
\end{equation}

The following lemma gives a presentation of the composition of 3-cobordisms between standard surfaces.
\begin{Lemma}[{\cite[Section~\RomanNumeralCaps{4}.2.3]{turaevqinvariants}}]\label{naivecompo}
If $\Gamma$ is a special ribbon~$(g,h)$-graph and $\Gamma'$ is a special ribbon~$(h,k)$-graph, then the composition of $3$-cobordisms $M_\Gamma \co S_g \to S_h$ and $M_{\Gamma'} \co S_{h} \to S_{k}$ is the~$3$-cobordism~$M_{\Gamma'\circ \Gamma} \co S_{g} \to S_{k}$, where~$\Gamma'\circ \Gamma$ is the special ribbon~$(g,k)$-graph obtained by stacking~$\Gamma'$ over~$\Gamma$.
\end{Lemma}
We note that the connectedness of top base of $M_\Gamma$ and bottom base of $M_\Gamma'$ is here important. For the general case of $3$-cobordisms with non-connected bases, which we do not need in what will follow, see \cite[Section~IV.2.8]{turaevqinvariants}.

\subsection{Extended Kirby calculus}\label{kirbycalculus}
Recall from Section~\ref{dehn} that any closed oriented {$3$-manifold} can be obtained by surgery of~$S^3$ along a framed link.
Kirby proved~\cite{kirby} that two framed links represent the same~$3$-manifold (up to an orientation preserving homeomorphism) if and only if they are related by a finite sequence of isotopies and of the \textit{Kirby moves}~$\textbf{K}1$ and~$\textbf{K}2$.
The move $\textbf{K}1$ consists in adding an unknot with the framing number (which is the self-linking number)~$1$ or~$-1$,
\[
L \sqcup \rsdraw{0.4}{0.7}{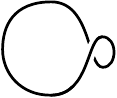} \longleftrightarrow L \longleftrightarrow L \sqcup \rsdraw{0.4}{0.7}{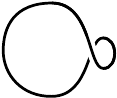}\hspace{0.1cm}.
\]
The move $\textbf{K}2$ consists in sliding a component over another component. More precisely, given two distinct components $L_i$ and $L_j$ of a framed link, this move replaces $L_i$ by the connected sum $L_i \# L_j'$ of $L_i$ with a copy $L'_j$ of $L_j$ obtained by slightly pushing $L_j$ along its framing. For example, sliding an unknot with framing number $0$ over an unknot with framing number $1$ can be depicted as
\[
\raisebox{-8.8mm}{\includegraphics[scale=0.7]{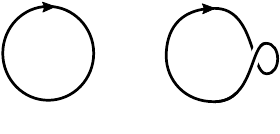}}
\put(-81.5,-23){\small $L_i$}
\put(-33,-1.5){\small $L_j$}
\,\longrightarrow\,
\raisebox{-8.8mm}{\includegraphics[scale=0.7]{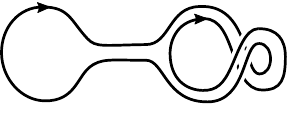}}\, .
\put(-89,-23){\small $L_i$}
\put(-40,-1.2){\small $L_j$}
\]
Building on \cite[Theorem~3 and Section~3.2.4]{kerler-bridged-links}, it follows that special ribbon graphs represent the same~$3$-cobordism (up to an orientation preserving homeomorphism) if and only if they are related by a finite sequence of isotopies and of the following moves: the move~$\textbf{K}1$, the generalized Kirby move~$\textbf{K}2'$, the move~\textbf{COUPON}, and the move~\textbf{TWIST}.
The move~$\textbf{K}2'$ consists in sliding an arc or circle component of a special ribbon graph over a distinct circle component.
The~\textbf{COUPON} move consists in changing the type of a crossing of a component passing over (or under) all its outputs
\[
\raisebox{-9mm}{\includegraphics[scale=0.75]{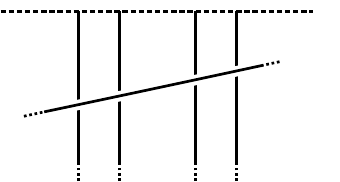}}
\put(-6,32){\small $\R^2\times \{1\}$}
\put(-71.5,25){\small $\cdots$}
\qquad\qquad \longleftrightarrow \quad
\raisebox{-9mm}{\includegraphics[scale=0.75]{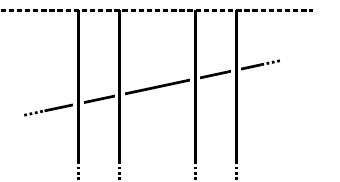}}
\put(-6,32){\small $\R^2\times \{1\}$}
\put(-71.5,25){\small $\cdots$}
\hphantom{\text{\small $\R^2\times \{1\}$}}\!,
\]
or all its inputs
\[
\raisebox{-9mm}{\includegraphics[scale=0.75]{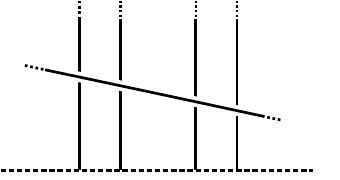}}
\put(-6,-25){\small $\R^2\times \{0\}$}
\put(-71.5,-18){\small $\cdots$}
\qquad\qquad \longleftrightarrow \quad
\raisebox{-9mm}{\includegraphics[scale=0.75]{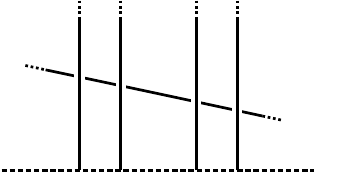}}
\put(-6,-25){\small $\R^2\times \{0\}$}
\put(-71.5,-18){\small $\cdots$}
\hphantom{\text{\small $\R^2\times \{0\}$}}\!,
\]
The move \textbf{TWIST} consists in a simultaneous twist of the output components
\[
\put(80, 58){\small $\R^2\times \{1\}$}
\put(33.5,45){\small $\cdots$}
\put(46.5,-52){\small $\cdots$}
\put(83,26){\small $\vdots$}
\put(83,-23){\small $\vdots$}
\raisebox{-20mm}{\includegraphics[scale=0.65]{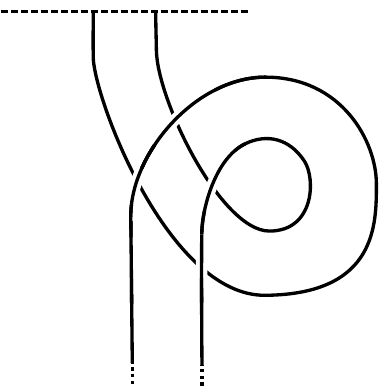}}\quad\longleftrightarrow
\put(81, 58){\small $\R^2\times \{1\}$}
\put(35,45){\small $\cdots$}
\,\raisebox{-20mm}{\includegraphics[scale=0.65]{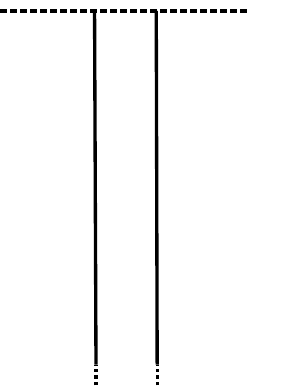}} \,\longleftrightarrow
\put(91, 58){\small $\R^2\times \{1\}$}
\put(44.5,45){\small $\cdots$}
\put(57.5,-52){\small $\cdots$}
\put(94,26){\small $\vdots$}
\put(94,-23){\small $\vdots$}
\quad\raisebox{-20mm}{\includegraphics[scale=0.65]{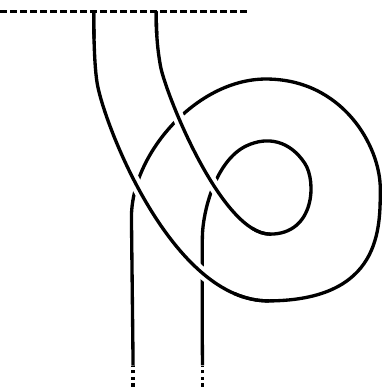}} \, ,
\]
or the input components
\[
\put(80, -57){\small $\R^2\times \{0\}$}
\put(34,-44){\small $\cdots$}
\put(47,53){\small $\cdots$}
\put(84,22){\small $\vdots$}
\put(84,-28){\small $\vdots$}
\raisebox{-20mm}{\includegraphics[scale=0.65]{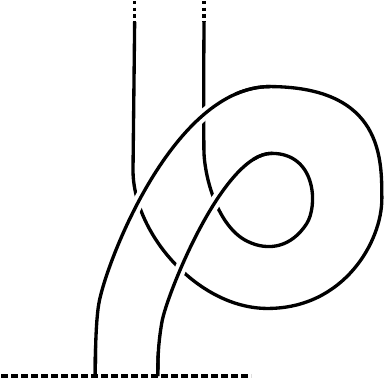}}\quad\longleftrightarrow
\put(81.5, -57){\small $\R^2\times \{0\}$}
\put(35.5,-44){\small $\cdots$}
\,\raisebox{-20mm}{\includegraphics[scale=0.65]{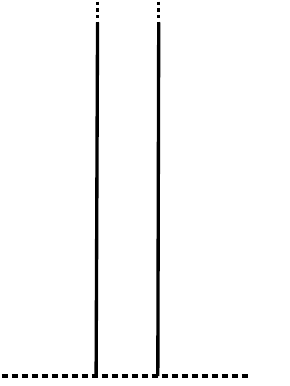}} \,\longleftrightarrow
\put(91, -57){\small $\R^2\times \{0\}$}
\put(44.5,-44){\small $\cdots$}
\put(58,53){\small $\cdots$}
\put(94,22){\small $\vdots$}
\put(94,-28){\small $\vdots$}
\quad\raisebox{-20mm}{\includegraphics[scale=0.65]{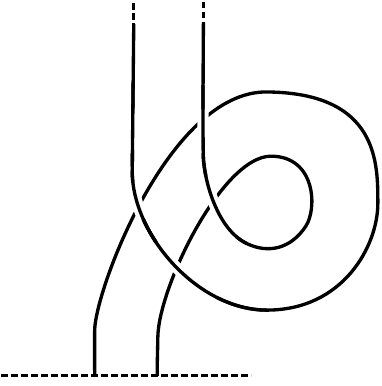}} \, .
\]
Note that in \cite[Theorem~3]{kerler-bridged-links} the move \textbf{COUPON} is called $\tau$-move. Also, the $\kappa$-move (see~\cite{Fenn-Rourke}) and the Hopf link move in \cite{kerler-bridged-links} can be replaced by the moves $\textbf{K}1$ and $\textbf{K}2'$ chosen here, while the~$\sigma$-move in \cite{kerler-bridged-links} can be omitted, as we work in the category of closed tangles in the sense of~\cite[Section~3.2.4]{kerler-bridged-links}.

\begin{Remark} A reviewer of the paper pointed out that the \textbf{TWIST} move is here redundant.
Indeed, the \textbf{TWIST} move follows by an application of \textbf{COUPON} move and the Fenn--Rourke move (see \cite{Fenn-Rourke}), which itself follows by application of $\textbf{K}1$, $\textbf{K}2'$, and isotopy. We keep it because it plays a prominent role in the proof of our main result.
\end{Remark}

\section{Cyclic objects from surfaces}\label{cobcyclic}
The first main result of this paper is that closed oriented surfaces can be organized in a (co)cyclic object in the category~$\textbf{3}\Cob_0$ of~$3$-dimensional cobordisms:

\begin{Theorem} \label{CYCINCOB}
For~$g\geq 1$, let $\Sigma_g$ be a closed oriented genus~$g$ surface. Then the fa\-mily~$\{\Sigma_g\}_{g\ge 1}$ has a structure of a cocyclic object~$X^\bullet$ in \emph{3\Cob}$_0$ and a structure of a cyclic object~$X_\bullet$ in \emph{3\Cob}$_0$.
\end{Theorem}
We construct (co)cyclic objects in~$\textbf{3}\Cob_0$ by means of surgery presentation of~$3$-cobordisms developed in~\cite{reshetikhin_invariants_1991, turaevqinvariants} and reviewed in Section~\ref{seccobs}.
We prove Theorem~\ref{CYCINCOB} in Sections~{\ref{construction}--\ref{generalsurf}}.
First, in Section~\ref{construction} we construct the functor~$Y^\bullet \co \Delta C\to \textbf{3}\Cob_0$.
Next, in Section~\ref{Stildesketch}, we construct the functor~$Y_{\bullet}\co \Delta C^{\opp}\to \textbf{3}\Cob_0$.
Cobordisms in both of these constructions have standard surfaces (see Section~\ref{cangrph}) as bases.
Finally, in Section~\ref{generalsurf} we pass from~$Y^\bullet$ and~$Y_\bullet$ to arbitrary~$X^\bullet$ and~$X_\bullet$, as stated in Theorem~\ref{CYCINCOB}.

A $3$-dimensional TQFT is a symmetric monoidal functor from $\textbf{3}\Cob_0$ to $\Mod_\kk$. By composition, we have the following.

\begin{Corollary} If $Z$ is a~$3$-dimensional~{TQFT}, then~$Z\circ X^\bullet$ is a cocyclic~$\kk$-module and $Z\circ X_\bullet$ is a cyclic~$\kk$-module.
\end{Corollary}
A fundamental construction of a~$3$-dimensional TQFT is the Reshetikhin--Turaev TQFT $\RT_{\mc B} \co \textbf{3}\Cob_0 \to \Mod_\kk$ associated to an anomaly free modular category~$\mc B$.
We postpone calculations of~$\RT_{\mc B}\circ X^\bullet$ and~$\RT_{\mc B}\circ X_\bullet$ to Section~\ref{tftsrel}.
Some algebraic preparations for it are given in Sections \ref{catprelimini} and \ref{algebraiccyclic}.

\subsection[The construction Y\^{}\{bullet\}]{The construction~$\boldsymbol{Y^\bullet}$}\label{construction}

Recall the standard surface~$S_g$ (see Section~\ref{cangrph}) of genus $g$. For~any~$n\in \N$, set~$Y^n=S_{n+1}$.
For~$n\in \N^*$~and~$0\le i\le n$, the faces~$Y^\bullet(\delta_i^n)\co S_{n} \to S_{n+1}$ are defined as follows.
The~morphism~$Y^\bullet(\delta_0^n)\co S_{n} \to S_{n+1}$~is~the~cobordism class presented by the special ribbon graph~$G_{Y^\bullet (\delta_0^n)}$:
\[
G_{Y^\bullet (\delta_0^n)}=\,
\put(39,-20){\small $\cdots$}
\put(9.5,-28){\small $1$}
\put(73,-28){\small $n$}
\raisebox{-9.9mm}{\includegraphics[scale=0.5]{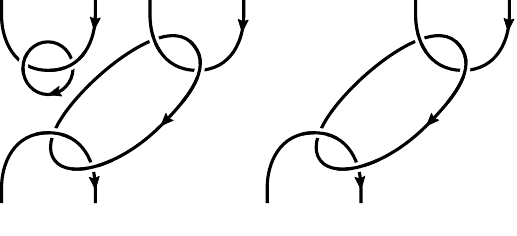}}\, .
\]
For~$1\le i \le n-1$, the morphism~$Y^\bullet(\delta_i^n)\co S_{n} \to S_{n+1}$ is defined as the cobordism class presented by the special ribbon graph~$G_{Y^\bullet (\delta_i^n)}$:
\[
G_{Y^\bullet (\delta_i^n)}=\,
\put(39,-20){\small $\cdots$}
\put(139,-20){\small $\cdots$}
\put(10,-28){\small $1$}
\put(74,-28){\small $i$}
\put(173,-28){\small $n$}
\raisebox{-9.9mm}{\includegraphics[scale=0.5]{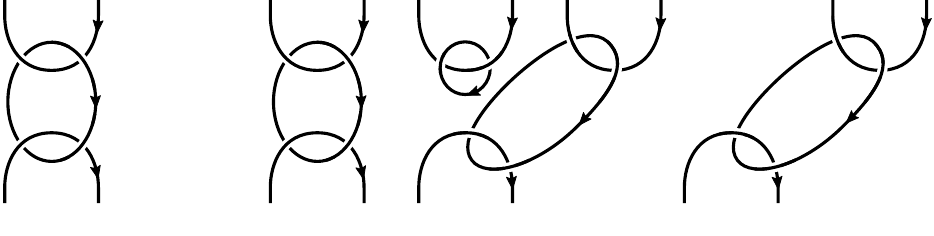}}\, .
\]
Finally, the morphism~$Y^\bullet(\delta_n^n)\co S_{n} \to S_{n+1}$ is defined as the cobordism class presented by the special ribbon graph~$G_{Y^\bullet (\delta_n^n)}$:
\[
G_{Y^\bullet (\delta_n^n)}=\,
\put(39,-20){\small $\cdots$}
\put(10,-28){\small $1$}
\put(73.5,-28){\small $n$}
\raisebox{-9.9mm}{\includegraphics[scale=0.5]{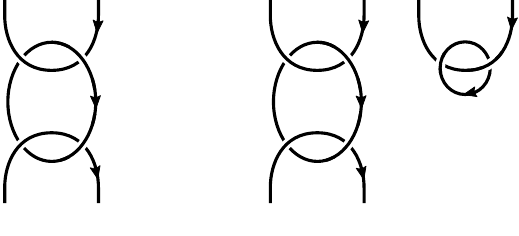}}\, .
\]
For $0\le j\le n$, the degeneracy $Y^\bullet\bigl(\sigma_j^n\bigr)\co S_{n+2} \to S_{n+1}$ is the cobordism class presented by the special ribbon graph $G_{Y^\bullet (\sigma^n_j)}$:
\[
G_{Y^\bullet (\sigma_j^n)}=\,
\put(39,-20){\small $\cdots$}
\put(210,-20){\small $\cdots$}
\put(10,-28){\small $0$}
\put(109,-28){\small $j$}
\put(137.5,-28){\small $j+1$}
\put(236.5,-28){\small $n+1$}
\raisebox{-9.9mm}{\includegraphics[scale=0.5]{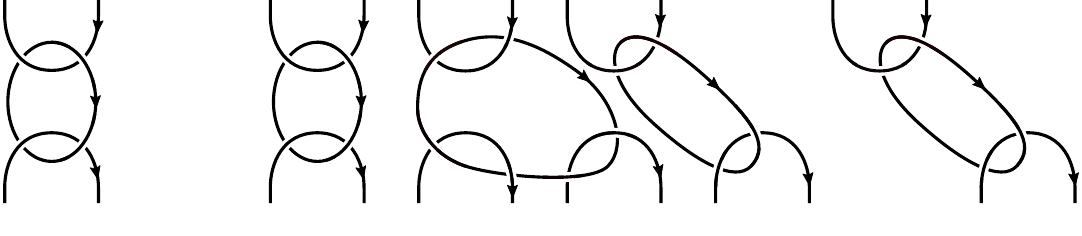}}\, .
\]
The morphism~$Y^\bullet(\tau_0)\co S_{1} \to S_{1}$ is the identity map~$\id_{S_1}$, which is represented by the special ribbon graph~$I_1$ depicted in~\eqref{eq-id-graph}.
For~$n\in \N^*$, the cocyclic operator~$Y^\bullet(\tau_n)\co S_{n+1} \to S_{n+1}$ is the cobordism class presented by the special ribbon graph~$G_{Y^\bullet (\tau_n)}$:
\[
G_{Y^\bullet (\tau_n)}=\,
\put(95,-32){\small $\cdots$}
\put(9,-43){\small $0$}
\put(45,-43){\small $1$}
\put(154,-43){\small $n$}
\raisebox{-15.2mm}{\includegraphics[scale=0.5]{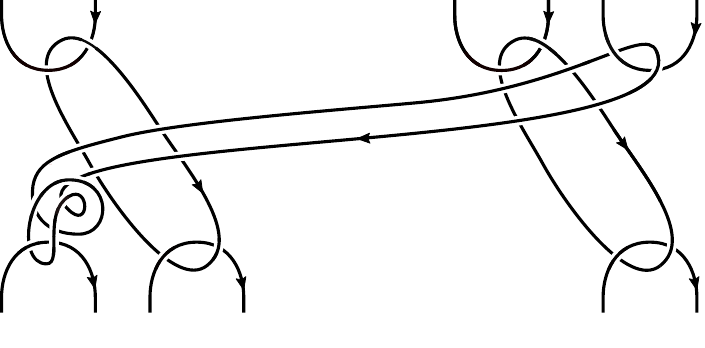}}\, .
\]

\begin{Lemma}
	\label{cobslem}
	The~fa\-mi\-ly~$Y^\bullet=\{S_{n+1}\}_{n\in \N}$,~equipped~with~the~co\-fa\-ces~$\{Y^{\bullet}(\delta_i^n)\}_{n\in \N^*, 0\le i \le n}$,~the co\-de\-ge\-ne\-ra\-cies~$\{Y^{\bullet}(\sigma_j^n) \}_{n\in \N, 0\le j \le n}$, and the co\-cy\-clic o\-pe\-ra\-tors~$\{Y^{\bullet}(\tau_n)\}_{n\in \N}$
	is a co\-cy\-clic object in $\textbf{\emph{3Cob}}_0$.
\end{Lemma}
To prove Lemma \ref{cobslem}, we intensively use some well-known consequences of the Kirby calculus from Section~\ref{kirbycalculus}, which we recollect in the following lemma.

\begin{Lemma} \label{kirbybash} One has the following moves on special ribbon graphs:
	\begingroup
	\allowdisplaybreaks
	\begin{align*}
	& (a) \quad \, \rsdraw{0.45}{0.5}{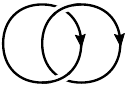}\; \lra \varnothing, \qquad
	(b) \quad \, \rsdraw{0.45}{0.5}{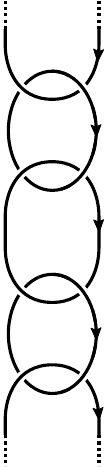} \;
	\lra \,\rsdraw{0.45}{0.5}{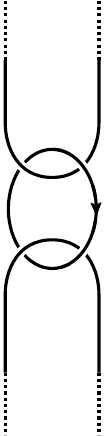}\;, \qquad
	(c) \quad\, \rsdraw{0.45}{0.5}{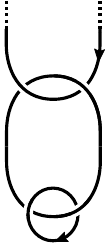} \;
	\lra \,\rsdraw{0.45}{0.5}{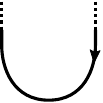}\;,& \\
\vspace{-10mm}\\
	& (d)\quad\, \rsdraw{0.45}{0.5}{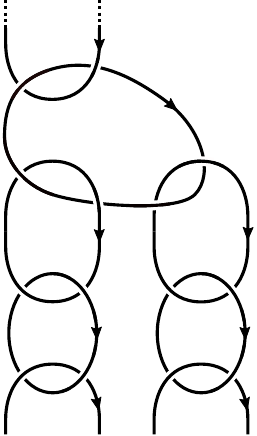} \;
	\lra \,\rsdraw{0.45}{0.5}{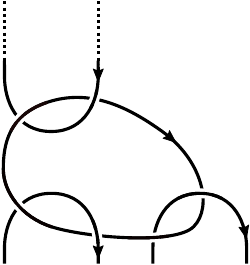}\;, \qquad
	(e)\quad\, \rsdraw{0.45}{0.5}{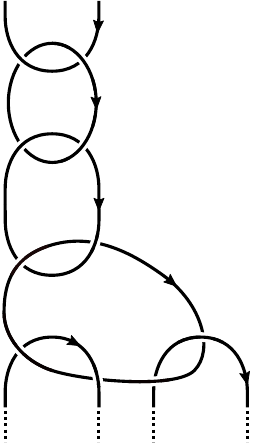} \;
	\lra \,\rsdraw{0.45}{0.5}{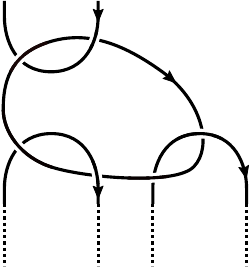}\;,& \\
\vspace{-10mm}\\
	& (f)\quad\,\rsdraw{0.45}{0.5}{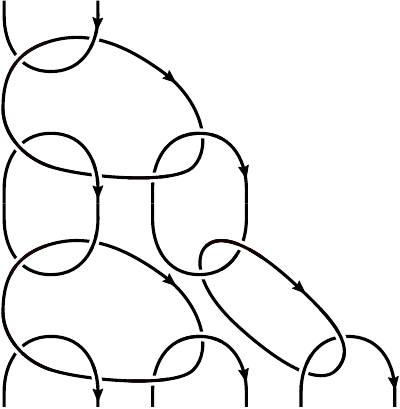}\; \lra \hspace{0.5cm} \,\rsdraw{0.45}{0.5}{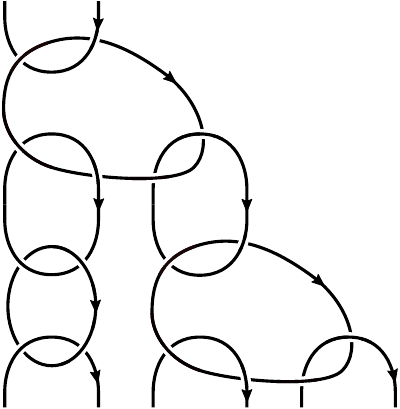}\;
	\lra \hspace{0.5cm} \,\rsdraw{0.45}{0.5}{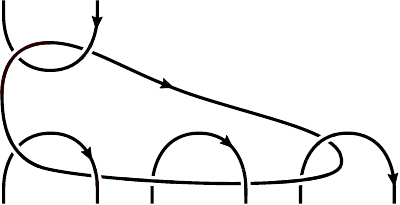}\;.&\\
\vspace{-10mm}
	\end{align*} \endgroup
\end{Lemma}
A proof of the part $(a)$ of the above lemma is given in \cite[Proposition~2]{kirby}. The other parts can be deduced from part $(a)$, isotopy and the generalized Kirby move $\textbf{K}2'$. One can also consult discussions in \cite[Section~3.1]{kerler2}.

\begin{proof}[Proof of Lemma \ref{cobslem}] In the proof, we shorten the notation by writing $\delta_i^n$, $\sigma_j^n$, and $\tau_n$ instead of $Y^\bullet(\delta_i^n)$, $Y^\bullet(\sigma_j^n)$, and $Y^\bullet(\tau_n)$.
	Let us verify the simplicial relation~\eqref{cofaces}.
	For $1\le i<j \le n$, we have
	\begingroup
	\allowdisplaybreaks
	\begin{align*}
	\delta_j^{n+1}\delta_i^n & \overset{(i)}{=}\,
 \put(27.1,48){\small $\cdots$}
 \put(140,48){\small $\cdots$}
 \put(253,48){\small $\cdots$}
 \put(27.1,-44.5){\small $\cdots$}
 \put(104,-44.5){\small $\cdots$}
 \put(182,-44.5){\small $\cdots$}
 \put(10.3,-53){\small $0$}
 \put(52.1,-53){\small $i$}
 \put(120,-53){\small $j-1$}
 \put(205.5,-53){\small $n$}
 \raisebox{-19mm}{\includegraphics[scale=0.5]{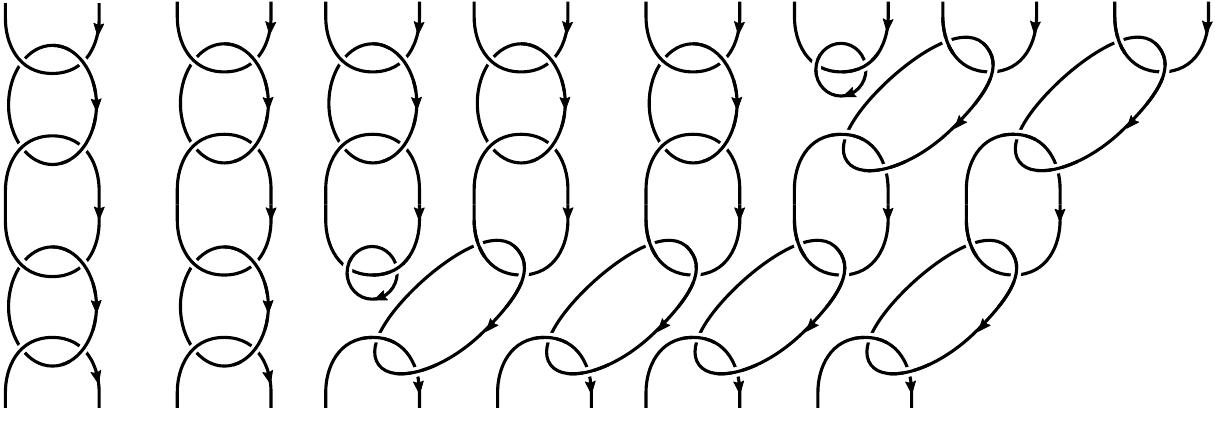}}\\
 &\overset{(ii)}{=} \,
 \put(27.5,70){\small $\cdots$}
 \put(140.3,70){\small $\cdots$}
 \put(253.6,70){\small $\cdots$}
 \put(27.5,-70){\small $\cdots$}
 \put(105,-70){\small $\cdots$}
 \put(182,-70){\small $\cdots$}
 \put(10.3,-78){\small $0$}
 \put(52.1,-78){\small $i$}
 \put(120,-78){\small $j-1$}
 \put(205.5,-78){\small $n$}
 \raisebox{-27.8mm}{\includegraphics[scale=0.5]{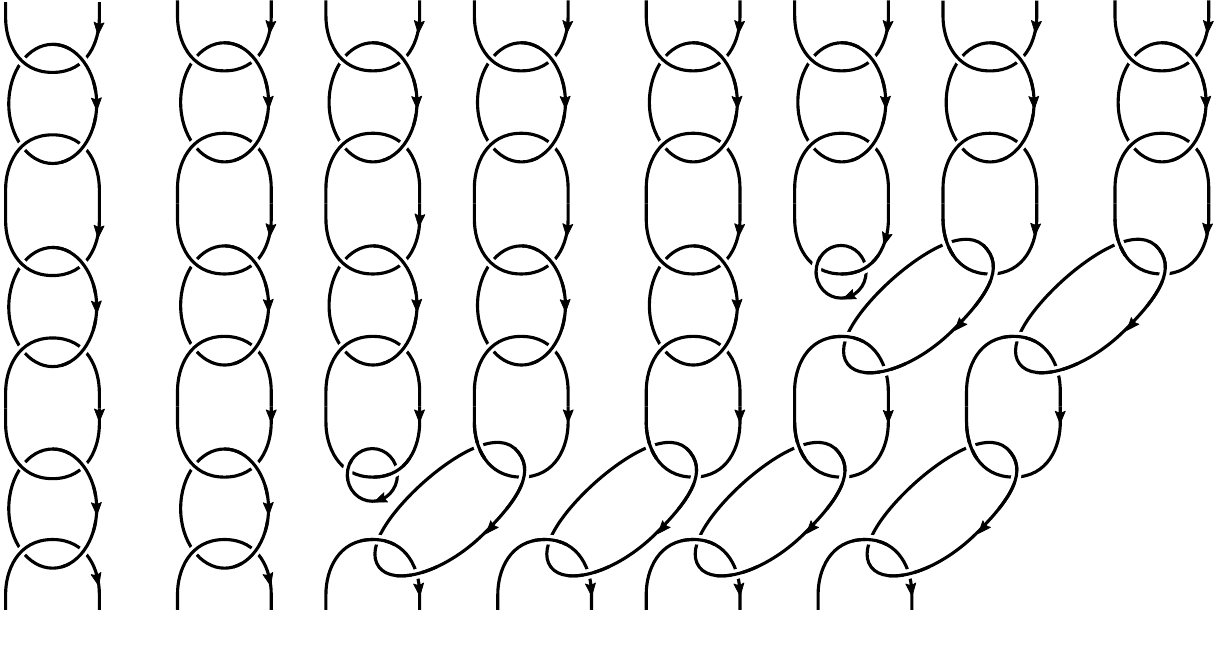}}\\
\vspace{-10mm}\\
 &\overset{(iii)}{=} \,
 \put(27.5,70){\small $\cdots$}
 \put(140.3,70){\small $\cdots$}
 \put(253.6,70){\small $\cdots$}
 \put(27.5,-70){\small $\cdots$}
 \put(105,-70){\small $\cdots$}
 \put(182,-70){\small $\cdots$}
 \put(10.3,-78){\small $0$}
 \put(52.1,-78){\small $i$}
 \put(120,-78){\small $j-1$}
 \put(205.5,-78){\small $n$}
 \raisebox{-28.5mm}{\includegraphics[scale=0.5]{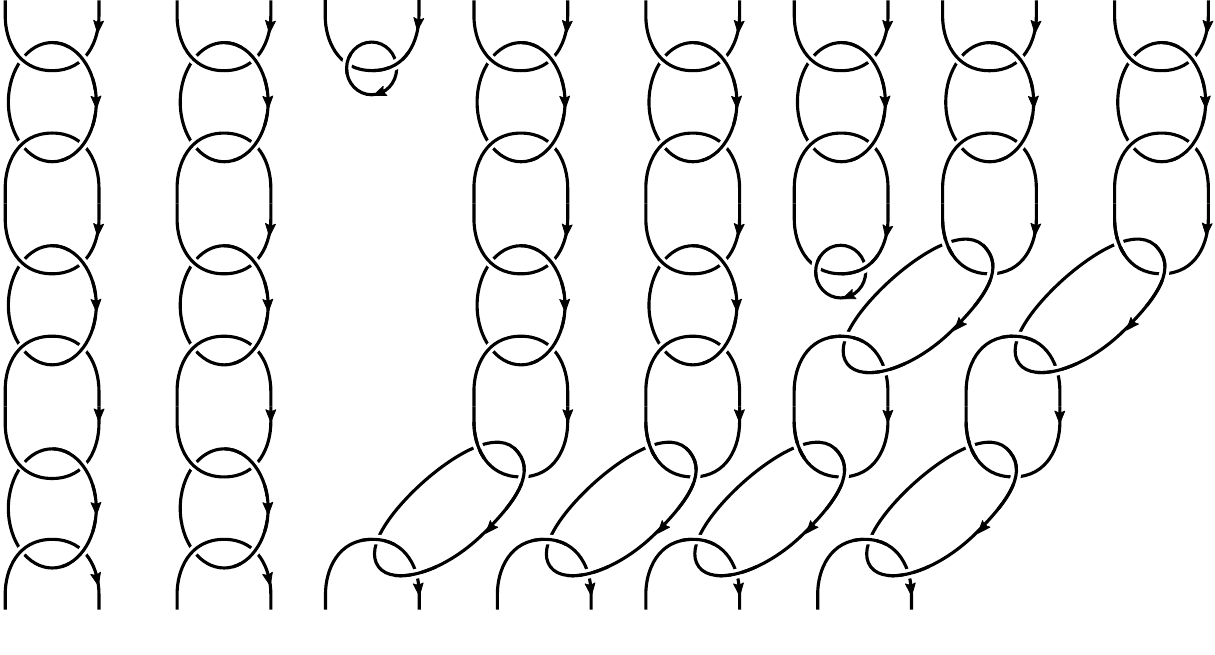}}\\
\vspace{-10mm}\\
 &\overset{(iv)}{=} \,
 \put(27.5,70){\small $\cdots$}
 \put(140.5,70){\small $\cdots$}
 \put(253.8,70){\small $\cdots$}
 \put(27.5,-70){\small $\cdots$}
 \put(105,-70){\small $\cdots$}
 \put(182,-70){\small $\cdots$}
 \put(10.3,-78){\small $0$}
 \put(52.1,-78){\small $i$}
 \put(120,-78){\small $j-1$}
 \put(205.5,-78){\small $n$}
 \raisebox{-28mm}{\includegraphics[scale=0.5]{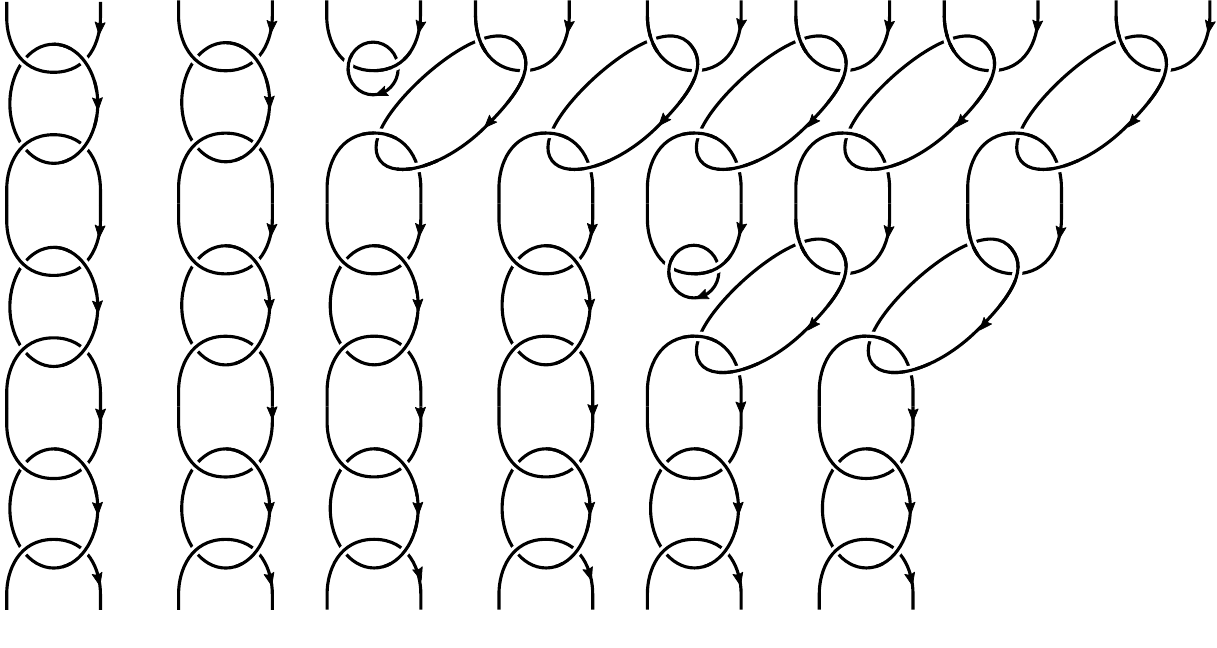}}\\
\vspace{-10mm}\\
 &\overset{(v)}{=} \,
 \put(27.2,48){\small $\cdots$}
 \put(140,48){\small $\cdots$}
 \put(253,48){\small $\cdots$}
 \put(27.2,-44.5){\small $\cdots$}
 \put(105,-44.5){\small $\cdots$}
 \put(182,-44.5){\small $\cdots$}
 \put(10.3,-54){\small $0$}
 \put(52.1,-54){\small $i$}
 \put(120,-54){\small $j-1$}
 \put(205.5,-54){\small $n$}
 \raisebox{-20mm}{\includegraphics[scale=0.5]{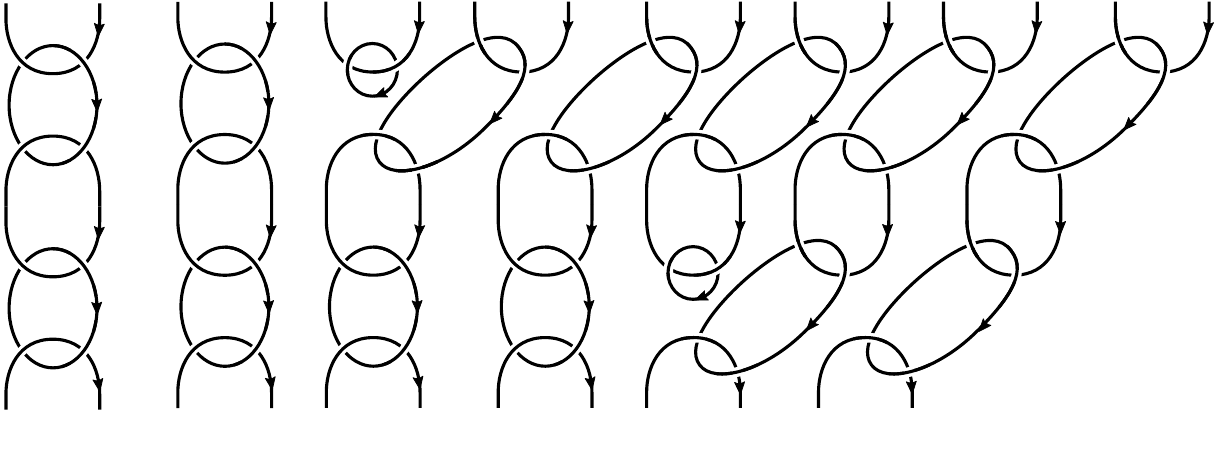}}
 \,\overset{(vi)}{=}\delta_i^{n+1}\delta_{j-1}^n.
	\end{align*}
	\endgroup
	Here~$(i)$ and~$(vi)$ follow from definitions and Lemma~\ref{naivecompo},~$(ii)$ and~$(v)$ follow from Lemma~\ref{naivecompo} and the fact that the graph from equation~\eqref{eq-id-graph} represents the identity cobordism,~$(iii)$ from Lemma~\ref{kirbybash}\,$(c)$, and~$(iv)$ by isotopy. The remaining cases are verified in a similar way.

	The relation~\eqref{codegeneracies} essentially follows from Lemma \ref{kirbybash}, parts $(d)$, $(e)$ and $(f)$.
We now verify the relation \eqref{compcofcod}. Suppose that $1\le i=j \le n-1$. We have
	\begingroup
	\allowdisplaybreaks
	\begin{align*}
 \sigma_i^n\delta_i^{n+1} &\overset{(i)}{=}\,
 \put(27.2,48){\small $\cdots$}
 \put(140.5,48){\small $\cdots$}
 \put(176.5,0){\small $\cdots$}
 \put(27.2,-48){\small $\cdots$}
 \put(140.5,-48){\small $\cdots$}
 \put(10.3,-55){\small $0$}
 \put(52.1,-55){\small $i$}
 \put(155.5,-55){\small $n+1$}
 \raisebox{-19.5mm}{\includegraphics[scale=0.5]{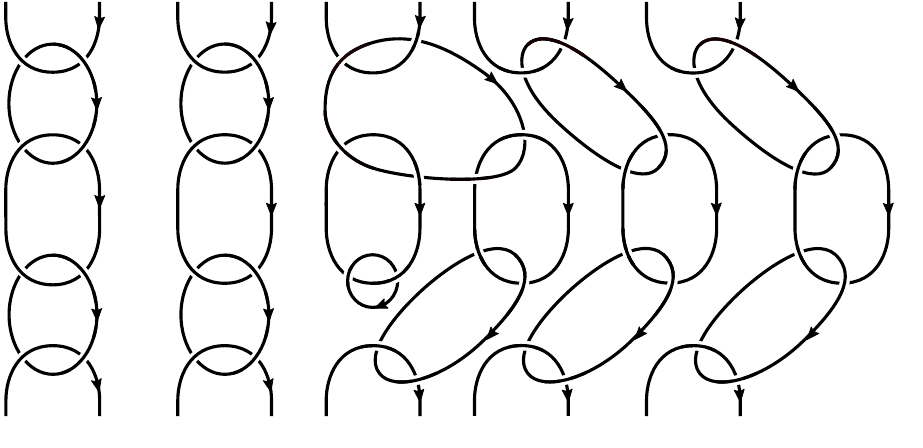}}	\\
\vspace{-10mm}\\
 &\overset{(ii)}{=}\,
 \put(27.2,48){\small $\cdots$}
 \put(140.5,48){\small $\cdots$}
 \put(176.5,0){\small $\cdots$}
 \put(27.2,-48){\small $\cdots$}
 \put(140.5,-48){\small $\cdots$}
 \put(10.3,-55){\small $0$}
 \put(52.1,-55){\small $i$}
 \put(155.5,-55){\small $n+1$}
 \raisebox{-19.5mm}{\includegraphics[scale=0.5]{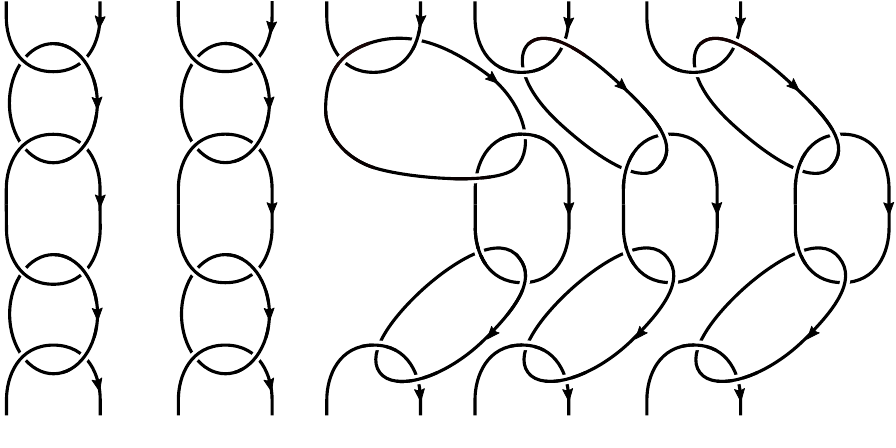}}	
 \,\overset{(iii)}{=}\id_{S_{n+1}}.
	\end{align*}\endgroup
	Here~$(i)$ follows from definition and Lemma~\ref{naivecompo},~$(ii)$ from Lemma~\ref{kirbybash}\,$(c)$, and~$(iii)$ by isotopy, Lemma~\ref{naivecompo}, and presentation of identity cobordism which is given in equation~\eqref{eq-id-graph}. The cases when~$i=0$ or $i=n$ are proven similarly.
The case when~$1\le i=j+1 \le n+1$ is verified in a~similar way as the case~$0 \le i=j \le n$. The cases when~$0 \le i<j \le n$ or~$1\le j+1<i \le n+1$ essentially follow from Lemma~\ref{kirbybash}\,$(c)$.
	
	It remains to show that relations~\eqref{compcoccof},~\eqref{compcoccod} and~\eqref{cocyclicity} hold.
	Indeed, according to~\cite[Section~6.1.1]{loday98}, these relations imply relations~\eqref{tndelta0} and~\eqref{tnsigma0}.
Let us verify the relation~\eqref{compcoccof}. In~the case when~$n\ge 3$ and~$2\le i \le n-1$, we have 	
\begingroup
	\allowdisplaybreaks
	\begin{align*}
	\tau_n\delta_i^n &\overset{(i)}{=}\,
 \put(28,71){\small $\cdots$}
 \put(141,71){\small $\cdots$}
 \put(176.5,0){\small $\cdots$}
 \put(63.5,0){\small $\cdots$}
 \put(63.5,-48){\small $\cdots$}
 \put(140.5,-48){\small $\cdots$}
 \put(10.3,-55){\small $1$}
 \put(46.2,-55){\small $2$}
 \put(88.5,-55){\small $i$}
 \put(164.5,-55){\small $n$}
 \raisebox{-19.5mm}{\includegraphics[scale=0.5]{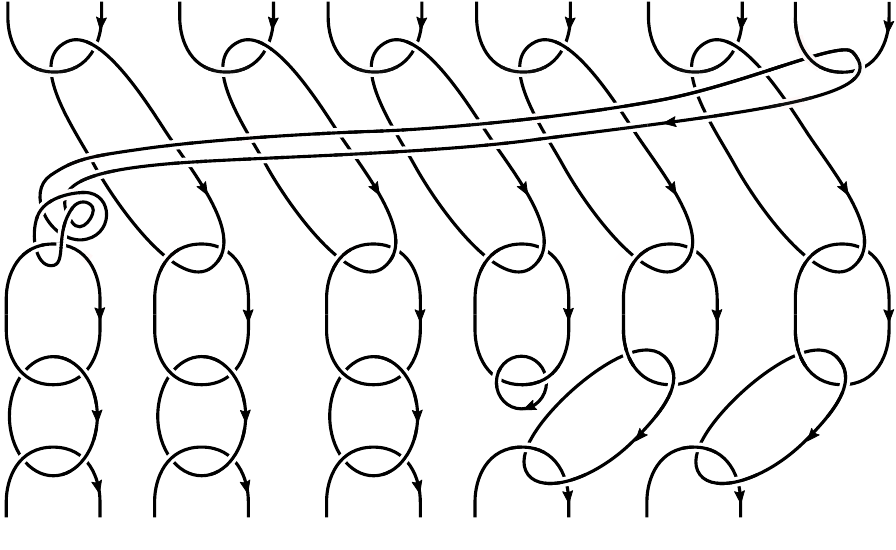}}	\\
\vspace{-10mm}\\
	&\overset{(ii)}{=}\,
 \put(27,71){\small $\cdots$}
 \put(140,71){\small $\cdots$}
 \put(175.5,0){\small $\cdots$}
 \put(62.5,0){\small $\cdots$}
 \put(62.5,-48){\small $\cdots$}
 \put(139.5,-48){\small $\cdots$}
 \put(9.3,-55){\small $1$}
 \put(45.2,-55){\small $2$}
 \put(87.5,-55){\small $i$}
 \put(163.5,-55){\small $n$}
 \raisebox{-19.5mm}{\includegraphics[scale=0.5]{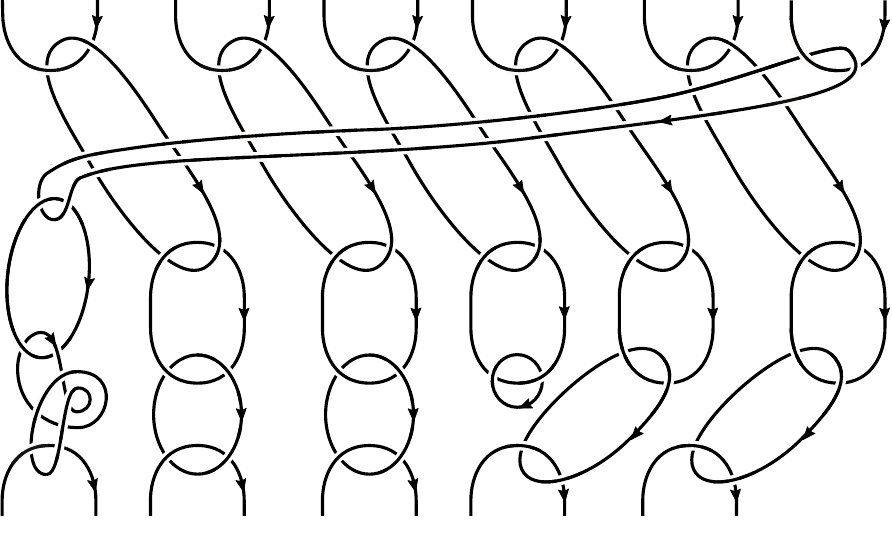}} \\
	&\overset{(iii)}{=}\,
 \put(27,58){\small $\cdots$}
 \put(140,58){\small $\cdots$}
 \put(27,10){\small $\cdots$}
 \put(104.5,10){\small $\cdots$}
 \put(62.5,-61){\small $\cdots$}
 \put(139.8,-61){\small $\cdots$}
 \put(9.3,-69){\small $1$}
 \put(45.2,-69){\small $2$}
 \put(87.5,-69){\small $i$}
 \put(163.5,-69){\small $n$}
 \raisebox{-25mm}{\includegraphics[scale=0.5]{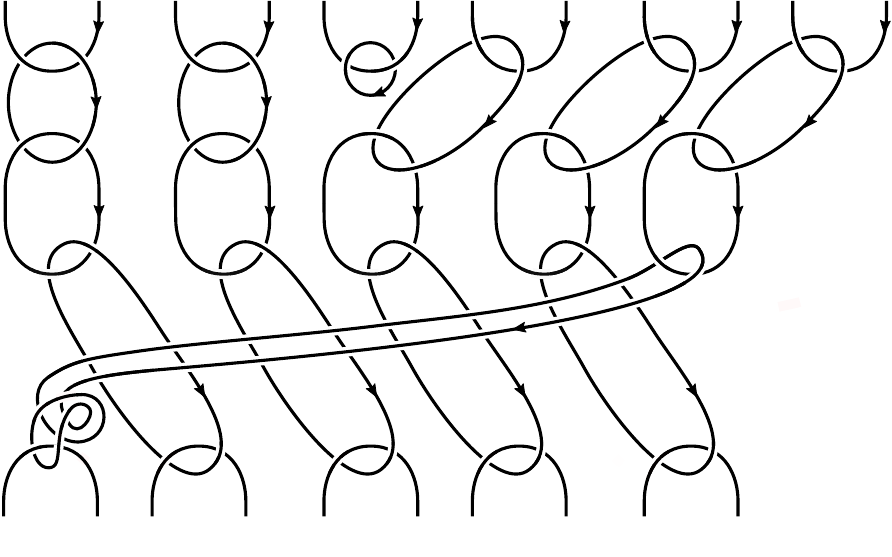}}
 \,\overset{(iv)}{=} \delta_{i-1}^n\tau_{n-1}.
	\end{align*}\endgroup
	Here~$(i)$ and~$(iv)$ follow from definitions and Lemma~\ref{naivecompo},~$(ii)$ by isotopy,~$(iii)$ follows from Lemma~\ref{kirbybash}\,$(c)$ and by isotopy. The remaining cases are proven similarly.
	
The relation~\eqref{compcoccod} essentially follows by applying the parts~$(b)$ and~$(e)$ of Lemma~\ref{kirbybash}.
	Finally, we check the relation~\eqref{cocyclicity} in the case~$n=1$.
	We have
\[	
	\tau_1^2 \overset{(i)}{=}\,	
 \put(9.3,-125){\small $0$}
 \put(45.3,-125){\small $1$}
 \raisebox{-44mm}{\includegraphics[scale=0.5]{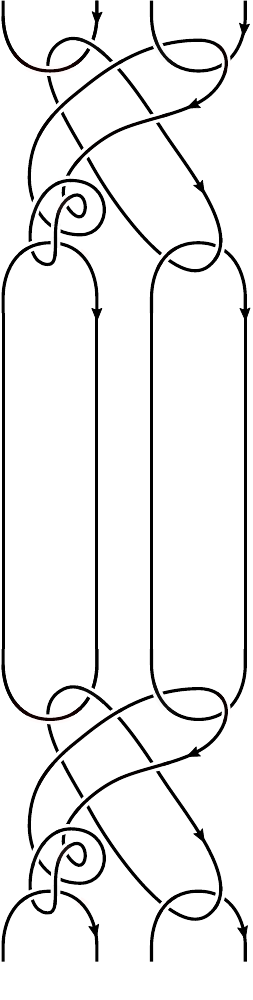}}
 \,\overset{(ii)}{=}\,	
 \put(9.3,-125){\small $0$}
 \put(45.3,-125){\small $1$}
 \raisebox{-44mm}{\includegraphics[scale=0.5]{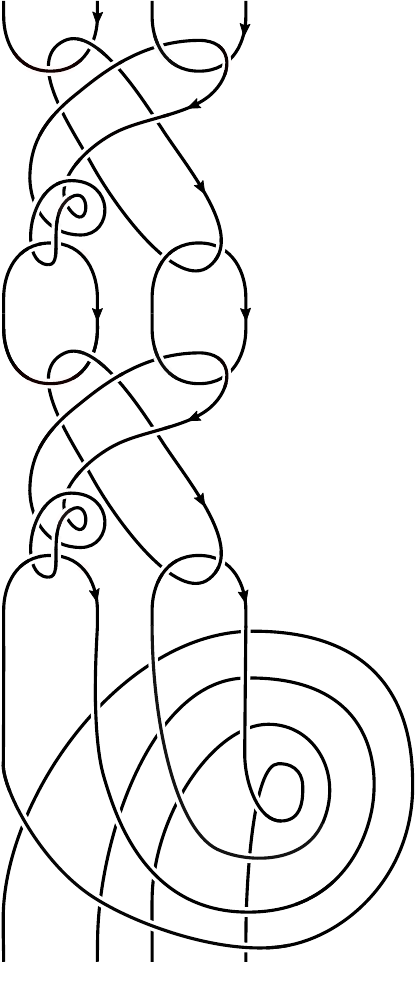}}
 \,\overset{(iii)}{=}\,
 \put(9.3,-125){\small $0$}
 \put(46.3,-125){\small $1$}
 \raisebox{-44mm}{\includegraphics[scale=0.5]{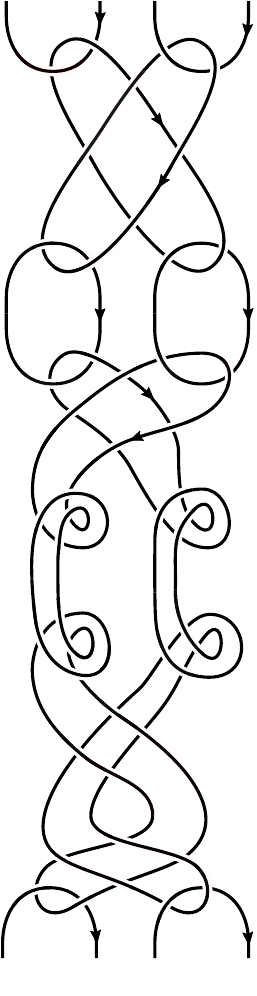}}
	\,\overset{(iv)}{=} \id_{S_2}.
\]
	Here~$(i)$ follows from definition and Lemma \ref{naivecompo},~$(ii)$ follows by the (negative)~\textbf{TWIST} move,~$(iii)$ follows by isotopy,~$(iv)$ follows by isotopy, Lemma~\ref{naivecompo}, and presentation of the identity cobordism, as depicted in~\eqref{eq-id-graph}. The general case is proven by a~similar reasoning. This finishes the proof of Lemma \ref{cobslem}. \end{proof}

\subsection[The construction Y\_\{bullet\}]{The construction~$\boldsymbol{Y_\bullet}$} \label{Stildesketch}
Recall the standard surface $S_g$ of genus $g$ from Section~\ref{cangrph}.
For~$n\in \N$, denote~$Y_n=S_{n+1}$.
For~$n\in \N^*$ and~$0\le i\le n$, the face~$Y_\bullet(d_i^n)\co S_{n+1} \to S_{n}$ is the cobordism class presented by the special ribbon graph~$G_{Y_\bullet (d_i^n)}$:
\[
G_{Y_\bullet (d_i^n)}=\,
\put(39,-20){\small $\cdots$}
\put(174.5,-20){\small $\cdots$}
\put(10,-29){\small $0$}
\put(110,-29){\small $i$}
\put(209,-29){\small $n$}
\raisebox{-11mm}{\includegraphics[scale=0.5]{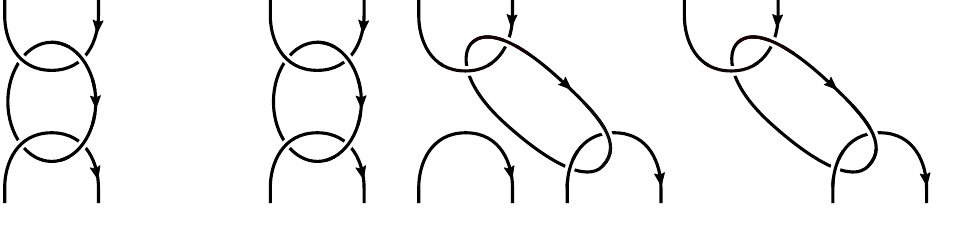}}\, .
\]
For~$0\le j\le n$, the morphism~$Y_\bullet\bigl(s_j^n\bigr)\co S_{n+1} \to S_{n+2}$ is the cobordism class presented by the special ribbon graph~$G_{Y_\bullet (s^n_j)}$:
\[
G_{Y_\bullet (s_j^n)}=\,
\put(39,-18){\small $\cdots$}
\put(177.5,-18){\small $\cdots$}
\put(10,-27){\small $0$}
\put(111,-27){\small $j$}
\put(213,-27){\small $n$}
\raisebox{-11mm}{\includegraphics[scale=0.5]{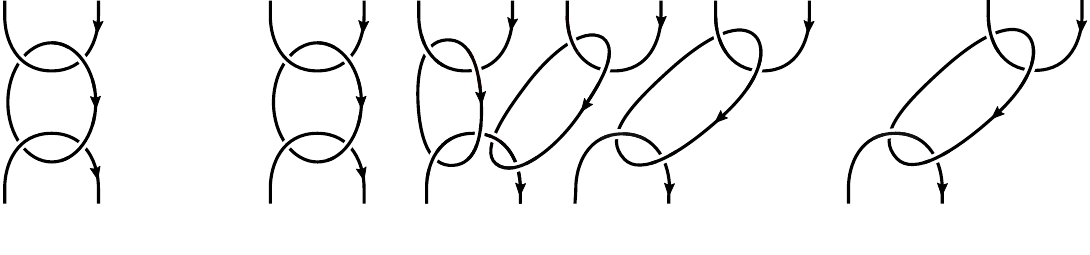}}\, .
\]
The morphism~$Y_\bullet(t_0)\co S_{1} \to S_{1}$ equals the identity map~$\id_{S_1}$, which is represented by the graph~$I_1$ depicted in~\eqref{eq-id-graph}.
For~$n\in \N^*$, the morphism~$Y_\bullet(t_n)\co S_{n+1} \to S_{n+1}$ is the cobordism class presented by the special ribbon graph~$G_{Y_\bullet (t_n)}$:
\[
G_{Y_\bullet (t_n)}=\,
\put(60,-30){\small $\cdots$}
\put(10,-42){\small $0$}
\put(110.5,-42){\small $n-1$}
\put(155.5,-42){\small $n$}
\raisebox{-15.5mm}{\includegraphics[scale=0.5]{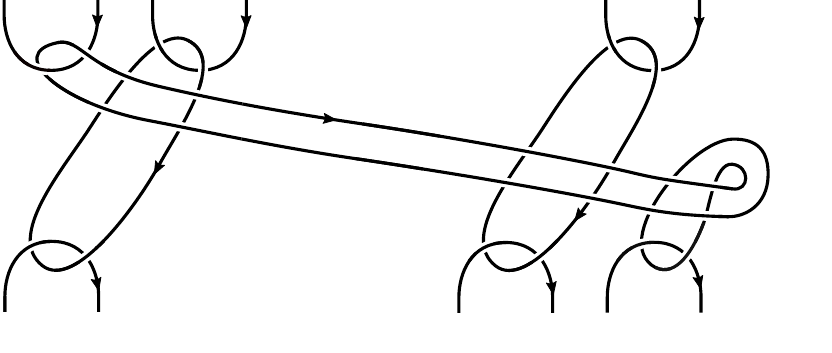}}\, .
\]
The proof of Theorem~\ref{CYCINCOB} in this case is similar to the proof of its version with~$X^\bullet$, which was detailed in Section~\ref{construction}. Namely, we have the following lemma.

\begin{Lemma}
	\label{cobslemv2}
	The~fa\-mi\-ly~$Y_\bullet=\{S_{n+1}\}_{n\in \N}$,~equipped~with~the~fa\-ces~$\{Y_\bullet(d_i^n)\}_{n\in \N^*, 0\le i \le n}$,~the de\-ge\-ne\-ra\-cies~$\bigl\{Y_\bullet\bigl(s_j^n\bigr) \bigr\}_{n\in \N, 0\le j \le n}$, and the cy\-clic o\-pe\-ra\-tors~$\{Y_\bullet(t_n)\}_{n\in \N}$
	is a cy\-clic object in~$\textbf{\emph{3Cob}}_0$.
\end{Lemma}
To prove Lemma~\ref{cobslemv2}, one uses the following result, which is analogous to Lemma~\ref{kirbybash}.

\begin{Lemma} \label{kirbybash2} One has the following moves on special ribbon graphs:
	\begingroup
	\allowdisplaybreaks
	\begin{align*}
	& (a) \quad\, \rsdraw{0.45}{0.5}{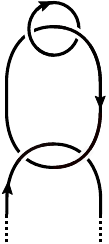} \;
	\lra \,\rsdraw{0.45}{0.5}{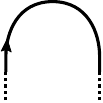}\;, \qquad (b) \quad \, \rsdraw{0.45}{0.5}{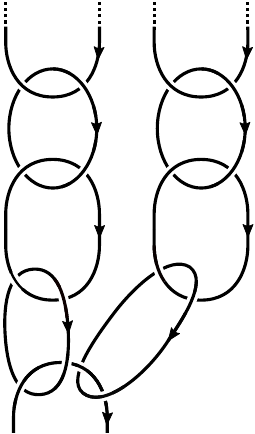} \;
	\lra \,\rsdraw{0.45}{0.5}{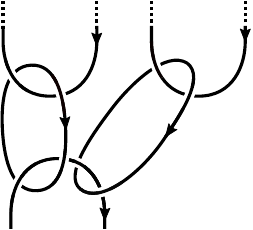}\;, &\\
	& (c) \quad\, \rsdraw{0.45}{0.5}{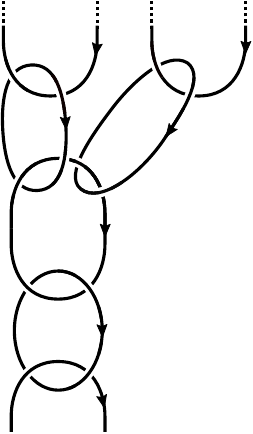} \;
	\lra \,\rsdraw{0.45}{0.5}{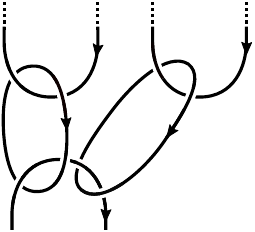}\;, &\\
	&(d)\quad\,\rsdraw{0.45}{0.5}{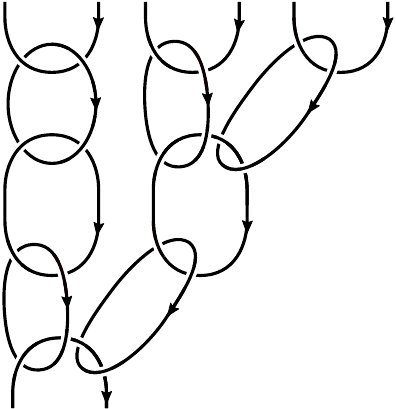}\; \lra \hspace{0.5cm} \,\rsdraw{0.45}{0.5}{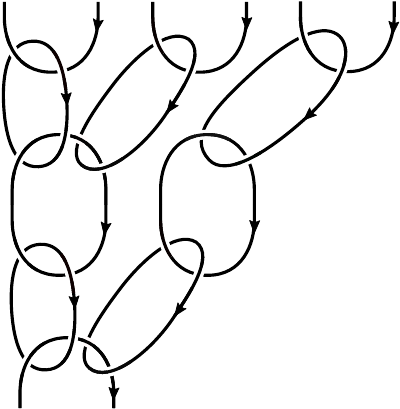}\;
	\lra \hspace{0.5cm} \,\rsdraw{0.45}{0.5}{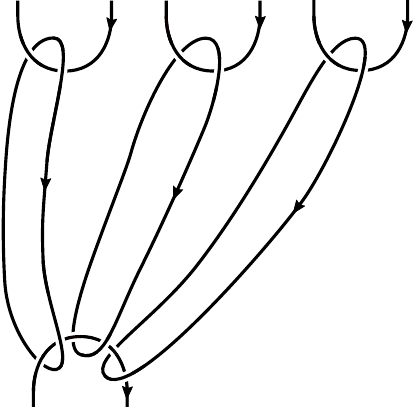}\;.&
	\end{align*} \endgroup
\end{Lemma}
The proof of all items in Lemma \ref{kirbybash2} follows by Lemma \ref{kirbybash}$(a)$, the generalized Kirby move~$\textbf{K}2'$ and isotopy.

\subsection[Passing to X\^{}\{bullet\} and X\_\{bullet\}]{Passing to $\boldsymbol{X^\bullet}$ and $\boldsymbol{X_\bullet}$} \label{generalsurf}
Let~$\{\Sigma_{n+1}\}_{n\ge 0}$ be any family of closed oriented surfaces. For each~$n$, there exists an orientation preserving homeomorphism~$f_n \co \Sigma_{n+1} \to S_{n+1}$. Denote by~$\Cyl({f_n}) \co \Sigma_{n+1} \to S_{n+1}$ the associated morphism in~$\textbf{3}\Cob_0$, given by the quadruple
\[
\Cyl({f_n}) = \bigl(C_{S_{n+1}}=S_{n+1}\times [0,1], h_n \co(-\Sigma_{n+1}) \sqcup S_{n+1}\to \partial\bigl(C_{S_{n+1}}\bigr), \Sigma_{n+1},S_{n+1}\bigr),
\]
where~$h_n(x)=(f_n(x),0)$, if~$x \in \Sigma_{n+1}$ and~$
h_n(x)=(x,1)$, if~$x \in S_{n+1}$. It follows from \cite[Section~\RomanNumeralCaps{4}.5.1]{turaevqinvariants}, that the cobordism $\Cyl({f_n})$ is determined up to isotopy.

We pass from $Y^\bullet$ to the cocyclic object $X^\bullet$ in $\textbf{3}\Cob_0$ as follows.
First, for any $n\in \N$, define $X^n=\Sigma_{n+1}$.
Next, the cofaces~$\{X^\bullet(\delta_i^n) \co \Sigma_n \to \Sigma_{n+1}\}_{n\in \N^*, 0\le i \le n}$, the codegeneracies \smash{$\bigl\{X^\bullet\bigl(\sigma_j^n\bigr) \co \Sigma_{n+2} \to \Sigma_{n+1}\bigr\}_{n\in \N, 0\le j \le n}$}, and the cocyclic operators~$\{X^\bullet(\tau_n) \co \Sigma_{n+1} \to \Sigma_{n+1}\}_{n\in\N}$ are defined by formulas
\begin{align*}
	&X^\bullet(\delta_i^n)=(\Cyl({f_n}))^{-1}Y^\bullet(\delta_i^n)\Cyl({f_{n-1}}), \\
	&X^\bullet\bigl(\sigma_j^n\bigr)=(\Cyl({f_n}))^{-1}Y^\bullet\bigl(\sigma_j^n\bigr)\Cyl({f_{n+1}}), \\
	&X^\bullet(\tau_n)=(\Cyl({f_n}))^{-1}Y^\bullet(\tau_n)\Cyl({f_n}).
\end{align*}
It follows from definitions that the family of cylinders~$\{\Cyl({f_n})\co \Sigma_{n+1} \to S_{n+1}\}_{n\in \N}$ is a natural isomorphism between cocyclic objects~$X^\bullet$ and $Y^\bullet$ in $\textbf{3}\Cob_0$. One similarly passes from~$Y_\bullet$ to a~cyclic object~$X_\bullet$ in~$\textbf{3}\Cob_0$.

\section{Preliminaries on monoidal categories and Hopf algebras}\label{catprelimini}

In this section, we recall some algebraic preliminaries on ribbon categories and their graphical calculus as well as categorical Hopf algebras and related concepts. We will mostly use conventions and notations of \cite{moncatstft}.

\subsection{Conventions} In what follows, we suppress in our formulas the associativity and unitality constraints of the monoidal category. This does not lead to ambiguity since by Mac Lane's coherence
theorem~\cite{maclane63}, all legitimate ways of inserting these constraints give the same results. We denote by $\tens$ and~$\uu$ the monoidal product and unit object of a~monoidal category. For any objects $X_1, \dots, X_n$ of a~monoidal category with $n\ge 2$, we set
\[
X_1 \tens X_2 \tens \cdots \tens X_n=(\cdots ((X_1\tens X_2)\tens X_3)\tens \cdots \tens X_{n-1} )\tens X_n
\]
and similarly for morphisms.
A monoidal category is $\kk$-linear, if its Hom sets have a structure of a $\kk$-module such that the composition and monoidal product of morphisms are $\kk$-bilinear. For shortness, we often use the term monoidal $\kk$-category.

\subsection{Braided categories and graphical calculus} \label{braidedcats}

In this section, we briefly recall some conventions on braided categories and their graphical calculus, which were introduced and developed by Joyal and Street in \cite{joyal1985braided, joyal1986braided, joyal1991geometry}.
A \textit{braiding} of a monoidal category $(\mc B, \tens, \uu)$ is a family~$\tau=\{\tau_{X,Y}\co X\tens Y \to Y\tens X\}_{X,Y\in \Ob{\mc B}}$ of natural isomorphisms such that $\tau_{X,Y\tens Z}= (\id_Y \tens \tau_{X,Z})(\tau_{X,Y} \tens \id_Z)$ and $\tau_{X\tens Y, Z}= (\tau_{X,Z} \tens \id_Y)(\id_X \tens \tau_{Y,Z})$ hold for all objects~$X$, $Y$, $Z$ in $\mc B$.
A~\textit{braided category} is a monoidal category endowed with a~braiding. A braiding~$\tau$ of~$\mc B$ is \textit{symmetric} if for all~$X,Y \in \Ob{\mc B}$, $\tau_{Y,X}\tau_{X,Y}=\id_{X\tens Y}$.
A~\textit{symmetric category} is a monoidal category endowed with a~symmetric braiding. For example, the category of left~$\kk$-modules is symmetric.
A \textit{twist} for a braided monoidal category~$\mc B$ (see~\cite{shum1994tortile}) is a natural isomorphism $\theta=\{\theta_X \co X \to X\}_{X \in \Ob{\mc B}}$ such that for all $X,Y \in \Ob{\mc B}$,
\begin{equation*} 
\theta_{X\tens Y}= \tau_{Y,X}\tau_{X,Y} (\theta_X \tens \theta_Y).
\end{equation*}
By a \textit{balanced category}, we mean a braided category endowed with a twist. For example, the family $\{\id_X\co X \to X\}_{X\in \Ob{\mc B}}$ is a twist for $\mc B$ if and only if $\mc B$ is symmetric. Also, any ribbon category (see Section~\ref{pivotal}) has a canonical twist.

In this paper, we intensively use the \textit{Penrose graphical calculus}, which allows us to avoid lengthy algebraic computations by using simple topological arguments.
The diagrams read from bottom to top.
In a monoidal category~$\mc B$, the diagrams are made of arcs colored by objects of~$\mc B$ and of boxes, colored by morphisms of~$\mc B$.
Arcs colored by~$\uu$ may be omitted in the pictures.
The identity morphism of an object~$X$, a morphism~$f\co X\to Y$ in~$\mc B$, and its composition with a morphism~$g\co Y\to Z$ in~$\mc B$ are represented respectively as
\[
\put(4,-34){\small $X$}
\raisebox{-13mm}{\includegraphics{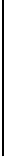}}\, ,
\put(36,-35){\small $X$}
\put(36,30){\small $Y$}
\put(30,-2){\small $f$}
\qquad
\raisebox{-13mm}{\includegraphics{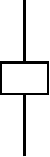}}\, ,
\put(36,-35){\small $X$}
\put(36,30){\small $Z$}
\put(36,-2.5){\small $Y$}
\put(31,14.5){\small $g$}
\put(30,-17){\small $f$}
\qquad
\raisebox{-13mm}{\includegraphics{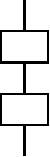}}\, .
\]
The tensor product of two morphisms~$f\co X\to Y$ and~$g\co U \to V$ is represented by placing a~picture of~$f$ to the left of the picture of~$g$.
Any diagram represents a~morphism and the latter depends only on the isotopy class of the diagram representing it.
When $\mc B$ is braided with a~braiding $\tau$, we exceptionally depict
\[
\tau_{X,Y}=\,
\put(-2,-25){\small $X$}
\put(19,23){\small $X$}
\put(-1,23){\small $Y$}
\put(19,-25){\small $Y$}
\raisebox{-10mm}{\includegraphics[scale=0.8]{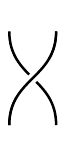}}\,
\qquad \text{and} \qquad
\tau_{X,Y}^{-1}=\,
\put(-2,-25){\small $Y$}
\put(19,23){\small $Y$}
\put(-2,23){\small $X$}
\put(19,-25){\small $X$}
\raisebox{-10mm}{\includegraphics[scale=0.8]{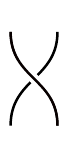}}\,
.
\]
When~$\mc B$ is a balanced category~$\theta=\{\theta_X\co X\to X\}_{X\in \Ob{\mc B}}$, we depict
\[
%
\theta_X = \,
\put(7,-14.5){\small $X$}
\raisebox{-5.5mm}{\includegraphics[scale=0.8]{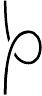}}\,
\qquad \text{and} \qquad
\left(\theta_X\right)^{-1} = \,
\put(7,-14.5){\small $X$}
\raisebox{-5.5mm}{\includegraphics[scale=0.8]{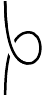}}\, .
\]
We warn the reader that this notation should not be confused with notation of a left twist in a~ribbon category (see Section~\ref{pivotal}). We made this choice of notation since any ribbon category is a~particularly important example of a balanced category.

\subsection{Pivotal categories and graphical calculus} \label{pivotal}

	A \textit{pivotal category} is a monoidal category~$\mc C$ such that to any object~$X$ of~$\mc C$ is associated a dual object $X^*\in \Ob{\mc C}$ and four morphisms
\begin{alignat*}{3}
	&\evl_X \co \ X^* \tens X \to \uu, \qquad&& \coevl_X \co \ \uu \to X\tens X^*,& \\
	&\evr_X \co \ X \tens X^* \to \uu, \qquad&& \coevr_X \co \ \uu \to X^* \tens X,&
\end{alignat*}
	satisfying several conditions and such that the so called left and right duality functors coincide as monoidal functors. The latter implies in particular that the dual morphism $f^* \co Y^* \to X^*$ of a morphism $f\co X\to Y$ in $\mc C$ is computed by
	\begin{align*}
\begin{split}
	f^*&=(\id_{X^*} \tens \evr_Y )(\id_{X^*} \tens f \tens \id_{Y^*})(\coevr_X \tens \id_{Y^*}) \\
	& =(\evl_Y \tens \id_{{X^*}})(\id_{{Y^*}} \tens f \tens \id_{{X^*}})(\id_{{Y^*}} \tens \coevl_X).
\end{split}
	\end{align*}
	The graphical calculus for monoidal categories from Section~\ref{braidedcats} is extended to pivotal cate\-gories by orienting arcs. If an arc colored by~$X$ is oriented upwards, the represented object in source/target of corresponding morphism is~$X^*$. For example,~$\id_{X}$, $\id_{X^*}$, and a morphism $f\co X\tens Y^* \tens Z \to U \tens V^*$ are depicted by
\[
\id_{X}=\,
\put(6,-18.5){\small $X$}
\raisebox{-7mm}{\includegraphics[scale=1]{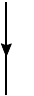}}\, ,
\qquad
\id_{X^*}=\,
\put(6,-18.5){\small $X$}
\raisebox{-7mm}{\includegraphics[scale=1]{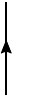}}\,
=\,
\put(6,-18.5){\small $X^*$}
\raisebox{-7mm}{\includegraphics[scale=1]{page18-01}}\ ,
\qquad
f=\,
\put(7,-18.5){\small $X$}
\put(31,-18.5){\small $Y$}
\put(52,-18.5){\small $Z$}
\put(15.5,17.5){\small $U$}
\put(36.5,17.5){\small $V$}
\put(24.5,0){\small $f$}
\raisebox{-7mm}{\includegraphics[scale=1]{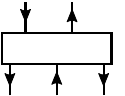}}\, .
\]
The morphisms $\evl_X$, $\evr_X$, $\coevl_X$, and $\coevr_X$ are respectively depicted by
\[
\put(25.5,-8){\small $X$}
\raisebox{-3mm}{\includegraphics[scale=1]{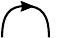}}\ ,
\put(47.5,-8){\small $X$}
\qquad
\raisebox{-3mm}{\includegraphics[scale=1]{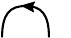}}\ ,
\put(47,2.3){\small $X$}
\qquad
\raisebox{-3mm}{\includegraphics[scale=1]{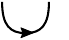}}\ ,
\put(47,2.3){\small $X$}
\qquad
\raisebox{-3mm}{\includegraphics[scale=1]{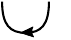}}\ .
\]

Let~$\mc B$ be a braided pivotal category.
The \emph{left twist} of an object~$X$ of $\mc B$ is defined by
\[
\theta_X^l=
\put(6,-16){\small $X$}
\,\raisebox{-6mm}{\includegraphics[scale=0.85]{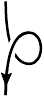}}\,
= (\id_X\tens \evr_X)(\tau_{X,X}\tens \id_{X^*})(\id_X\tens \coevl_X)\co\ X\to X,\]
while the \emph{right twist} of $X$ is defined by
\[
\theta_X^r=
\put(19,-16){\small $X$}
\, \raisebox{-6mm}{\includegraphics[scale=0.85]{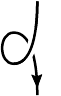}}\,
=(\evl_X\tens \id_{X})(\id_{X^*}\tens \tau_{X,X})(\coevr_X\tens \id_X)\co\ X\to X.\]
The left and the right twist are natural isomorphisms with inverses
\[
\bigl(\theta_X^l\bigr)^{-1}=
\put(6.5,-16){\small $X$}
\,\raisebox{-6mm}{\includegraphics[scale=0.85]{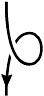}}\,
\qquad \text{and} \qquad
\bigl(\theta_X^r\bigr)^{-1}=
\put(19,-16){\small $X$}
\,\raisebox{-6mm}{\includegraphics[scale=0.85]{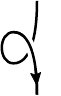}}\, .
\]
A \textit{ribbon category} is a braided pivotal category $\mc B$ such that $\theta_X^l=\theta_X^r$ for all $X \in \Ob{\mc B}$.
In this case, the family~\smash{$\theta=\bigl\{\theta_X=\theta_X^l=\theta_X^r \co X\to X\bigr\}_{X\in \Ob{\mc{B}}}$} is a twist in the sense of Section~\ref{braidedcats} and is called
the {\it twist} of $\mc B$.

Finally, let~$\mc B$ be a ribbon category. A ribbon graph~$\Gamma$ (recall Section~\ref{ribbongrph}) is~\emph{$\mc B$-colored} if each arc and circle of~$\Gamma$ is endowed with an object of~$\mc B$ and each coupon of~$\Gamma$ is endowed with a~morphism in~$\mc B$ from the object determined by its bottom base (as in Section~\ref{braidedcats}) to the object determined by its top base. By~\cite[Theorem~2.5]{turaevqinvariants}, any~$\mc B$-colored ribbon graph determines a morphism in~$\mc B$.
The colors of~$\mc B$-colored ribbon graphs are shown on their diagrams. For example, given objects~$X,Y,Z,T \in \Ob{\mc B}$ and a morphism~$f \colon Y^* \otimes X \to Z$ in ${\mc B}$,
the diagram
\[
\put(-1,-32){\small $X$}
\put(30,-32){\small $Y$}
\put(51,-32){\small $T$}
\put(22,13){\small $Z$}
\put(15,-2){\small $f$}
\raisebox{-12mm}{\includegraphics[scale=0.85]{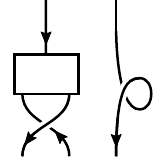}}
\]
represents a $\mc B$-colored ribbon $(3,2)$-graph whose associated morphism in $\mc B$ is
\[
(f \circ \tau_{X,Y^*}) \otimes \theta_T \colon \ X \otimes Y^* \otimes T \to Z \tens T.
\]

\subsection{Categorical traces and dimensions} 
Let~$\mc C$ be a pivotal category.
The~\textit{left trace} $\ltr{f}$ and the~\textit{right trace}~$\rtr{f}$ of a morphism $f\co X \to X$ in~$\mc C$, are respectively defined as
\begin{align*}
&\ltr{f}=\evl_{X}(\id_{X^*}\tens f)\coevr_X \in \End(\uu) \qquad \text{and} \\
&\rtr{f}=\evr_{X}(f\tens \id_{X^*})\coevl_X \in \End(\uu).
\end{align*}
The left dimension $\dim_l(X)$ and the right dimension $\dim_r(X)$ of an object $X$ in $\mc C$ are respectively defined as
\[\dim_l(X)=\ltr{\id_X} \qquad \text{and} \qquad \dim_r(X)=\rtr{\id_X}.\]
A~\textit{spherical category} is a pivotal category~$\mc C$ such that left and right trace of any endomorphism in~$\mc C$ coincide. In particular,~$\dim_l(X)=\dim_r(X)$ holds for all~$X \in\Ob{\mc C}$. In a spherical category, we omit~$l$ and~$r$ from the notation for trace and dimension. Any ribbon category is spherical, see~\cite[Corollary 3.4]{moncatstft}.

\subsection{Modular categories} \label{modular} An object~$X$ of a~$\kk$-linear category $\mc C$ is \textit{simple} if the~$\kk$-module~$\End_{\mc C}(X)$ is free of rank $1$. In that case, the map $\kk \to \End_{\mc C}(X)$, $k\mapsto k\id_X$ is an isomorphism of $\kk$-algebras. A \textit{prefusion}~$\kk$-\textit{category} is a monoidal~$\kk$-linear category~$\mc C$ such that there is a set $I$ of simple objects of $\mc C$ satisfying the following conditions:
\begin{itemize}\itemsep=0pt
	\item[$(a)$] For any distinct elements $i$, $j$ of $I$, $\Hom_{\mc C}(i,j)=0$.
	\item[$(b)$] The unit object $\uu$ of $\mc C$ is an element of $I$.
	\item[$(c)$] Any object of $\mc C$ is a finite direct sum of elements of $I$.
\end{itemize}
The set $I$ is called a \textit{representative set of simple objects} of $\mc C$.
A~\textit{fusion} $\kk$-category is a rigid prefusion $\kk$-category such that the set of isomorphism classes of simple objects is finite.

Now let $\mc C$ be a pivotal fusion $\kk$-category.
We identify $\kk$ and $\End_{\mc C}(\uu)$ via the $\kk$-linear isomorphism $k\mapsto k\id_{\uu}$.
Pick a representative set~$I$ of simple objects of~$\mc C$.
The \emph{dimension of the category}~$\mc C$ is an element of~$\End_{\mc C}(\uu)\cong \kk$ defined by
\[\dim (\mc C) = \sum_{i\in I} \dim_l(i) \dim_r(i).\]
The dimension of $\mc C$ does not depend on the choice of $I$ since isomorphic objects of $\mc C$ have the same left/right dimensions.
Note that if $\mc C$ is spherical, then
\[
\dim (\mc C)= \sum_{i\in I} (\dim (i))^2.
\]

Finally, let~$\mc C$ be a ribbon fusion $\kk$-category and~$I$ a representative set of simple objects of~$\mc C$. The scalars
\[
\Delta_{\pm}=\sum_{i\in I}\dim(i) \, \mathrm{tr}\bigl(\theta_i^{\pm 1}\bigr) \in \End_{\mc C}(\uu) \cong \kk,
\]
where $\theta$ is the twist of $\mc C$, do not depend on the choice of $I$.
The~$S$-\textit{matrix}~$[S_{i,j}]_{i,j\in I}$ of~$\mc C$ is defined by $S_{i,j}=\text{tr}(\tau_{i,j}\tau_{j,i})\in \End_{\mc C}(\uu) \cong \kk$.
Note that the invertibility of $S$ does not depend on the choice of $I$.
A \textit{modular $\kk$-category} is a ribbon fusion $\kk$-category whose $S$-matrix is invertible. The scalars~$\Delta_+$, $\Delta_-$, and~$\dim(\mc C)$ associated with a modular $\kk$-category $\mc C$ are invertible in~$\kk$ and are related by $\dim(\mc C)=\Delta_-\Delta_+$ (see~\cite[p.~89]{turaevqinvariants}).
A modular $\kk$-category is \textit{anomaly free} if~$\Delta_+=\Delta_-$.

\subsection{Categorical Hopf algebras, pairings, and integrals} \label{braidedHopf}
In this section, we review categorical algebras and bialgebra pairings. See \cite{majid1994algebras, majid2000foundations} for details.
An~\textit{algebra} in a monoidal category~$\mc C$ is a triple~$(A,m,u)$, where~$A$ is an object of~$\mc C$,~$m\co A\tens A \to A$ and~$u \co \uu \to A$ are morphisms in~$\mc C$, called~\textit{multiplication} and~\textit{unit} respectively, which satisfy
$
m(m\tens \id_A)=m(\id_A \tens m)$ and $m(u \tens \id_A)= \id_A= m(\id_A\tens u).
$
The multiplication and the unit are depicted by
\[
m=
\put(27,-15){\small $A$}
\,\raisebox{-5.5mm}{\includegraphics[scale=1]{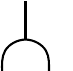}}\, ,
\qquad
u =
\put(12.5,-14.8){\small $A$}
\,\raisebox{-6mm}{\includegraphics[scale=1]{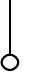}}\, .
\]
An \textit{algebra morphism} between two algebras $(A,m_A,u_A)$ and $(B,m_B,u_B)$ in a monoidal category~$\mc C$ is a morphism $f\co A \to B$ in~$\mc C$ such that $fm_A=m_B(f\tens f)$ and $fu_A=u_B$.

Dually, a~\textit{coalgebra} in a monoidal category~$\mc C$ is a triple~$(C, \Delta, \varepsilon)$, where~$C$ is an object of~$\mc C$,~$\Delta \co C\to C\tens C$ and~$\varepsilon\co C \to \uu$ are morphisms in~$\mc C$, called~\textit{comultiplication} and~\textit{counit} respectively, which satisfy
$
(\Delta\tens \id_C)\Delta= (\id_C \tens \Delta)\Delta$ and $(\id_C \tens \varepsilon)\Delta=\id_C= (\varepsilon \tens \id_C)\Delta.
$
The comultiplication and the counit are depicted by
\[
\Delta=
\put(17,-11.5){\small $C$}
\,\raisebox{-4.5mm}{\includegraphics[scale=1]{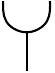}}\, ,
\qquad
\varepsilon=
\put(9.5,-11.5){\small $C$}
\,\raisebox{-4.5mm}{\includegraphics[scale=1]{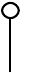}}\, .
\]
A \textit{coalgebra morphism} between two coalgebras $(C,\Delta_C,\vareps_C)$ and $(D,\Delta_D,\vareps_D)$ in a monoidal cate\-gory~$\mc C$ is a morphism $f\co C \to D$ in $\mc C$ such that $\Delta_D f=(f\tens f)\Delta_C$ and $\vareps_D f=\vareps_C$.

From now on, let $\mc B$ be a braided monoidal category.
A \textit{bialgebra} in a $\mc B$ is a quintuple~$(A,m,u, \Delta, \varepsilon)$ such that $(A,m,u)$ is an algebra in $\mc B$, $(A,\Delta, \varepsilon)$ is a coalgebra in $\mc B$, and such that the following additional relations hold: $\Delta m= (m\tens m)(\id_A\tens \tau_{A,A} \tens \id_A)(\Delta\tens \Delta)$, $\varepsilon m =\varepsilon
\tens \varepsilon$, $\Delta u = u\tens u$, and $\varepsilon u= \id_{\uu}$.
A \textit{bialgebra morphism} between two bialgebras~$A$ and~$B$ in a braided monoidal category~$\mc B$ is a morphism $A \to B$ in $\mc B$, which is both an algebra and a~coalgebra morphism.

A \textit{Hopf algebra} in~$\mc B$ is a sextuple~$(A,m,u, \Delta, \varepsilon, S)$, where~$(A,m,u, \Delta, \varepsilon)$ is a bialgebra in~$\mc B$ and~$S\co A \to A$ is an isomorphism in~$\mc B$, called the \textit{antipode}, which satisfies~$m(S\tens \id_A)\Delta=u \varepsilon = m(\id_A \tens S)\Delta$.
The antipode and its inverse are depicted by
\[
S =
\put(8.5,-11.5){\small $A$}
\,\raisebox{-4.5mm}{\includegraphics[scale=1]{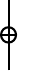}}\, ,
\qquad
S^{-1}=
\put(8.5,-11.5){\small $A$}
\,\raisebox{-4.5mm}{\includegraphics[scale=1]{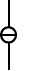}}\, .
\]
A \textit{Hopf algebra morphism} between two Hopf algebras is a bialgebra morphism between them.

A~\textit{bialgebra pairing} for a bialgebra $(A,m,u, \Delta, \varepsilon)$ in $\mc B$ is a morphism~$\omega\co A \tens A \to \uu$ in $\mc B$ such that
\begin{alignat*}{3}
&\omega(m\tens \id_A)=\omega(\id_A\tens \omega\tens \id_A)(\id_{A\tens A}\tens \Delta), \qquad&& \omega(u\tens \id_A)= \varepsilon,&\\
&\omega(\id_A\tens m)=\omega(\id_A \tens \omega \tens \id_A)(\Delta \tens \id_{A\tens A}), \qquad&& \omega (\id_A \tens u) = \varepsilon.&
\end{alignat*}
A bialgebra pairing~$\omega$ for~$A$ is \textit{non-degenerate} if there exists a morphism~$\Omega\co \uu \to A\tens A$ in~$\mc B$ such that $(\omega\tens \id_A)(\id_A\tens\Omega)=\id_A$ and $(\id_A \tens \omega)(\Omega \tens \id_A)=\id_A$.
The morphism~$\Omega$ is called the \textit{inverse} of the pairing~$\omega$.
If~$A$ is Hopf algebra, the pairing~$\omega$ for the underlying bialgebra is called a~\textit{Hopf pairing}.

Finally, a \textit{left} (respectively \textit{right}) \emph{integral} of a bialgebra~$(A,m, \Delta,u, \varepsilon)$ in~$\mc B$ is a morphism~$\Lambda\co \uu \to A$ in $\mc B$ such that
\[m(\id_A\tens \Lambda)= \Lambda \varepsilon, \qquad \text{respectively} \ m(\Lambda \tens \id_A)=\Lambda \varepsilon.\]
Dually, a \textit{left} (respectively \textit{right}) \textit{cointegral} of $A$ is a morphism $\lambda \co A \to \uu$ in $\mc B$ such that
\[
(\id_A\tens \lambda )\Delta= u\lambda, \qquad \text{respectively} \ (\lambda \tens \id_A)\Delta= u \lambda.
\]

\subsection{Coend of a category}\label{Coendprelimini}
Let~$\mc B$ be a pivotal category. Let~$F_{\mc B} \co \mc B^{\opp} \times \mc B \to \mc B$ be the functor defined by~$F_{\mc B}(X,Y)=X^*\tens Y$.
	A \textit{dinatural transformation} from~$F_{\mc B}$ to an object~$D$ of~$\mc B$ is a function~$d$ that assigns to any object~$X$ of~$\mc B$ a morphism~$d_X\co X^*\tens X \to D$ such that for all morphisms~$f\co X \to Y$ in $\mc B$,
	\[d_X(f^* \tens \id_{X})=d_Y(\id_{Y^*}\tens f).\]
	\textit{The coend of} $\mc B$, if it exists, is a pair~$(\coend,i)$ where~$\coend$ is an object of~$\mc B$ and~$i$ is a dinatural transformation from~$F_{\mc B}$ to~$\coend$, which is universal among all dinatural transformations.
	More precisely, for any dinatural transformation~$d$ from~$F_{\mc B}$ to~$D$, there exists a unique morphism~$\phi\co \coend \to D$ in~$\mc B$ such that~$d_X=\phi i_X$ for all~$X \in \Ob{\mc B}$.
	A coend~$(\coend,i)$ of a category~$\mc B$, if it exists, is unique up to a unique isomorphism commuting with the dinatural transformation.
We depict the dinatural transformation $i=\{i_X\co X^*\tens X \to \coend \}_{X\in \Ob{\mc B}}$ as
\[
i_X=
\put(26.5,-21){\small $X$}
\put(17,5){\small $\coend$}
\,\raisebox{-8mm}{\includegraphics[scale=1]{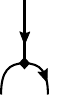}}\, .
\]
An important and well-known factorization property is given in the following lemma.
\begin{Lemma}[Fubini theorem for coends, \cite{catswork}] \label{CoendFubini} Let~$(\coend,i)$ be a coend of a braided pivotal ca\-te\-go\-ry~$\mc B$.
	If~$d=\{d_{X_1,\dots, X_n}\co X_1^*\tens X_1\tens \cdots \tens X_n^* \tens X_n \to D\}_{X_1,\dots, X_n\in \Ob{\mc B}}$ is a family
	of morphisms in~$\mc B$, which is dinatural in each~$X_i$ for~$1\leq i \leq n$, then there exists a unique morphism~$\phi\co \coend^{\tens n} \to D$ in~$\mc B$ such that~$d_{X_1,\dots,X_n}=\phi(i_{X_1}\tens \cdots \tens i_{X_n})$
	for all~$X_1,\dots, X_n \in \Ob{\mc B}$.
\end{Lemma}
According to~\cite{LYUBASHENKO1995279, {MAJID1993187}}, coend of a braided pivotal category~$\mc B$ is a Hopf algebra in~$\mc B$ endowed with a canonical Hopf pairing. Its unit is~$u=(\id_\uu \tens i_{\uu})(\coevl_{\uu} \tens \id_\uu) \co \uu \to \coend$.
Multiplication~$m \co \coend \tens \coend \to \coend$ and canonical pairing~$\omega \co \coend \tens \coend \to \uu$ are unique morphisms such that for all~$X,Y \in \Ob{\mc B}$,
\[
\put(36,31){\small $\coend$}
\put(15,-13){\small $\coend$}
\put(57,-13){\small $\coend$}
\put(24,-30){\small $X$}
\put(67,-30){\small $Y$}
\put(28.3,9.5){\small $m$}
\raisebox{-11mm}{\includegraphics[scale=0.95]{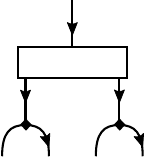}}\, =
\put(68,33){\small $\coend$}
\put(31,-30){\small $X$}
\put(9,-30){\small $X$}
\put(73,-30){\small $Y$}
\put(52,-30){\small $Y$}
\put(18,-1){\small $\id_{Y\tens X}$}
\put(61,-4.5){\small $\id_{Y\tens X}$}
\ \raisebox{-11mm}{\includegraphics[scale=0.95]{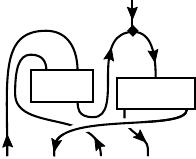}}\, ,
\put(35,-9){\small $\coend$}
\put(78,-9){\small $\coend$}
\put(45,-29){\small $X$}
\put(88,-29){\small $Y$}
\put(51,19){\small $\omega$}
\qquad
\raisebox{-11mm}{\includegraphics[scale=0.95]{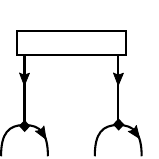}}\, =
\put(29,-30){\small $X$}
\put(69,-30){\small $Y$}
\ \raisebox{-11mm}{\includegraphics[scale=0.95]{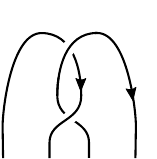}}\, .
\]
Its comultiplication~$\Delta \co \coend\to \coend \tens \coend$, counit~$\varepsilon \co \coend \to \uu$, and antipode~$S\co \coend\to \coend$ are unique morphisms such that for all~$X \in \Ob{\mc B}$,
\[
\put(32,31){\small $\coend$}
\put(4,31){\small $\coend$}
\put(24,-9){\small $\coend$}
\put(33,-30){\small $X$}
\put(17,16){\small $\Delta$}
\raisebox{-11mm}{\includegraphics[scale=0.95]{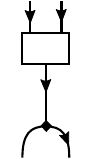}}\, =
\put(16,-30){\small $X$}
\put(59,-3){\small $X$}
\put(81,-30){\small $X$}
\put(30,18){\small $\coend$}
\put(72.2,18){\small $\coend$}
\raisebox{-11mm}{\includegraphics[scale=0.95]{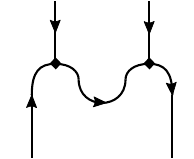}}\, ,
\put(36.5,-9){\small $\coend$}
\put(45.5,-30){\small $X$}
\put(29.9,15.7){\small $\varepsilon$}
\qquad
\raisebox{-11mm}{\includegraphics[scale=0.95]{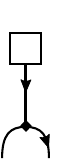}}\, =
\put(27,-30){\small $X$}
\ \raisebox{-11mm}{\includegraphics[scale=0.95]{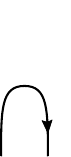}}\, ,
\put(37,31){\small $\coend$}
\put(37,-10){\small $\coend$}
\put(45,-30){\small $X$}
\put(29.2,8.2){\small $S$}
\qquad
\raisebox{-11mm}{\includegraphics[scale=0.95]{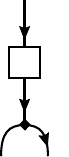}}\, =
\put(2,-30){\small $X$}
\put(36,-30){\small $X$}
\put(28,18){\small $\coend$}
\ \raisebox{-11mm}{\includegraphics[scale=0.95]{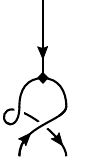}}\, .
\]
Two useful properties of antipode of coend $\coend$ are
\begin{align}
&S^2=\theta_{\coend}^r \qquad \text{and} \label{invol} \\ &\omega(S\tens \id_\coend)=\omega(\id_\coend \tens S). \label{omegaS}
\end{align}
Also, we have by definitions of $\omega$ and $S$ that $X,Y \in \Ob{\mc B}$,
\begin{equation*}
 \put(15.5,-6){\small $\coend$}
 \put(58.5,-6){\small $\coend$}
 \put(25,-27){\small $X$}
 \put(68,-27){\small $Y$}
 \put(29.6,22){\small $\omega$}
 \raisebox{-10mm}{\includegraphics[scale=0.95]{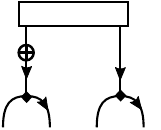}}\, =
 \put(15.5,-6){\small $\coend$}
 \put(58.5,-6){\small $\coend$}
 \put(25,-27){\small $X$}
 \put(68,-27){\small $Y$}
 \put(29.9,21.5){\small $\omega$}
 \raisebox{-10mm}{\includegraphics[scale=0.95]{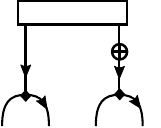}}\, =
 \put(27.5,-27){\small $X$}
 \put(67.5,-27){\small $Y$}
 \ \raisebox{-10mm}{\includegraphics[scale=0.95]{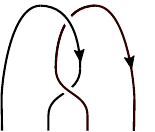}}\, .
\end{equation*}

\section{Cyclic modules from (co)algebras} \label{algebraiccyclic}

Let~$\mc B$ be a balanced~$\kk$-linear category. In this section, we first review the construction of (co)cyclic~$\kk$-modules from coalgebras and algebras in~$\mc B$. Our setting is only a particular case of constructions by Akrami and Majid~\cite{cycliccocycles}, since any algebra in a balanced category is a ribbon algebra in the sense of~\cite{cycliccocycles}. However, the cocyclic $\kk$-module from \cite{cycliccocycles} is here viewed as cyclic dual of a certain cyclic $\kk$-module.
Proofs are for convenience of context also given in \cite{bartulovic2022cyclicsetsribbonstring}, where these were merely (co)cyclic sets, since one dropped the hypothesis of~$\kk$-linearity of~$\mc B$. Further, we outline some basic computations of cyclic (co)homology of some of the introduced (co)cyclic~$\kk$-modules.
Finally, in Section~\ref{hopfexplicit}, we explicit the (co)faces, (co)degeneracies, and (co)cyclic operators of (co)cyclic~$\kk$-modules associated to the coend of the representation category of a~finite dimensional ribbon Hopf algebra. For more details on Hochschild and cyclic (co)homology of (co)cyclic $\kk$-modules, see \cite[Section~9.6]{weibel}.

\subsection{Cocyclic modules from coalgebras} \label{objetdedebut}

Any coalgebra~$C$ in~$\mc B$ gives rise to a cocyclic~$\kk$-module~$C^\bullet$ as follows. For any~$n \in \N$, define ${C^n=\Hom_\mc{B}\bigl(C^{\tens n+1}, \uu\bigr)}$. Next, define the cofaces \smash{$\bigl\{\delta_i^n \co C^{n-1} \to C^n\bigr\}{}_{n\in \N^*, 0\le i \le n}$}, the codegeneracies \smash{$\bigl\{\sigma_j^n \co C^{n+1} \to C^n\bigr\}_{n\in \N, 0\le j \le n}$}, and the cocyclic operators~$\{\tau_n \co C^n \to C^n\}_{n\in \N}$ by set\-ting
	\begin{gather*}
	\delta_i^n (f) =
 \put(10.5,-15){\small $\cdots$}
 \put(53.5,-15){\small $\cdots$}
 \put(35.5,15){\small $f$}
 \put(3.5,-25){\small $0$}
 \put(35.5,-25){\small $i$}
 \put(67,-25){\small $n$}
 \,\raisebox{-10mm}{\includegraphics[scale=0.95]{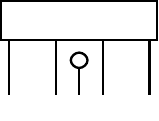}}\, ,
	\qquad \!
 \sigma_j^n (f) =
 \put(17,-15){\small $\cdots$}
 \put(54,-15){\small $\cdots$}
 \put(38.5,15){\small $f$}
 \put(3.5,-25){\small $0$}
 \put(38.5,-25){\small $j$}
 \put(74,-25){\small $n$}
 \,\raisebox{-10.7mm}{\includegraphics[scale=0.95]{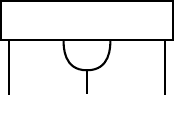}}\, ,
 \qquad \!
 \tau_n(f)=
 \put(25,-15){\small $\cdots$}
 \put(35,15){\small $f$}
 \put(3.5,-25){\small $0$}
 \put(41,-25){\small $n-1$}
 \put(69,-25){\small $n$}
 \,\raisebox{-9.5mm}{\includegraphics[scale=0.95]{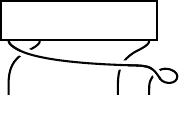}}\, .
	\end{gather*}
	An integer~$k$ below an arc denotes the~$k$-th tensorand of a tensor power of~$C$. This construction is functorial in~$C$, that is, a morphism between coalgebras in $\mc B$ induces the morphism of corresponding cocyclic $\kk$-modules.

	\subsection{Cyclic modules from algebras}
	\label{sec-alg}
	Any algebra~$A$ in~$\mc B$ gives rise to a cyclic $\kk$-module~$A_\bullet$ as follows. For any~$n \in \N$, define $A_n=\Hom_\mc{B}\bigl(A^{\tens n+1}, \uu\bigr)$.
	Next, define the faces~\smash{$\{d_i^n\co A_n \to A_{n-1} \}_{n\in \N^*, 0\le i \le n}$}, the dege\-neracies~\smash{$\bigl\{s^n_j\co A_n \to A_{n+1}\bigr\}_{n \in \N, 0\le j \le n}$}, and cyclic operators~$\{t_n\co A_n \to A_n\}_{n\in \N}$ by setting
\begingroup
\allowdisplaybreaks
	\begin{gather*}
	d_0^n(f)=
 \put(22,-16){\small $\cdots$}
 \put(18.8,13){\small $f$}
 \put(14.4,-26){\small $0$}
 \put(26,-26){\small $n-1$}
 \,\raisebox{-10.5mm}{\includegraphics[scale=0.95]{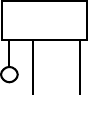}}\, ,
 \qquad
 d_i^n(f)=
 \put(10.5,-16){\small $\cdots$}
 \put(53.5,-16){\small $\cdots$}
 \put(35,13){\small $f$}
 \put(3.5,-26){\small $0$}
 \put(17.4,-26){\small $i-1$}
 \put(46.7,-26){\small $i$}
 \put(58.2,-26){\small $n-1$}
 \,\raisebox{-10.5mm}{\includegraphics[scale=0.95]{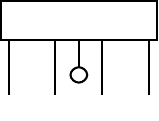}}\, ,
 \qquad
	d_n^n(f)=
 \put(10.5,-16){\small $\cdots$}
 \put(18.8,13){\small $f$}
 \put(3.5,-26.3){\small $0$}
 \put(15.5,-26.3){\small $n-1$}
 \,\raisebox{-10.4mm}{\includegraphics[scale=0.95]{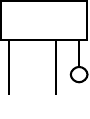}}\, , \\
	s_j^n(f)=
 \put(12.5,-16){\small $\cdots$}
 \put(59.5,-16){\small $\cdots$}
 \put(38,15){\small $f$}
 \put(3.5,-24){\small $0$}
 \put(28,-24){\small $j$}
 \put(40,-24){\small $j+1$}
 \put(66,-24){\small $n+1$}
 \,\raisebox{-10.4mm}{\includegraphics[scale=0.95]{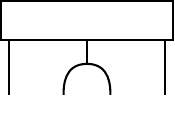}}\, ,
 \qquad
 t_n(f)=
 \put(39,-16){\small $\cdots$}
 \put(36,15){\small $f$}
 \put(3.5,-24){\small $0$}
 \put(17.5,-24){\small $1$}
 \put(67,-24){\small $n$}
 \,\raisebox{-9.5mm}{\includegraphics[scale=0.95]{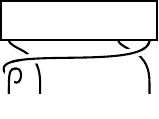}}\, .
	\end{gather*}
\endgroup
This construction is functorial in~$A$, that is, a morphism between algebras in $\mc B$ induces the morphism of corresponding cyclic $\kk$-modules.

	\subsection{Cyclic duals} \label{passing}
	The cyclic duality~$L$ from Section~\ref{connes loday dual} transforms the cocyclic $\kk$-module~$C^\bullet$ from Section~\ref{objetdedebut} into the cyclic~$\kk$-module~$C^\bullet\circ L$. For any~$n\in \N$, $C^\bullet\circ L(n)=C^n=\Hom_{\mc B}\bigl(C^{\tens n+1}, \uu\bigr)$.
	The~faces~\smash{$\bigl\{\tilde{d}_i^n\bigr\}{}_{n\in \N^*, 0\le i \le n}$}, the degeneracies \smash{$\bigl\{\tilde{s}_j^n\bigr\}{}_{n \in \N, 0\le j \le n}$}, and the cyclic ope\-rators \smash{$\bigl\{\tilde{t}_n\bigr\}{}_{n \in \N}$} are computed by formulas
\begin{gather*}
	\tilde{d}_i^n (f) =	
 \put(17,-15){\small $\cdots$}
 \put(54,-15){\small $\cdots$}
 \put(38.5,15){\small $f$}
 \put(3.5,-25){\small $0$}
 \put(32,-25){\small $i+1$}
 \put(74,-25){\small $n$}
 \,\raisebox{-11mm}{\includegraphics[scale=0.95]{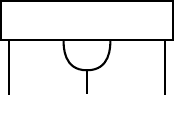}}\, ,
 \qquad
 \tilde{d}_n^n(f)=
 \put(39,-16){\small $\cdots$}
 \put(40,15){\small $f$}
 \put(11,-25){\small $1$}
 \put(74,-25){\small $n$}
 \,\raisebox{-10mm}{\includegraphics[scale=0.95]{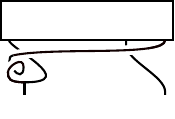}}\, \\
	\tilde{s}_j^n(f)=
 \put(10,-17){\small $\cdots$}
 \put(54,-17){\small $\cdots$}
 \put(34.5,13){\small $f$}
 \put(3.5,-26){\small $0$}
 \put(27,-26){\small $j+1$}
 \put(58,-26){\small $n+1$}
 \,\raisebox{-10.5mm}{\includegraphics[scale=0.95]{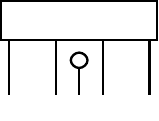}}\, ,
 \qquad
	\tilde{t}_n(f)=
 \put(39,-17){\small $\cdots$}
 \put(36,13){\small $f$}
 \put(3.5,-26){\small $0$}
 \put(17.5,-26){\small $1$}
 \put(67,-26){\small $n$}
 \,\raisebox{-10.5mm}{\includegraphics[scale=0.95]{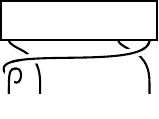}}\, .
	\end{gather*}
	
Similarly, the functor~$L^{\opp}$ transforms the cyclic $\kk$-module~$A_\bullet$ from Section~\ref{sec-alg} into the cocyclic~$\kk$-module~$A_\bullet \circ L^{\opp}$. By definitions,~{$A_\bullet \circ L^\opp(n)=A_n=\Hom_{\mc B}\bigl(A^{\tens n+1}, \uu\bigr)$} for all~${n\in \N}$.
	The cofaces~\smash{$\bigl\{\tilde{\delta}_i^n\bigr\}{}_{n\in \N^*, 0\le i \le n}$}, the codegeneracies~\smash{$\bigl\{\tilde{\sigma}_j^n \bigr\}{}_{n\in \N, 0\le j \le n}$}, and the cocyclic operators~$\{\tilde{\tau}_n\}{}_{n\in \N}$ are computed by formulas
	\begin{gather*}
	\tilde{\delta}_i^n(f) =
 \put(12.5,-16){\small $\cdots$}
 \put(59.5,-16){\small $\cdots$}
 \put(38,14){\small $f$}
 \put(3.5,-25){\small $0$}
 \put(29,-25){\small $i$}
 \put(43,-25){\small $i+1$}
 \put(75,-25){\small $n$}
 \,\raisebox{-11mm}{\includegraphics[scale=0.95]{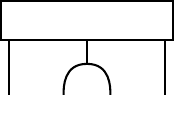}}\, ,
 \qquad
	\tilde{\delta}_n^n(f) =
 \put(34,-8){\small $\cdots$}
 \put(38,14){\small $f$}
 \put(3.5,-25){\small $0$}
 \put(46,-25){\small $n-1$}
 \put(74.5,-25){\small $n$}
 \,\raisebox{-9.9mm}{\includegraphics[scale=0.95]{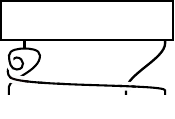}}\, \\
 \tilde{\sigma}_j^n(f)=
 \put(11,-16){\small $\cdots$}
 \put(54,-16){\small $\cdots$}
 \put(35.5,14){\small $f$}
 \put(3.5,-26){\small $0$}
 \put(24.5,-26){\small $j$}
 \put(37.6,-26){\small $j+1$}
 \put(67.5,-26){\small $n$}
 \,\raisebox{-10.5mm}{\includegraphics[scale=0.95]{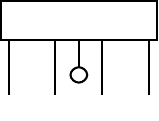}}\, ,
 \qquad
	\tilde{\tau}_n(f)=
 \put(25,-15){\small $\cdots$}
 \put(35,15){\small $f$}
 \put(3.5,-25){\small $0$}
 \put(41,-25){\small $n-1$}
 \put(69,-25){\small $n$}
 \,\raisebox{-10mm}{\includegraphics[scale=0.95]{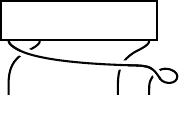}}\, .
	\end{gather*}
	Note that the construction~$A_\bullet \circ L^\opp$ is a particular case of the work of Akrami and Majid~\cite{cycliccocycles}, since any algebra in a balanced category is a ribbon algebra in the sense of~\cite{cycliccocycles}.

\subsection[On the (co)homology of C\^{}\{bullet\} and C\_\{bullet\}]{On the (co)homology of $\boldsymbol{C^\bullet}$ and $\boldsymbol{C_\bullet}$} 
Let~$\kk$ be a commutative ring,~$\mc{B}$ a balanced~$\kk$-category and~$C$ any coalgebra in~$\mc B$. Let~$\alpha \co \uu \to C$ be a morphism such that~$\varepsilon\alpha=\id_\uu$. For example, if~$C$ is a bialgebra in $\mc B$, then we can take $\alpha$ to be the unit of~$C$.
By expanding and writing in the categorical setting~\cite[Remark~1]{khalkhali-rangipour}, we obtain that Hochschild (co)homology of underlying (co)simplicial modules of $C^\bullet$ and $C_\bullet$ appears only in degree zero. Indeed, the~$n$-th Hochschild cohomology $HH^n(C^\bullet)$ of $C^\bullet$ is the $n$-th cohomology of the cochain complex
\[
\xymatrix{\Hom_{\mc B}(C,\uu) \ar[r]^{\beta_1} &\Hom_{\mc B}\bigl(C^{\tens 2},\uu\bigr) \ar[r]^{\beta_2} & \Hom_{\mc B}\bigl(C^{\tens 3},\uu\bigr) \ar[r]^{\hspace{0.8cm}\beta_3} & \cdots },
\]
where \smash{$\beta_{n}= \sum_{i=0}^n (-1)^i \delta_i^n$}. Then the family \smash{$\bigl\{h_n\co \Hom_{\mc B}\bigl(C^{\tens n+1},\uu\bigr)\to \Hom_{\mc B}\bigl(C^{\tens n},\uu\bigr)\bigr\}_{n\in \N^*}$}, defined by setting for any~$f \in \Hom_{\mc B}\bigl(C^{\tens n+1},\uu\bigr)$,
\begin{equation*}
	h_n(f)=
 \put(25.5,-17){\small $\cdots$}
 \put(20.5,12.5){\small $f$}
 \put(6.5,-8.5){\small $\alpha$}
 \put(18.5,-26.5){\small $1$}
 \put(39,-26.5){\small $n$}
 \, \raisebox{-12mm}{\includegraphics[scale=0.95]{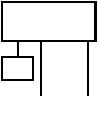}}\, ,
\end{equation*}
satisfies the equalities
\begin{align*}
\begin{split}
&\beta_{n}h_{n}+h_{n+1} \beta_{n+1}= \id_{\Hom_{\mc B}(C^{\tens n+1},\uu)} \quad \text{for } n\ge 1 \qquad \text{and}\\
&h_1\beta_1+\Hom_{\mc B}(\alpha\varepsilon,\uu)= \id_{\Hom_{\mc B}(C,\uu)}.
\end{split}
\end{align*}
As a corollary,~$HH^0(C^\bullet)=\ker(\beta_1)$ and~$HH^n(C^\bullet)=0$ for~$n>0$. From the cohomological form of the Connes' long exact sequence~\cite[Proposition~9.6.11]{weibel} and the fact that Hochschild and cyclic (co)homology always agree in degree~$0$, we easily obtain that $HC^{n}(C^\bullet)\cong \ker(\beta_1)$ for even~$n$ and~$HC^{n}(C^\bullet)\cong 0$ for odd~$n$. A similar statement can be derived for Hochschild and cyclic homologies of~$C_\bullet$. This calculation shows that in a sense,~$C^\bullet$ and~$C_\bullet$ are not interesting from the homological point of view. Therefore, we focus on their cyclic duals in sections that follow.

\subsection{Internal characters}
Let~$\mc B$ be a ribbon~$\kk$-category with a coend~$(\coend,i)$. For any object~$X \in \Ob{\mc B}$, the morphism $\chi_X=i_X \coevr_X\co \uu \to \coend$, also known as internal character~\cite{FGR, shimizu}, enjoys the following trace-like property:
\begin{equation} \label{tracelike}
 \put(4.1,-19.3){\small $\chi_X$}
 \raisebox{-9mm}{\includegraphics[scale=0.85]{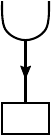}}\, \overset{(i)}{=}
 \put(5.5,-19.8){\small $\chi_X$}
 \,\raisebox{-9.3mm}{\includegraphics[scale=0.85]{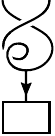}}\, \overset{(ii)}{=}
 \put(7.4,-20){\small $\chi_X$}
 \,\raisebox{-9mm}{\includegraphics[scale=0.85]{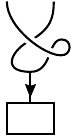}}\, .
\end{equation}
Here~$(i)$ follows by definition of comultiplication of~$\coend$, naturality and definition of twists and braidings, and the isotopy invariance of graphical calculus. Note that a pictorial proof of this fact is given in~\cite[p.~27]{braidedconnes}. The equality~$(ii)$ is obtained by composing both sides of~$(i)$ with~\smash{$\tau_{\coend, \coend}^{-1}\bigl(\id_\coend \tens {\theta}_{\coend}^{-1}\bigr)$}.

Next, the morphism $\psi_X=\omega(\chi_X \tens \id_{\coend})$ satisfies
\begin{equation} \label{tracelike2}
 \put(4,20.5){\small $\psi_X$}
 \raisebox{-9mm}{\includegraphics[scale=0.85]{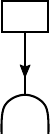}}\, =
 \put(7,20.5){\small $\psi_X$}
 \,\raisebox{-8.8mm}{\includegraphics[scale=0.85]{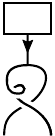}}\, .
\end{equation}
Indeed, we have
\begin{gather*}
\put(3.7,8){\small $\chi_X$}
\put(28,44){\small $\omega$}
\raisebox{-17.5mm}{\includegraphics[scale=0.82]{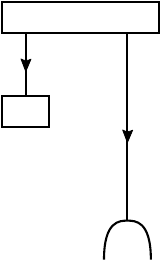}}\, \overset{(i)}{=}
\put(8,-32){\small $\chi_X$}
\put(32,44){\small $\omega$}
\put(30,19){\small $\omega$}
\,\raisebox{-17.5mm}{\includegraphics[scale=0.82]{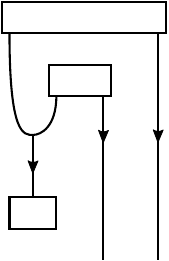}}\, \overset{(ii)}{=}
\put(8,-32){\small $\chi_X$}
\put(32,44){\small $\omega$}
\put(30,19){\small $\omega$}
\,\raisebox{-17.5mm}{\includegraphics[scale=0.82]{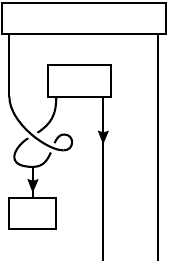}}\, \overset{(iii)}{=}
\put(8,-32){\small $\chi_X$}
\put(32,44){\small $\omega$}
\put(30,19){\small $\omega$}
\,\raisebox{-17.5mm}{\includegraphics[scale=0.82]{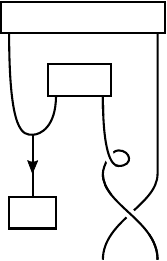}}\, \overset{(iv)}{=}
\put(6,8){\small $\chi_X$}
\put(31,44){\small $\omega$}
\,\raisebox{-17.5mm}{\includegraphics[scale=0.82]{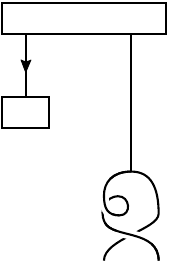}}\, .
\end{gather*}
Here $(i)$ and $(iv)$ follow by the axioms of Hopf pairing $\omega$, $(ii)$ from equation~\eqref{tracelike}, $(iii)$ by equations~\eqref{invol} and \eqref{omegaS}, naturality of braiding, and isotopy.

In the case when~$\mc B$ is a ribbon fusion category with the representative set $I$ of simple objects of~$\mc B$, the family~$\{\chi_i\}_{i\in I}$ is a basis for~$\Hom_{\mc B}(\uu,\coend)$.
Moreover, if the pairing~$\omega$ is non-degenerate,~$\Hom_{\mc B}(\uu,\coend)$ is isomorphic to~$\Hom_{\mc B}(\coend, \uu)$ as $\kk$-module. Combining this with equation~\eqref{tracelike2}, we get that $ HC^0(\coend_\bullet\circ L^\opp) = HH^0(\coend_\bullet\circ L^\opp) = \Hom_{\mc B}(\coend, \uu)$.

\begin{Remark}
Let~$\kk$ be an algebraically closed field and~$\mc B$ a finite tensor~$\kk$-category in the sense of~\cite{etingof-ostrik} with the representative set~$I$ of simple objects (in particular,~$I$ is finite). By~\cite[Theorem~4.1]{shimizu}, the characters~$\{\chi_i\}_{i\in I}$ form a linearly independent set in~$\Hom_{\mc B}(\uu, \coend)$. Under additional unimodularity hypothesis, it is also shown that~$\{\chi_i\}_{i\in I}$ is a basis for~$\Hom_{\mc B}(\uu, \coend)$ if~and only if~$\mc B$ is semisimple. Note that Shimizu works with ends and not coends. However, these theories are essentially the same (see~\cite[Remark 3.12]{shimizu}).
\end{Remark}

\subsection{The coend of the re\-pre\-senta\-tion ca\-te\-go\-ry of a Hopf al\-ge\-bra} \label{hopfexplicit}

Let~$\kk$ be a field and~$H$ a finite dimensional ribbon Hopf algebra over~$\kk$. For the comultiplication, we will use the usual Sweedler notation, that is, one writes~$\Delta(h)= h_{(1)}\tens h_{(2)}$ for any~$h\in H$. Denote by~$R$ and~$\theta$ the~$R$-matrix and the twist element of~$H$. To recall these notions and their properties see \cite{drinfeld1989almost} and \cite{reshetikhin1990ribbon}.
Here we closely follow \cite[Sections 4.2--4.6]{virelizier_2006}. We write
\[
R=\sum_i a_i \tens b_i \in H\tens H \qquad \text{and}\qquad R^{-1}=\sum_i \alpha_i\tens \beta_i \in H\tens H.
 \]
The category $\rep{H}$ of finite-dimensional left $H$-modules is a ribbon category. The coend~$\coend$ of~$\rep{H}$ is a categorical Hopf algebra (see Sections \ref{braidedHopf} and \ref{Coendprelimini}), which is notably studied from topological point of view by Lyubashenko \cite{lyubashenko1995invariants, lyubashenko1995tangles}. Related to this is the Majid's transmutation procedure, which is used to obtain a categorical Hopf algebra from a quasitriangular Hopf algebra. For more details, see \cite{majid1994algebras, majid2000foundations}.
As a $\kk$-module, $\coend$ is equal to~$H^*$ and as a left~$H$-module, it is given by the coadjoint action, that is, for all~$h,k\in H$ and~$f\in H^*$,
\[(h\rtrn f)(k)=f\bigl(S\bigl(h_{(1)}\bigr) k h_{(2)}\bigr).\]

For~$n\in \N^*$, consider the evaluation~$\text{ev}\co H^{\tens n} \to \Hom_\kk\bigl({{H^{*}}^{\tens n}}, \kk\bigr)$ defined by setting for all $X\in H^{\tens n}$ and $f \in {{H^{*}}^{\tens n}}$,
\[
\text{ev}(X)(f)=\langle f,X \rangle.
\]
According to \cite[Lemma 4.5\,$(d)$]{virelizier_2006}, this evaluation induces an isomorphism between $\kk$-modules $\Hom_{\rep{H}}\bigl(\coend^{\tens n}, \kk\bigr)$ and
\[
V_n(H) = \bigl\{X \in H^{\tens n} \mid X\ltrn h = \varepsilon(h)X \text{ for any } h\in H\bigr\}.
\]
Here the right~$H$-action~$\ltrn$ on~$H^{\tens n}$ is defined by setting for any $h\in H$ and any elementary tensor~$X=x_1\tens \cdots \tens x_n \in H^{\tens n}$,
\[
X\ltrn h = S\bigl(h_{(1)}\bigr)x_1h_{(2)} \tens S\bigl(h_{(3)}\bigr)x_2h_{(4)}\tens \cdots \tens S\bigl(h_{(2n-1)}\bigr)x_nh_{(2n)}.
\]
Remark that the $\kk$-module~$V_n(H)$ is equal to the~$0$-th Hochschild homology $HH_0\bigl(H, H^{\tens n}\bigr)$ of $H$ with coefficients in $H^{\tens n}$, where~$H^{\tens n}$ is the bimodule over~$H$ (the left action is given by trivial action via counit).

Under the above isomorphism between~$\Hom_{\rep{H}}\bigl(\coend^{\tens n}, \kk\bigr)$ and~$V_n(H)$, the cyclic~$\kk$-module $\coend^\bullet\circ L$ is identified with the cyclic~$\kk$-module~$\textbf{W}_\bullet$ which is defined as follows. For any~$n\in \N$, set~$\textbf{W}_n=V_{n+1}(H)$.
The faces~$\{d_i^n \co V_{n+1}(H) \to V_{n}(H)\}_{n \in \N^*, 0\le i \le n}$ are given by setting for any elementary tensor~$h_1\otimes\cdots\otimes h_{n+1} \in V_{n+1}(H)$,
\begin{align*}
	&d_i^n(h_1\otimes\cdots\otimes h_{n+1})=
	h_1\otimes h_2\otimes\cdots \otimes h_{i+1}h_{i+2}\otimes\cdots \otimes h_{n+1} \qquad \text{and} \\
	&d_n^n(h_1\otimes \cdots\otimes h_{n+1})=\sum_i\bigl(h_{n+1}\ltrn \bigl(a_i\theta^{-1}\bigr)\bigr)\bigl(h_1\ltrn (b_i)_{(1)}\bigr)\otimes h_2\ltrn (b_i)_{(2)}\otimes \cdots\otimes h_n\ltrn (b_i)_{(n)}.
\end{align*}
The degeneracies \smash{$\bigl\{s_j^n \co V_{n+1}(H) \to V_{n+2}(H)\bigr\}_{n\in \N, 0\le j \le n}$} are given by setting for any elementary tensor $h_1\otimes\cdots\otimes h_{n+1} \in V_{n+1}(H)$,
\[
s_j(h_1\otimes \cdots\otimes h_{n+1})=
	h_1\otimes\cdots\otimes h_{j+1}\otimes 1_H\otimes h_{j+2}\otimes \cdots \otimes h_{n+1}.
\]
The cyclic operators $\{t_n \co V_{n+1}(H) \to V_{n+1}(H)\}_{n\in \N}$ are given by setting for any elementary tensor $h_1\otimes\cdots\otimes h_{n+1} \in V_{n+1}(H)$,
\[
t_n(h_1\tens \cdots \tens h_{n+1})= \sum_i h_{n+1}\ltrn \bigl(a_i\theta^{-1}\bigr)\tens h_1\ltrn (b_i)_{(1)} \otimes \cdots\otimes h_n\ltrn (b_i)_{(n)}.\]

Similarly, the cocyclic~$\kk$-module~$\coend_\bullet\circ L^\opp$ is identified with the cocyclic~$\kk$-module~$\textbf{W}^\bullet$ which is defined as follows.
For any~$n\in \N$, set~$\textbf{W}^n = V_{n+1}(H)$.
The cofaces~$\{\delta_i^n \co V_{n}(H) \to V_{n+1}(H)\}_{n\in \N^*, 0\le i \le n}$ are given by setting for any elementary tensor $h_1\otimes\cdots\otimes h_{n} \in V_n(H)$,
\begin{align*}
&\delta_i^n(h_1\otimes\cdots\otimes h_{n})=
	h_1\otimes h_2\otimes\cdots \otimes \Delta^{\text{Bd}}(h_{i+1}) \otimes\cdots \otimes h_{n} \qquad \text{and} \\
	&\delta_n^n(h_1\otimes \cdots\otimes h_{n})=\sum_i (h_1)_{(2)}^{\text Bd} \ltrn (\beta_i)_{(1)}\tens h_2 \ltrn (\beta_i)_{(2)} \tens \cdots \tens h_{n}\ltrn (\beta_i)_{(n)}\tens (h_1)_{(1)}^{\text Bd} \ltrn(\alpha_i\theta),
\end{align*}
where $\Delta^{\text Bd}$ is a comultiplication on the braided Hopf algebra $H^{\text Bd}$ (see \cite[Lemma 4.4]{virelizier_2006}) asso\-ciated to $H$. As algebras $H^{\text Bd}=H$, but the comultiplication and antipode in $H^{\text Bd}$ are different. Explicitly, for each $h\in H$,
\[
\Delta^{\text Bd}(h)=h_{(2)}a_i\tens S\bigl((b_i)_{(1)}\bigr)h_{(1)}(b_i)_{(2)}.
\]
The codegeneracies \smash{$\bigl\{\sigma_j^n \co V_{n+2}(H) \to V_{n+1}(H)\bigr\}_{n \in \N, 0\le j \le n}$} are given by setting for any elementary tensor~$h_1\otimes\cdots\otimes h_{n+2} \in V_{n+2}(H)$,
\[
\sigma_j^n(h_1\otimes \cdots \otimes h_{n+2})=
h_1\otimes\cdots\otimes h_{j+1}\otimes \varepsilon(h_{j+2})\otimes \cdots \otimes h_{n+2}.\]
The cocyclic operators $\{\tau_n \co V_{n+1}(H) \to V_{n+1}(H)\}_{n\in \N}$ are given by setting for any elementary tensor $h_1\otimes\cdots\otimes h_{n+1} \in V_{n+1}(H)$,
\[\tau_n(h_1\tens \cdots \tens h_{n+1})= \sum_i h_2\ltrn (\beta_i)_{(1)}\tens \cdots \tens h_{n+1}\ltrn (\beta_i)_{(n)}\tens h_1 \ltrn(\alpha_i\theta).\]

\section{Cyclic modules from TQFTs}\label{tftsrel}

In Theorem~\ref{CYCINCOB}, we prove existence of (co)cyclic objects in the category of $3$-cobordisms. By composition, any $3$-dimensional TQFT induces a (co)cyclic $\kk$-module.
In this section we compute it for the Reshetikhin--Turaev TQFT $\RT_{\mc B} \co \textbf{3}\Cob_0 \to \Mod_\kk$ associated to an anomaly free modular category $\mc B$. Note that the coend~$\coend$ of~$\mc B$ exists and is a Hopf algebra in~$\mc B$. Recall the cocyclic~$\kk$-module~$\coend^\bullet$ and the cyclic $\kk$-module~$\coend_\bullet$ (see Section~\ref{algebraiccyclic}) associated to $\coend$, as well as the reindexing involution~$\Phi \co \Delta C \to \Delta C$ (see Section~\ref{connes loday dual}).
The second main result of this paper is the following.

\begin{Theorem}\label{CYCINCOBRT}
	The cocyclic~$\kk$-modules~$\emph{\RT}_{\mc B}\circ X^\bullet$ and~$\coend^\bullet \circ \Phi$ are isomorphic.
	The cyclic~$\kk$-modules~$\emph{\RT}_{\mc B}\circ X_\bullet$ and~$\coend_\bullet \circ \Phi^\opp$ are isomorphic.
\end{Theorem}
We note that the $n$-th Hochschild cohomology of $\coend^\bullet \circ \Phi$ and $\coend^\bullet$ are equal. Indeed, the Hochschild differentials of the associated cochain complexes are equal. The cohomolo\-gical form of the Connes' long exact sequence~\cite[Proposition~9.6.11]{weibel} then implies that this is also the case for their cyclic cohomology.
A proof of Theorem~\ref{CYCINCOBRT} is provided in Sections~\ref{sect-modul-coend}--\ref{RTtildeS}.

The cyclic duality~$L\co \Delta C^{\opp}\to \Delta C$ and reindexing involution automorphism $\Phi$ transform the cocyclic object~$X^\bullet$ in~$\textbf{3}\Cob_0$ into a cyclic object~$X^\bullet \circ \Phi \circ L$ in~$\textbf{3}\Cob_0$.
Similarly, the functors~$L^\opp \co \Delta C\to \Delta C^{\opp}$ and~$\Phi^\opp$ transform the cyclic object~$X_\bullet$ in~$\textbf{3}\Cob_0$ into a cocyclic object~$X_\bullet \circ \Phi^{\opp} \circ L^\opp$ in~$\textbf{3}\Cob_0$.
By Theorem~\ref{CYCINCOBRT} and the fact that $\Phi$ is involutive, we obtain the following.

\begin{Corollary}\label{cycincobmajidsolotar}
The cyclic~$\kk$-modules $\emph{\RT}_{\mc B}\circ X^\bullet \circ \Phi \circ L$ and $\coend^\bullet \circ L$ are isomorphic.
The cocyclic~$\kk$-modules $\emph{\RT}_{\mc B}\circ X_\bullet \circ \Phi^{\opp} \circ L^\opp$ and $\coend_\bullet \circ L^\opp$ are isomorphic.
\end{Corollary}

Recall that the cyclic~$\kk$-module $\coend^\bullet \circ L$ and the cocyclic~$\kk$-module $\coend_\bullet \circ L^\opp$ (see Section~\ref{passing}) are cyclic duals of $\coend^\bullet$ and the cyclic $\kk$-module~$\coend_\bullet$.

Another fundamental construction of a~$3$-dimensional TQFT is the Turaev--Viro TQFT $\TV_{\mc C} \co \textbf{3}\Cob_0\to \Mod_\kk$ associated to a spherical fusion~$\kk$-category~$\mc C$ with invertible dimension (for details, see~\cite{moncatstft}). Moreover, in the case when~$\mc C$ is additive and~$\kk$ is an algebraically closed field, the center~$\Zen{\mc C}$ of~$\mc C$ is an anomaly free modular category (see~\cite[Theorems~5.3 and~5.4]{moncatstft}).
In this case, according to~\cite[Theorem~17.1]{moncatstft}, the TQFTs~$\RT_{\Zen{\mc C}}$ and~$\TV_{\mc C}$ are isomorphic.
Denote by~$\coendzc$ the coend of~$\Zen{\mc C}$.
These results and Theorem~\ref{CYCINCOBRT} imply the following corollary.

\begin{Corollary}
	The cocyclic $\kk$-modules $\emph{\TV}_{\mc C} \circ X^\bullet$ and $\coendzc^\bullet \circ \Phi$ are isomorphic. The cyclic~$\kk$-modules $\emph{\TV}_{\mc C} \circ X_\bullet$ and~$\coendzc_\bullet \circ \Phi^{\opp}$ are isomorphic.
\end{Corollary}
By using cyclic duality $L$ and reindexing involution $\Phi$, we obtain the following.

\begin{Corollary}
	The cyclic $\kk$-modules $\emph{\TV}_{\mc C} \circ X^\bullet \circ \Phi \circ L$ and $\coendzc^\bullet \circ L$ are isomorphic. The cocyclic $\kk$-modules $\emph{\TV}_{\mc C} \circ X_\bullet \circ \Phi^{\opp} \circ L^\opp$ and $\coendzc_\bullet \circ L^\opp$ are isomorphic.
\end{Corollary}

To show the claim from Theorem~\ref{CYCINCOBRT}, it follows from Section~\ref{generalsurf} that it suffices to compute~$\RT_\mc B\circ Y^\bullet$ and~$\RT_\mc B\circ Y_\bullet$. In Section~\ref{sect-modul-coend}, we give some algebraic preliminaries on coend of a~modular category. Next, in Section~\ref{sect-calc-RT-coend}, we describe the Reshetikhin--Turaev TQFT~$\RT_{\mc B}$ via the coend~$\coend$ of $\mc B$. Then, in Section~\ref{RTS}, we compute the cocyclic~$\kk$-module~$\RT_{\mc B}\circ Y^\bullet$.
In~Section~\ref{sec: finaltouch}, we prove that the latter is isomorphic to~$\coend^\bullet \circ \Phi$, as stated in Theorem~\ref{CYCINCOBRT}.
Finally, in Section~\ref{RTtildeS}, we sketch the computation of~$\RT_{\mc B}\circ Y_\bullet$ and the proof of the fact that it is isomorphic to the cyclic~$\kk$-module~$\coend_\bullet \circ \Phi^{\opp}$.

Remember that~$\coend$ is a Hopf algebra in~$\mc B$ endowed with a Hopf pairing~$\omega \co \coend \otimes \coend \to \uu$.
Here, we denote by~$m$, $u$, $\Delta$, $\varepsilon$, and $S$ multiplication, unit, comultiplication, counit, and antipode of~$\coend$, respectively. In what follows, we will often drop the notation of~$\coend$ while using graphical calculus.

\subsection{Modularity and pairing of a coend}\label{sect-modul-coend}

In this section, we provide some algebraic preliminaries needed for computation of the Reshe\-ti\-khin--Turaev TQFT~$\RT_{\mc B}$ via the coend of anomaly free modular category.
In the following lemma, we compute the inverse, under some conditions, of the pairing of the coend. Note that the statement of Lemma \ref{intcointexchange}\,$(a)$ is only a particular case of \cite[Theorem~5]{kerler}.
Also, the statement and the proof of Lemma \ref{intcointexchange}\,$(b)$ is similar to \cite[Lemma 6.2]{moncatstft}.
\begin{Lemma} \label{intcointexchange}
	Let~$\mc B$ be a rib\-bon $\kk$-ca\-te\-go\-ry with a co\-end~$\coend$ and suppose that the canonical pair\-ing~$\omega\co \coend\tens \coend \to \uu$ associated to the coend is non-degenerate.
	\begin{itemize}\itemsep=0pt
		\item[$(a)$] If~$\Lambda \co \uu \to \coend$ is a right integral of a coend, then~$\omega(\Lambda \tens \idrm_\coend)\co \coend\to \uu$ is a left cointegral of~$\coend$.
		\item[$(b)$] Let~$\Lambda$ be a right integral of the coend~$\coend$ of~$\mc B$.
		Suppose that the element $\omega(\Lambda \tens \Lambda)$ is invertible.
		The inverse of the pairing~$\omega$ is given by the morphism~$\Omega\co \uu \to \coend\tens \coend$, which is computed by
		\[
		\Omega=\left(\,		
 \put(3,-11.7){\small $\Lambda$}
 \put(23,-11.7){\small $\Lambda$}
 \put(13.5,11.5){\small $\omega$}
 \raisebox{-5mm}{\includegraphics[scale=0.85]{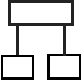}}
		\,\right)^{-1}
 \put(19.7,-37.5){\small $\Lambda$}
 \put(19.7,-11.5){\small $\Lambda$}
 \put(33,11.7){\small $\omega$}
 \raisebox{-14mm}{\includegraphics[scale=0.85]{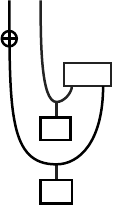}}\,.
	\]
	\end{itemize}
\end{Lemma}

\begin{proof} We prove the above statements graphically.
	\begin{itemize}\itemsep=0pt
		\item[$(a)$] By the property of the bialgebra pairing and the fact that $\Lambda$ is a right integral for~$C$,\vspace{-1mm} 
		\[
 \put(38.5,-6.5){\small $\Lambda$}
 \put(10,27){\small $\omega$}
 \put(48,27){\small $\omega$}
 \raisebox{-9.5mm}{\includegraphics[scale=0.85]{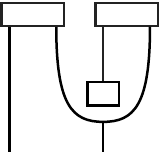}}\, =
 \put(5.5,-12.5){\small $\Lambda$}
 \put(28,11){\small $\omega$}
 \put(33,27){\small $\omega$}
 \,\raisebox{-9.5mm}{\includegraphics[scale=0.85]{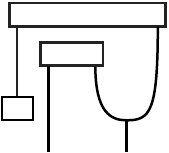}}\, =
 \put(5.5,-19){\small $\Lambda$}
 \put(41,27){\small $\omega$}
 \,\raisebox{-9.5mm}{\includegraphics[scale=0.85]{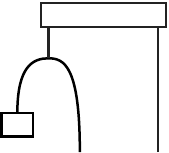}}\, =
 \put(5.5,4){\small $\Lambda$}
 \put(28,27){\small $\omega$}
 \,\raisebox{-9.5mm}{\includegraphics[scale=0.85]{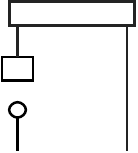}}\, .
\]
Since $\omega$ is non-degenerate, the conclusion follows.
		\item[$(b)$] We have\vspace{-1mm}
		\begingroup
		\allowdisplaybreaks
		\begin{gather*}
 \put(31,34){\small $\omega$}
 \put(34.5,14){\small $\omega$}
 \put(20,-12){\small $\Lambda$}
 \put(20,-34.5){\small $\Lambda$}
 \raisebox{-15mm}{\includegraphics[scale=0.85]{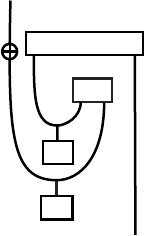}}\, \overset{(i)}{=}
 \put(38,34){\small $\omega$}
 \put(21,-8.5){\small $\Lambda$}
 \put(21,-34){\small $\Lambda$}
 \,\raisebox{-15mm}{\includegraphics[scale=0.85]{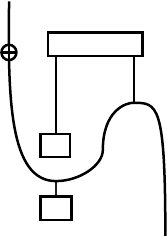}}\, \overset{(ii)}{=}
 \put(34.5,17.5){\small $\omega$}
 \put(19,-9){\small $\Lambda$}
 \put(19,-35){\small $\Lambda$}
 \,\raisebox{-15mm}{\includegraphics[scale=0.85]{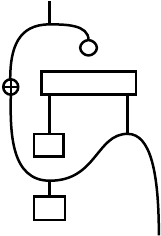}}\, \overset{(iii)}{=}
 \put(55,40.5){\small $\omega$}
 \put(46.5,17){\small $\Lambda$}
 \put(18,-34.5){\small $\Lambda$}
 \,\raisebox{-15mm}{\includegraphics[scale=0.85]{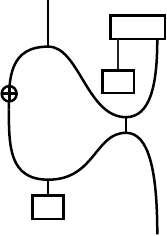}}\, \\
		\qquad\overset{(iv)}{=}
 \put(55.5,38.5){\small $\omega$}
 \put(48,15){\small $\Lambda$}
 \put(21,-34.5){\small $\Lambda$}
 \,\raisebox{-15mm}{\includegraphics[scale=0.85]{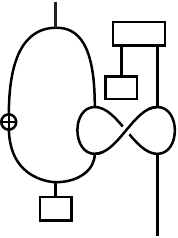}}\, \overset{(v)}{=}
 \put(58,38.5){\small $\omega$}
 \put(50.5,15){\small $\Lambda$}
 \put(23.5,-34.5){\small $\Lambda$}
 \,\raisebox{-15mm}{\includegraphics[scale=0.85]{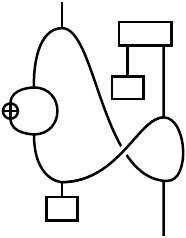}}\, \overset{(vi)}{=}
 \put(34,38.5){\small $\omega$}
 \put(26,14.5){\small $\Lambda$}
 \put(9,-35){\small $\Lambda$}
 \,\raisebox{-15mm}{\includegraphics[scale=0.85]{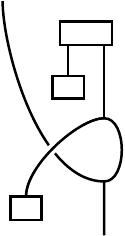}}\, \overset{(vii)}{=}		
 \put(34,14.5){\small $\omega$}
 \put(24,-8.5){\small $\Lambda$}
 \put(43.5,-8.5){\small $\Lambda$}
 \,\raisebox{-15mm}{\includegraphics[scale=0.85]{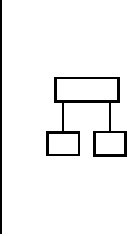}}\,.
		\end{gather*}
		\endgroup and therefore the claim follows by invertibility of $\omega(\Lambda \tens \Lambda)$.
		Here $(i)$ follows by pro\-perty of bialgebra pairing, $(ii)$ by (co)unitality, $(iii)$ by the part $(a)$, $(iv)$ from bialgebra axioms, $(v)$ by (co)associativity, $(vi)$ by the antipode axiom and (co)unitality,~$(vii)$ by the naturality of the braiding, the fact that $\Lambda$ is a right integral for $\coend$ and counitality.\hfill \qed
	\end{itemize}
\renewcommand{\qed}{}
\end{proof}

\begin{Corollary} \label{leminvpetitcalcul}
	Let $\mc B$ be a ribbon $\kk$-category with a coend $\coend$ and $\Lambda \co \uu \to \coend$ a right integral for the coend $\coend$.
	If the pairing $\omega\co \coend \tens \coend \to \uu$ is non-degenerate, then
	\begin{equation}\label{invpetitcalcul}
 \put(32.5,35){\small $\omega$}
 \put(33,13){\small $\omega$}
 \put(20,-14.2){\small $\Lambda$}
 \put(20,-36){\small $\Lambda$}
 \raisebox{-18mm}{\includegraphics[scale=0.85]{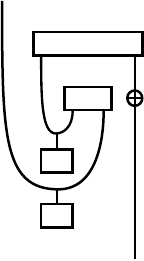}}\,	= \omega(\Lambda \tens \Lambda)
 \,\raisebox{-18mm}{\includegraphics[scale=0.85]{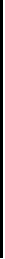}}\,	.
	\end{equation}
\end{Corollary}

\begin{proof}\samepage
	By Lemma~\ref{intcointexchange}\,$(b)$ and invertibility of the antipode, we have
	\[
	\omega(\Lambda \tens \Lambda) \id_\coend \overset{(i)}{=}
 \put(5,-11.8){\small $\Lambda$}
 \put(25,-11.8){\small $\Lambda$}
 \put(15.5,11.5){\small $\omega$}
 \put(55,-0.5){\small $\id_{\coend}$}
 \,\raisebox{-17.3mm}{\includegraphics[scale=0.85]{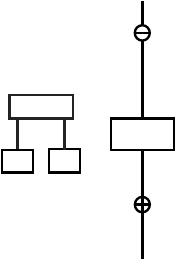}}	\,\overset{(ii)}{=}
 \put(37.5,37){\small $\omega$}
 \put(38,14.5){\small $\omega$}
 \put(25,-12.2){\small $\Lambda$}
 \put(25,-34.2){\small $\Lambda$}
 \,\raisebox{-17.3mm}{\includegraphics[scale=0.85]{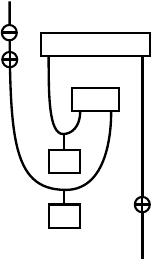}}	\,\overset{(iii)}{=}
 \put(34.5,37){\small $\omega$}
 \put(35,14.5){\small $\omega$}
 \put(22,-12.2){\small $\Lambda$}
 \put(22,-34.2){\small $\Lambda$}
 \,\raisebox{-17.3mm}{\includegraphics[scale=0.85]{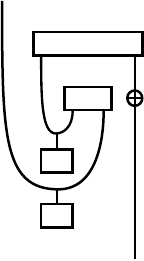}}\,	.
\]
Here $(i)$ and $(iii)$ follow from invertibility of the antipode, $(ii)$ follows by Lemma \ref{intcointexchange}\,$(b)$.
\end{proof}

\pagebreak

If~$\mc B$ is an additive ribbon fusion~$\kk$-category with a coend~$\coend$, then, according to~\cite[Theorem~6.6]{moncatstft}, the category~$\mc B$ is modular (in the sense of Section~\ref{modular}) if and only if the canonical pairing~$\omega\co \coend \tens \coend \to \uu$ associated to~$\coend$ is non-degenerate.
Moreover, the coend~$\coend$ is given by~\smash{$\coend= \bigoplus_{i\in I} i^*\tens i$}, where $I$ is a representative set of simple objects of $\mc B$. For $i\in I$, we denote the projection associated with the direct sum decomposition by $p_i \co \coend\to i^*\tens i$. However, we drop inclusions $i^*\tens i \to \coend$ in our notation. By \cite[Theorem~6.4]{moncatstft}, any integral of $\coend$ is a scalar multiple of the universal integral
\begin{equation} \label{univint}
\Lambda =\sum_{i\in I} \dim(i)\widetilde{\text{coev}}_i \co\ \uu \to \coend.
\end{equation}
Similarly, one can show that any cointegral of $\coend$ is a scalar multiple of the universal cointegral
\begin{equation} \label{univcoint}
\lambda =\text{ev}_{\uu}p_{\uu} \co\ \coend \to \uu.
\end{equation}
Universal (co)integrals $\lambda$ and $\Lambda$ satisfy $\lambda \Lambda= \id_\uu$.

\begin{Remark}\label{omegaalphaalpha} Let $\mc B$ be an additive ribbon fusion~$\kk$-category such that the canonical pairing~$\omega$ associated to the coend $\coend$ is non-degenerate.
	Recall the universal integral~$\Lambda$ and the universal cointegral~$\lambda$ of the coend~$\coend$, defined in equations~\eqref{univint} and~\eqref{univcoint}, respectively. Let us calculate~$\omega(\Lambda \tens \Lambda)$.
	By Lemma~\ref{intcointexchange}\,$(a)$, $\omega(\Lambda \tens \id_\coend)$ is a left cointegral of~$\coend$. By universality of~$\lambda$,~$\omega(\Lambda \tens \id_\coend)=k\lambda$, for some~$k \in \kk\cong \End(\uu)$.
	This further im\-plies that
	\[
	\omega(\Lambda \tens \Lambda)= k\lambda \Lambda =k\id_{\uu}.
\]
	Now, re\-mark also that
	\[
	\omega(\Lambda \tens u)= \varepsilon \Lambda = \dim(\mc B).
\]
	These two properties together with the fact that~$p_\uu i_\uu = \id_{\uu^* \tens \uu}$, with definition of unit~$u$ of the coend $\coend$, and definition of~$\lambda$, give that
	\begin{align*}
\omega(\Lambda \tens \Lambda)=&k\id_{\uu}=k(\id_\uu \tens \evl_\uu )(\coevl_\uu \tens \id_\uu) = k(\id_\uu \tens \evl_\uu p_\uu )(\id_\uu \tens i_\uu )(\coevl_\uu \tens \id_\uu) \\
= & k(\id_\uu \tens \lambda) u = k \lambda u = \omega(\Lambda \tens u)= \dim(\mc B).
	\end{align*}

\end{Remark}

\subsection{The Reshetikhin--Turaev TQFT via coends}
\label{sect-calc-RT-coend}
In~\cite{turaevqinvariants}, Turaev associates to any modular category~${\mc B}$ a 3-dimensional TQFT~$\RT_{\mc B}$. There, a~precise definition of a $3$-dimensional TQFT involves Lagrangian spaces in ho\-mo\-lo\-gy of surfaces and~$p_1$-structures in cobordisms. However, if the modular category~${\mc B}$ is anomaly free (see Section~\ref{modular}), then the TQFT~$\RT_{\mc B}$ does not depend on this additional data and is a genuine symmetric monoidal functor~$\RT_{\mc B} \co \textbf{3}\Cob_0 \to \Mod_\kk$.

Let ${\mc B}$ be an anomaly free modular $\kk$-category. Recall from Section~\ref{modular} that the scalar $\Delta=\Delta_+=\Delta_-$ is invertible and satisfies~$\Delta^2=\dim(\mc B)$.
By Section~\ref{specialribbon}, any special ribbon~$(g,h)$-graph~$\Gamma$ represents a~$3$-cobordism~$M_\Gamma \co S_g \to S_{h}$. Our goal is to compute the~$\kk$-linear homomorphism
\[
\RT_{\mc B}(M_\Gamma) \co\ \RT_{\mc B}(S_g) \to \RT_{\mc B}(S_h)
\]
in terms of the coend~$\coend$ of~${\mc B}$ (which always exists, see Section~\ref{sect-modul-coend}). First, it follows from the definition of~$\RT_{\mc B}$ and the computation of the coend~$\coend$ in terms of a representative set~$I$ of simple objects of~${\mc B}$ that
\[
\RT_{\mc B}(S_g)=\Hom_{\mc B}\bigl(\uu,\coend^{\otimes g}\bigr) \qquad \text{and} \qquad \RT_{\mc B}(S_h)=\Hom_{\mc B}\bigl(\uu,\coend^{\otimes h}\bigr).	
\]
Next, the formula $(2.3)\,(a)$ from \cite[Section~IV.2.3]{turaevqinvariants}, which computes $\RT_{\mc B}(M_\Gamma)$, rewrites in our setting as
\begin{equation}
\label{rtformula}
\RT_{\mc B}(M_\Gamma)=\Delta^{-n-h} \,
\Hom_{\mc B}(\uu, \lvert \Gamma \rvert),
\end{equation}
where~$\lvert \Gamma \rvert \co \coend^{\otimes g} \to \coend^{\otimes h}$ is a morphism in~$\mc B$ defined as follows.
By pulling down some part of each circle component of~$\Gamma$ and of each arc connecting the outputs of~$\Gamma$, we obtain that the ribbon graph~$\Gamma$ is isotopic to
\[
\put(52.5,33){\small $\cdots$}
\put(52.5,-36){\small $\cdots$}
\put(196.5,-23){\small $\cdots$}
\put(341,-23){\small $\cdots$}
\put(1.5,53){\small $1$}
\put(1.5,-58){\small $1$}
\put(38,53){\small $2$}
\put(38,-58){\small $2$}
\put(59,53){\small $2h-1$}
\put(59,-58){\small $2h-1$}
\put(107.5,53){\small $2h$}
\put(107.5,-58){\small $2h$}
\put(182,-40){\small $L_1$}
\put(254,-40){\small $L_n$}
\put(327,-40){\small $C_1$}
\put(399,-40){\small $C_{h}$}
\put(190,-3.5){\small $\widetilde{\Gamma}$}
\raisebox{-23mm}{\includegraphics[scale=0.8]{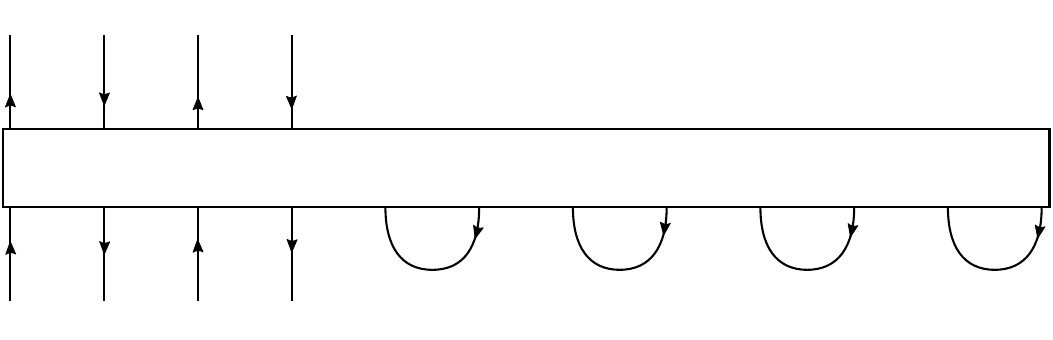}}	
\]
where the cups~$L_1, \dots, L_n$ correspond to the circle components of~$\Gamma$ and the cups~$C_1, \dots, C_h$ correspond to the upper arcs of~$\Gamma$.
Here,~\smash{$\widetilde{\Gamma}$} is a ribbon graph with~$(g+n+2h)$ arcs, $(2g+2n+2h)$ inputs,~$2h$ outputs, no coupons, no circle components, and such that:
\begin{itemize}\itemsep=0pt
	\item for all $1 \le i \le g+n$, an arc $a_i$ connects the $(2i-1)$-th input to the $(2i)$-th input of~$\widetilde{\Gamma}$,
	\item for all $1 \le j \le h$, an arc $u_j$ connects the $(2g+2n+2j-1)$-th input to the $(2j-1)$-th output of $\widetilde{\Gamma}$,
	and an arc $v_j$ connects the $(2j)$-th output to the $(2g+2n+2j)$-th input of~$\widetilde{\Gamma}$.
\end{itemize}
The ribbon graph ${\Gamma}$ is called \textit{closure} of the ribbon graph $\widetilde{\Gamma}$.
Coloring the arc~$a_i$ by an object~$X_i$ of~$\mc B$ and coloring both the arcs~$u_j$, $v_j$ by an object~$Y_j$ of~$\mc B$, we obtain a~$\mc B$-colored ribbon graph representing a morphism~$\phi_{X_1, \dots, X_{g+n}, Y_1, \dots, Y_h}$.
Let $i=\{i_X \co X^*\otimes X \to \coend\}_{X \in \Ob{\mc B}}$ be the universal dinatural transformation associated to the coend~$\coend$.
Then the family of morphisms
\[
\bigl(i_{Y_1} \otimes \cdots \otimes i_{Y_h}\bigr) \circ \phi_{X_1, \dots, X_{g+n}, Y_1, \dots, Y_h}
\]
from $X_1^* \otimes X_1 \otimes \cdots \otimes X_{g+n}^* \otimes X_{g+n} \otimes
Y_1^* \otimes Y_1 \otimes \cdots \otimes Y_h^* \otimes Y_h$ to $\coend^{\otimes h}$ is dinatural in each variable and so, by Lemma~\ref{CoendFubini}, it factorizes as
\[
\phi_{\Gamma} \circ \bigl(i_{X_1} \otimes \cdots \otimes i_{X_{g+n}} \otimes i_{Y_1} \otimes \cdots \otimes i_{Y_h}\bigr)
\]
for a unique morphism $\phi_{\Gamma} \co \coend^{\otimes g+n+h} \to \coend^{\otimes h}$. Then
\[
\lvert \Gamma \rvert=\phi_{\Gamma} \circ\bigl(\id_{\coend^{\otimes g}} \otimes \Lambda^{\otimes (n+h)}\bigr) \co\ \coend^{\otimes g} \to \coend^{\otimes h},
\]
where $\Lambda$ is the universal integral defined in equation~\eqref{univint}.
It follows from the fact that~$\Lambda$ is a~right integral for~$\coend$ (see Section~\ref{sect-modul-coend}) that the morphism~$\lvert \Gamma \rvert$ is an isotopy invariant of~$\Gamma$. This invariant is multiplicative,
\begin{equation*}
\lvert \Gamma \sqcup \Gamma'\rvert=\lvert \Gamma \rvert \otimes \lvert \Gamma' \rvert
\end{equation*}
for all special ribbon graphs, where~$\Gamma \sqcup \Gamma'$ is obtained by concatenating $\Gamma'$ to the right of~$\Gamma$.
\begin{Remark} If we normalize the above used $\lambda$ and $\Lambda$ and instead use $\lambda'=\Delta \lambda$ and $\Lambda'= \Delta^{-1} \Lambda$, then we still have $\lambda'\Lambda'= \id_\uu$, but the factor $\Delta^{-n-h}$ in formula \eqref{rtformula} disappears, which makes it and some of our next computations look simpler. This is a matter of convention.
\end{Remark}

\subsection[Computation of RT\_B circ Y\^{}\{bullet\}]{Computation of $\boldsymbol{\RT_{\mc B}\circ Y^\bullet}$} \label{RTS}
Recall that $\mc B$ denotes an anomaly free modular category.
In the following lemma, we calculate the isotopy invariant $\lvert \cdot\rvert$ (see Section~\ref{sect-calc-RT-coend}) for the particular special graphs.
\begin{Lemma} \label{simpletangles} Let $T_1$, $T_2$, and $T_3$ be the following special ribbon graphs:
	\[T_1 =\,
	\rsdraw{0.45}{0.70}{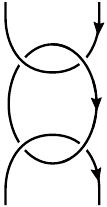}\;, \qquad T_2 =\,
	\rsdraw{-0.13}{0.70}{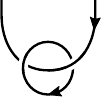}\;, \qquad T_3 = \, \rsdraw{0.45}{0.70}{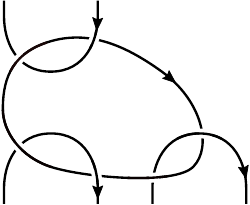}\;.\]
We have
\[
(a) \ \lvert T_1 \rvert= \dim(\mc B) \idrm_{\coend},\qquad
(b) \ \lvert T_2 \rvert= \dim(\mc B) u, \qquad
(c) \ \lvert T_3 \rvert= \dim(\mc B) m.
\]
\end{Lemma}

\begin{proof}
$(a)$ A ribbon graph whose closure is isotopic to $T_{1}$ is
		\[\widetilde{T_1} = \,
		\rsdraw{0.45}{0.73}{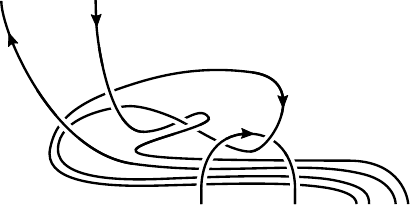}\;.
\]
		For all objects $X$, $Y$, and $Z$ in $\mc B$, we have
		\begingroup
		\allowdisplaybreaks
		\begin{align*}
 \put(39,14){\small $Z$}
 \put(107,-21){\small $Y$}
 \put(102,-62){\small $X$}
		\raisebox{-23mm}{\includegraphics[scale=0.75]{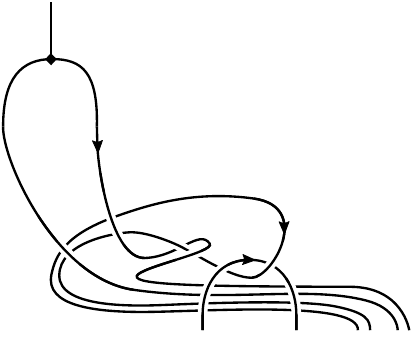}}\,&=	
 \put(40,14){\small $Z$}
 \put(75,-10){\small $Y$}
 \put(100,-62){\small $X$}
 \put(77.5,31.5){\small $\omega$}
		\,\raisebox{-23mm}{\includegraphics[scale=0.75]{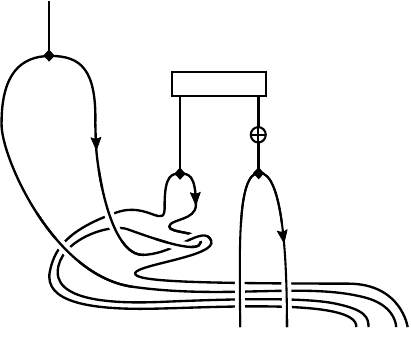}}\\
 &=
 \put(51,-37){\small $Z$}
 \put(11,-43){\small $Y$}
 \put(87,-68){\small $X$}
 \put(62,26.5){\small $\omega$}
 \put(60,6.5){\small $\omega$}
	 \, \raisebox{-25mm}{\includegraphics[scale=0.75]{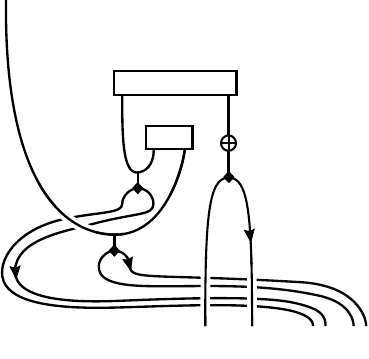}}\, .
		\end{align*}
	\endgroup
		Hence, by definition of~$\lvert \cdot \rvert$ given in Section~\ref{sect-calc-RT-coend}, Corollary~\ref{leminvpetitcalcul}, and Remark~\ref{omegaalphaalpha}, we have
		\[
		\lvert T_1 \rvert =
 \put(34.5,37){\small $\omega$}
 \put(35,14.5){\small $\omega$}
 \put(22,-12.2){\small $\Lambda$}
 \put(22,-34.2){\small $\Lambda$}
	 \, \raisebox{-17.3mm}{\includegraphics[scale=0.85]{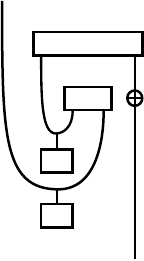}}\,
 = \omega(\Lambda \tens \Lambda)\id_{\coend} = \dim(\mc B) \id_{\coend}.
\]

$(b)$ A ribbon graph whose closure is isotopic to $T_2$ is \[
		\widetilde{T_2} = \, \rsdraw{0.45}{0.75}{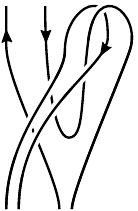}\;.
\]
		For all objects $X$, $Y$ in $\mc B$, we have
		\[
 \put(3,-59){\small $X$}
 \put(23,-59){\small $Y$}
	 \raisebox{-21mm}{\includegraphics[scale=0.75]{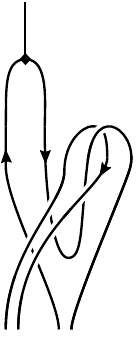}}\,=		
 \put(5,-59){\small $X$}
 \put(26,-59){\small $Y$}
 \put(35,33){\small $\omega$}
	 \,\raisebox{-21mm}{\includegraphics[scale=0.75]{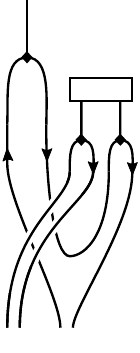}}\,=		
 \put(4.5,-59){\small $X$}
 \put(29,-59){\small $Y$}
 \put(34.5,33){\small $\omega$}
	 \,\raisebox{-21mm}{\includegraphics[scale=0.75]{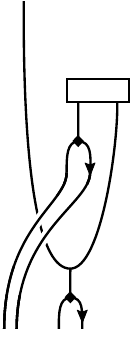}}\, .	
\]
		Hence, by definition of~$\lvert \cdot \rvert$ given in Section~\ref{sect-calc-RT-coend}, by Lemma~\ref{intcointexchange}\,$(a)$, and Remark~\ref{omegaalphaalpha}, we have
		\[
		\lvert T_2 \rvert =
 \put(11.5,-30){\small $\Lambda$}
 \put(11.5,-1.5){\small $\Lambda$}
 \put(18.5,22){\small $\omega$}
	 \,\raisebox{-11.5mm}{\includegraphics[scale=0.85]{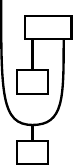}}\, =		
 \put(18,-31){\small $\Lambda$}
 \put(5.5,-15){\small $\Lambda$}
 \put(14,8.5){\small $\omega$}
	 \,\raisebox{-11.5mm}{\includegraphics[scale=0.85]{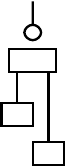}}\, = \omega(\Lambda \tens \Lambda)u = \dim(\mc B) u.		
\]

$(c)$ A ribbon graph whose closure is isotopic to $T_3$ is
		\[\widetilde{T_3}=\,
		\rsdraw{0.45}{0.75}{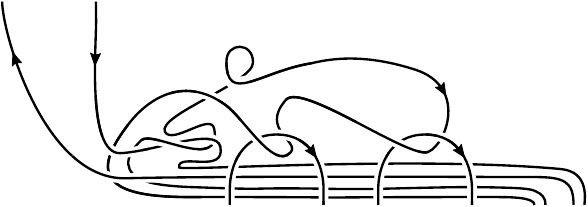}\;.\]
		For all objects $X$, $Y$, and $Z$ in $\mc B$, we have
		\begin{gather*}
 \put(38,0){\small $W$}
 \put(160,-20){\small $Z$}
 \put(162,-68.5){\small $Y$}
 \put(109,-68.5){\small $X$}
	 \raisebox{-25mm}{\includegraphics[scale=0.728]{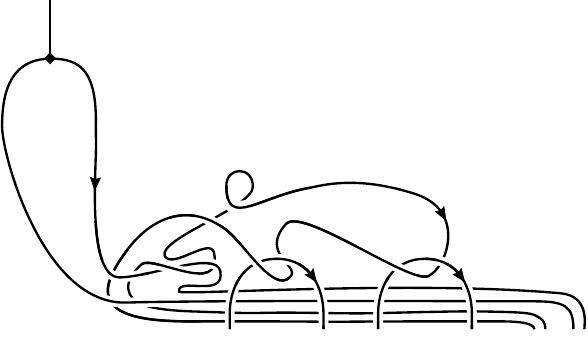}} 	
=
 \put(39,-4){\small $W$}
 \put(80,-14){\small $Z$}
 \put(164,-68.5){\small $Y$}
 \put(111,-68.5){\small $X$}
 \put(88,32.5){\small $\omega$}
 \put(106,46.5){\small $\omega$}
 \,\raisebox{-25mm}{\includegraphics[scale=0.728]{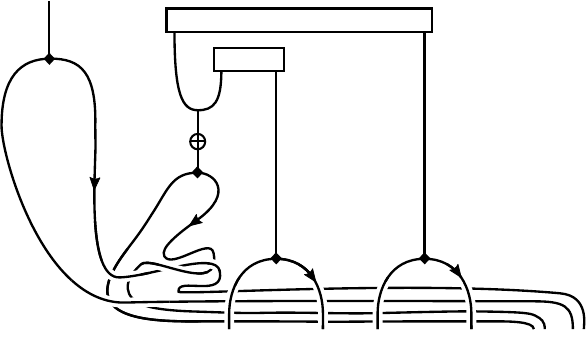}}\\				
 \qquad\quad=
 \put(180,-68.5){\small $W$}
 \put(165,-68.5){\small $Z$}
 \put(138,-68.5){\small $Y$}
 \put(84,-68.5){\small $X$}
 \put(56,4.5){\small $\omega$}
 \put(57.5,35.5){\small $\omega$}
 \put(75,49.5){\small $\omega$}
	 \,\raisebox{-25mm}{\includegraphics[scale=0.75]{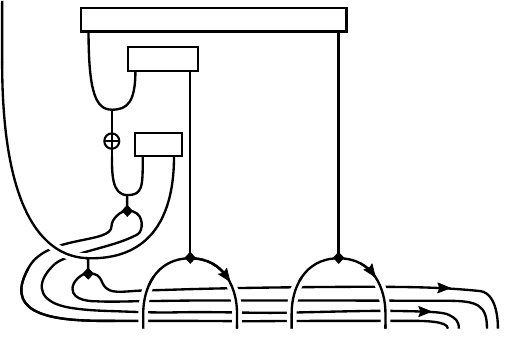}}\,.
		\end{gather*}
Hence, by definition of~$\lvert \cdot \rvert$ given in Section~\ref{sect-calc-RT-coend}, axioms of a Hopf pairing, the equation~\eqref{omegaS}, Corollary~\ref{leminvpetitcalcul}, and Remark~\ref{omegaalphaalpha}, we have
		\[
 \lvert T_3 \rvert =
 \put(42,49.5){\small $\omega$}
 \put(40.5,35.5){\small $\omega$}
 \put(38.3,5){\small $\omega$}
 \put(26.5,-24.8){\small $\Lambda$}
 \put(26.5,-47.2){\small $\Lambda$}
	 \,\raisebox{-21mm}{\includegraphics[scale=0.75]{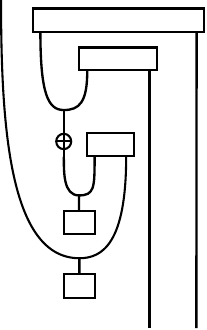}}\,=
 \put(30,49.5){\small $\omega$}
 \put(29.5,25.5){\small $\omega$}
 \put(19,-4.5){\small $\Lambda$}
 \put(19,-27.5){\small $\Lambda$}
	 \,\raisebox{-21mm}{\includegraphics[scale=0.75]{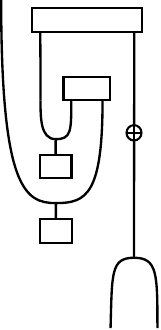}}\,= \omega(\Lambda \tens \Lambda)m = \dim(\mc B)m.
 \tag*{\qed}
 \]
\renewcommand{\qed}{}
\end{proof}
The previously proved Lemma \ref{simpletangles} implies the following result.

\begin{Lemma} \label{calculinvisotopie} We have
\begin{itemize}\itemsep=0pt
		\item[$(a)$] $\bigl\lvert G_{Y^\bullet (\delta_i^n)} \bigr\rvert=
\begin{cases}
\dim(\mc B)^{n+1}
\put(35,0){\small $\cdots$}
\put(20,-22.5){\small $0$}
\put(47,-22.5){\small $n-1$}
\,\raisebox{-9.3mm}{\includegraphics[scale=0.75]{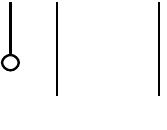}}\, , &	\mbox{if } i=0, \\[4.5mm]
\dim(\mc B)^{n+1}
\put(7.2,0){\small $\cdots$}
\put(61.1,0){\small $\cdots$}
\put(1,-22.5){\small $0$}
\put(13,-22.5){\small $i-1$}
\put(54.5,-22.5){\small $i$}
\put(65,-22.5){\small $n-1$}
\,\raisebox{-9.3mm}{\includegraphics[scale=0.75]{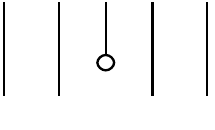}}\, , & \mbox{if } 1\le i \le n-1, \\[4.5mm]
\dim(\mc B)^{n+1}
\put(16.5,0){\small $\cdots$}
\put(0.8,-22.5){\small $0$}
\put(28,-22.5){\small $n-1$}
\,\raisebox{-9.3mm}{\includegraphics[scale=0.75]{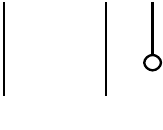}}\, , & \mbox{if } i=n,
\end{cases}$

\item[$(b)$]
$
\bigl\lvert G_{Y^\bullet (\sigma_j^n)} \bigr\rvert=\dim(\mc B)^{n+1} 
\put(7.3,0){\small $\cdots$}
\put(72.5,0){\small $\cdots$}
\put(1,-23.5){\small $0$}
\put(33,-23.5){\small $j$}
\put(43,-23.5){\small $j+1$}
\put(76,-23.5){\small $n+1$}
\,\raisebox{-9.8mm}{\includegraphics[scale=0.75]{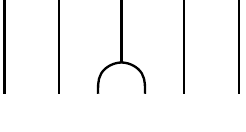}}\,,
$

\item[$(c)$] $\bigl\lvert G_{Y^\bullet (\tau_n)} \bigr\rvert =
\begin{cases}
\dim(\mc B)\idrm_{\coend} & \mbox{if } n=0, \\
\dim(\mc B)^{n+1}
\put(40,-15){\small $\cdots$}
\put(1.3,-24){\small $0$}
\put(18.1,-24){\small $1$}
\put(59.6,-24){\small $n-1$}
\,\raisebox{-9.2mm}{\includegraphics[scale=0.75]{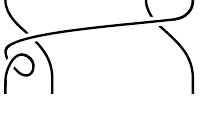}}	& \mbox{if } n\ge 1.	
\end{cases} $
	\end{itemize}
\end{Lemma}

\begin{proof}
$(a)$ Recall the special ribbon graphs $T_1$ and $T_2$ from Lemma \ref{simpletangles}.
		Let $n\in \N^*$ and $1 \le i \le n-1$. By Lemma \ref{simpletangles} and multiplicativity of $\lvert \cdot \rvert$,
		\begin{align*}
		\lvert G_{Y^\bullet (\delta_i^n)} \rvert&= \left(\dim(\mc B) \id_\coend \right)^{\tens i} \tens \dim(\mc B) u \tens \left(\dim(\mc B)\id_\coend \right)^{\tens n-i}= \hspace{0.65cm} \\
&=\dim(\mc B)^{n+1}
\put(7.2,0){\small $\cdots$}
\put(61.1,0){\small $\cdots$}
\put(1,-22.5){\small $0$}
\put(13,-22.5){\small $i-1$}
\put(54.5,-22.5){\small $i$}
\put(65,-22.5){\small $n-1$}
\,\raisebox{-9.3mm}{\includegraphics[scale=0.7]{page34-02}}\,.
		\end{align*}
The cases $i=0$ and $i=n$ are proven in a similar way.

$(b)$ This follows from parts $(a)$ and $(c)$ of~Lem\-ma~\ref{simpletangles} and multiplicativity of the isotopy invariant $\lvert \cdot \rvert$.		

$(c)$ Recall that~$Y^\bullet (\tau_0)$ is the identity morphism $\id_{S_1}$. Hence $\bigl\lvert G_{Y^\bullet (\tau_0)} \bigr\rvert= \dim({\mc B}) \id_{\coend}$, by Lemma~\ref{simpletangles}\,$(a)$. Let us show the statement in the case~$n=1$. The general case is verified similarly.
		The special ribbon graph~$G_{Y^\bullet (\tau_1)}$ depicts as
		\[
		\,\rsdraw{0.55}{0.75}{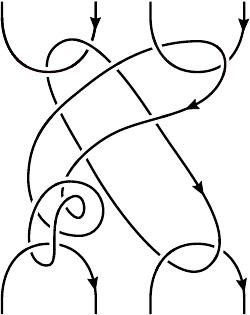} \;.
		\]
		A ribbon graph whose closure is isotopic to the special ribbon graph~$G_{Y^\bullet (\tau_1)}$ depicts as
		\[
		\widetilde{G_{Y^\bullet (\tau_1)}}= \, \rsdraw{0.55}{0.8}{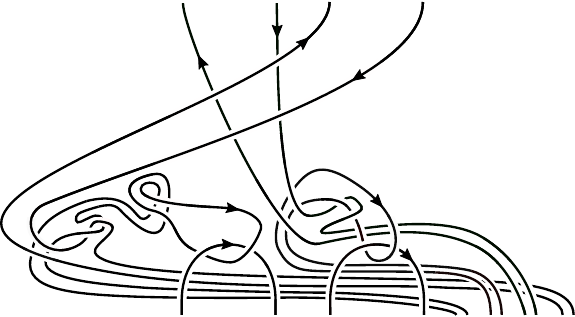} \;.
		\]
		
		For all objects $X$, $Y$, $Z$, $U$, $V$, and $W$ in $\mc B$, we have
		\begingroup
		\allowdisplaybreaks
		\begin{gather*}
 \put(103,38){\small $U$}
 \put(156,38){\small $V$}
 \put(95,-88){\small $X$}
 \put(151,-88){\small $Y$}
 \put(165,-88){\small $Z$}
 \put(178,-88){\small $W$}
 \raisebox{-32mm}{\includegraphics[scale=0.75]{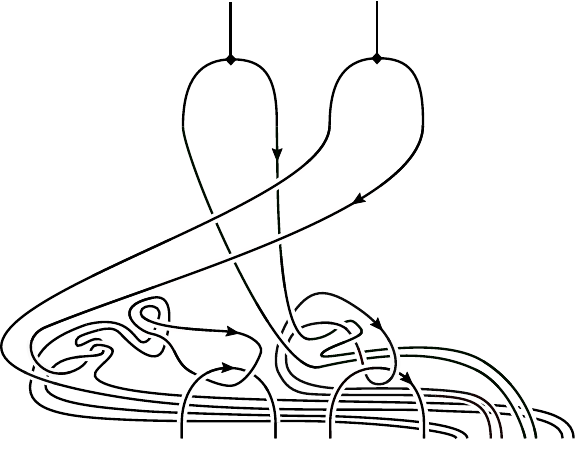}}\\
 \qquad\quad =
 \put(105,38){\small $U$}
 \put(158,38){\small $V$}
 \put(87,-88){\small $X$}
 \put(154,-88){\small $Y$}
 \put(167,-88){\small $Z$}
 \put(180,-88){\small $W$}
 \put(35,-8){\small $\omega$}
 \put(73,-33.5){\small $\omega$}
 \put(137.3,-10.5){\small $\omega$}
 \,\raisebox{-32mm}{\includegraphics[scale=0.75]{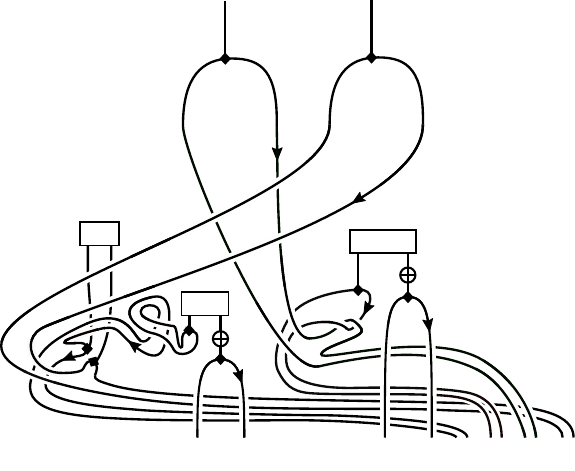}}\\
 \qquad\quad =
 \put(79,-88){\small $X$}
 \put(148,-88){\small $Y$}
 \put(160.5,-88){\small $Z$}
 \put(173.5,-88){\small $W$}
 \put(187.5,-88){\small $U$}
 \put(200,-88){\small $V$}
 \put(44,5){\small $\omega$}
 \put(35.5,-23){\small $\omega$}
 \put(123,-15){\small $\omega$}
 \put(123,2){\small $\omega$}
 \,\raisebox{-32mm}{\includegraphics[scale=0.75]{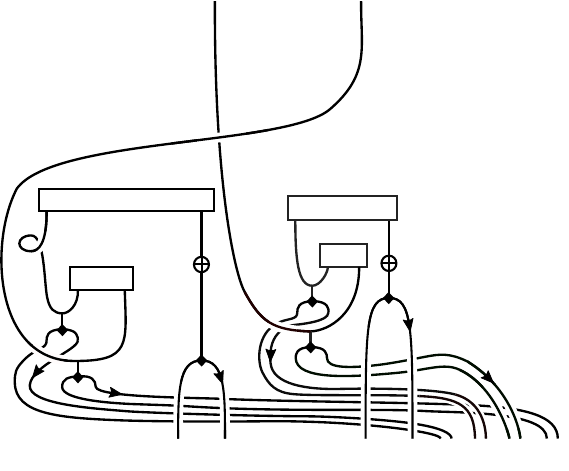}}
		\end{gather*}
		\endgroup
		Hence, by definition of~$\lvert \cdot \rvert$ given in Section~\ref{sect-calc-RT-coend}, the equations~\eqref{omegaS} and~\eqref{invol}, Corollary~\ref{leminvpetitcalcul} and Remark~\ref{omegaalphaalpha}, we have
		\begin{align*}
\lvert G_{Y^\bullet (\tau_1)} \rvert=&
 \put(44,-4.5){\small $\omega$}
 \put(35,-32.5){\small $\omega$}
 \put(20.8,-55.3){\small $\Lambda$}
 \put(26.5,-73.6){\small $\Lambda$}
 \put(123.5,-7){\small $\omega$}
 \put(122.5,-24.5){\small $\omega$}
 \put(110.7,-44.7){\small $\Lambda$}
 \put(110.2,-62.6){\small $\Lambda$}
 \,\raisebox{-31mm}{\includegraphics[scale=0.75]{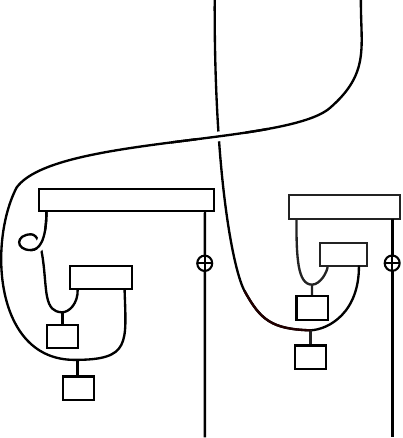}}\,=
 \put(44,-4.5){\small $\omega$}
 \put(35,-32.5){\small $\omega$}
 \put(20.8,-55.3){\small $\Lambda$}
 \put(26.5,-73.6){\small $\Lambda$}
 \put(123.5,-7){\small $\omega$}
 \put(122.5,-24.5){\small $\omega$}
 \put(110.7,-44.7){\small $\Lambda$}
 \put(110.2,-62.6){\small $\Lambda$}
 \,\raisebox{-31mm}{\includegraphics[scale=0.75]{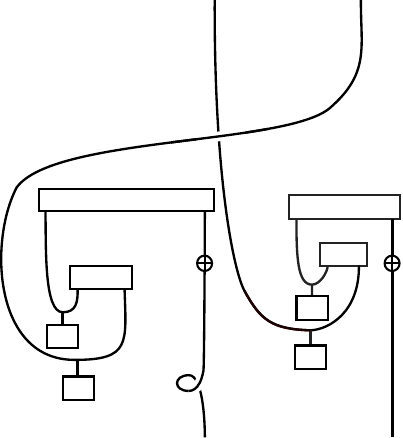}}\\
=& \left(\omega(\Lambda \tens \Lambda)\right)^{2} \, 
		\rsdraw{0.55}{0.75}{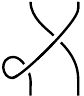} \; = \dim(\mc B)^{2} \, 
		\rsdraw{0.55}{0.75}{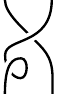} \;,
		\end{align*}
which finishes the proof.
\end{proof}

The explicit computation of $\mathrm{RT}_{\mc B}\circ Y^\bullet$ essentially follows from Lemma~\ref{calculinvisotopie}:
\begin{Lemma}\label{finaltouch}
	The cocyclic $\kk$-module $\mathrm{RT}_{\mc B}\circ Y^\bullet$ equals to the cocyclic $\kk$-module given by the family \smash{$\bigl\{\Hom_{\mc B}\bigl(\uu,\coend^{\tens n+1}\bigr)\bigr\}_{n\in\N}$}, equipped with the cofaces $\{\mathrm{RT}_{\mc B}(Y^\bullet(\delta_i^n))\}_{n\in \N^*, 0\le i \le n}$, codegene\-racies \smash{$\bigl\{\mathrm{RT}_{\mc B}\bigl(Y^\bullet\bigl(\sigma_j^n\bigr)\bigr)\bigr\}_{n\in \N, 0 \le j \le n}$}, and cocyclic operators $\{\mathrm{RT}_{\mc B}(Y^\bullet(\tau_n))\}_{n\in \N}$ given by formulas\looseness=1
	\begin{gather*}
\delta_0^n(f)=
\put(23,-2.5){\small $\cdots$}
\put(19,-15){\small $f$}
\put(15,25){\small $0$}
\put(28,25){\small $n-1$}
\,\raisebox{-8mm}{\includegraphics[scale=1]{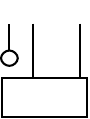}}\,,
\qquad
\delta_i^n(f)=
\put(11,-2.5){\small $\cdots$}
\put(56,-2.5){\small $\cdots$}
\put(36,-15){\small $f$}
\put(3.5,25){\small $0$}
\put(19,25){\small $i-1$}
\put(49.5,25){\small $i$}
\put(61.5,25){\small $n-1$}
\,\raisebox{-8.2mm}{\includegraphics[scale=1]{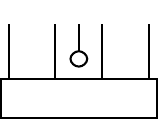}}\,,
\qquad
\delta_n^n(f)=
\put(11.5,-2.5){\small $\cdots$}
\put(19,-15){\small $f$}
\put(2,25){\small $0$}
\put(16,25){\small $n-1$}
\,\raisebox{-8mm}{\includegraphics[scale=1]{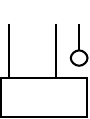}}\,,\\
\sigma_j^n(f)=
\put(23,18){\small $\cdots$}
\put(65,18){\small $\cdots$}
\put(45.5,-15){\small $f$}
\put(3.5,1){\small $0$}
\put(29,1){\small $j$}
\put(63,1){\small $j+1$}
\put(90,1){\small $n+1$}
\,\raisebox{-8mm}{\includegraphics[scale=1]{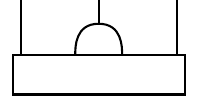}}\hspace{0.62cm} ,
\qquad
\tau_n(f)=
\put(49,-1){\small $\cdots$}
\put(45,-15){\small $f$}
\put(3,1){\small $0$}
\put(31,1){\small $1$}
\put(86,1){\small $n$}
\,\raisebox{-8.2mm}{\includegraphics[scale=1]{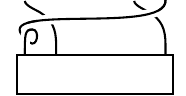}}\ .
	\end{gather*}
\end{Lemma}
\begin{proof}
All computations follow by using formula~\eqref{rtformula} from Section~\ref{sect-calc-RT-coend}, the construction of~$Y^\bullet$ from Section~\ref{construction}, and Lemma~\ref{calculinvisotopie}. Here we only provide a computation for cofaces $\{\RT_{\mc B}(Y^\bullet (\delta_i^n))\}_{n\in \N^*, 0\le i \le n}$. If $n\in \N^*$ and $1 \le i \le n-1$, then
\begin{align*}
\RT_{\mc B}(Y^\bullet (\delta_i^n))&=\Delta^{-(n+1)-(n+1)}\Hom_{\mc B}\bigl(\uu, \bigl\lvert G_{Y^\bullet (\delta_i^n)} \bigr\rvert \bigr)
 \\
&=\dim(\mc B)^{-(n+1)}\Hom_{\mc B}\left(\uu,\dim(\mc B)^{n+1}
\put(7.2,4){\small $\cdots$}
\put(61.1,4){\small $\cdots$}
\put(1,-18.5){\small $0$}
\put(13,-18.5){\small $i-1$}
\put(54.5,-18.5){\small $i$}
\put(65,-18.5){\small $n-1$}
\,\raisebox{-8mm}{\includegraphics[scale=0.75]{page34-02}}\hspace{0.3cm}\right)\\
&= \Hom_{\mc B}\left(\uu,\,
\put(7.2,4){\small $\cdots$}
\put(61.1,4){\small $\cdots$}
\put(1,-18.5){\small $0$}
\put(13,-18.5){\small $i-1$}
\put(54.5,-18.5){\small $i$}
\put(65,-18.5){\small $n-1$}
\,\raisebox{-8mm}{\includegraphics[scale=0.75]{page34-02}}\hspace{0.3cm}\right).
\end{align*}
The cases $i=0$ or $i=n$ are verified analogously.
\end{proof}

\subsection{The final step} \label{sec: finaltouch}
In Section~\ref{construction}, we constructed the cocyclic object~$Y^\bullet$ in the category of~$3$-cobordisms.
Next, in Lemma~\ref{finaltouch}, we computed the cocyclic~$\kk$-module~$\RT_{\mc B} \circ Y^\bullet$.
To prove the first part of Theorem~\ref{CYCINCOBRT}, it suffices to show that the latter is isomorphic to~$\coend^\bullet \circ \Phi$, where~$\Phi \co \Delta C \to \Delta C$ is the reindexing involution from Section~\ref{connes loday dual}. Recall that the coend~$\coend$ is a Hopf algebra in~$\mc B$, equipped with a non-degenerate Hopf pairing $\omega$. Denote its inverse by $\Omega$. The isomorphism between cocyclic~$\kk$-modules~$\RT_{\mc B} \circ Y^\bullet$ and~$\coend^\bullet \circ \Phi$ is provided by the family
\[
\omega^\bullet = \bigl\{\omega^n \co \Hom_{\mc B} \bigl(\uu, \coend^{\tens n+1}\bigr) \to \smash{\Hom_{\mc B} \bigl(\coend^{\tens n+1}, \uu\bigr)\bigr\}_{n\in \N}},
\]
 which is defined by setting for any $n\in \N$ and $f \in \Hom_{\mc B}\bigl(\uu, \coend^{\tens n+1}\bigr)$,
\begin{equation}
\label{omeganat}
\omega^n(f)=
\put(12,-12){\small $\cdots$}
\put(44,-44){\small $\cdots$}
\put(15,-25){\small $f$}
\put(1,-12){\small $0$}
\put(29,-12){\small $n$}
\,\raisebox{-16.5mm}{\includegraphics[scale=0.85]{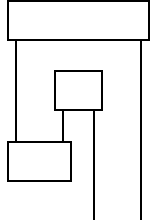}}\,.
\put(-37,4){\small $\omega$}
\put(-37,32.5){\small $\omega$}
\end{equation}
The inverse of $\omega^\bullet$ is given by the family \smash{$\Omega^\bullet = \bigl\{\Omega_n \co \Hom_{\mc B} \bigl(\coend^{\tens n+1}, \uu\bigr) \to \Hom_{\mc B} \bigl(\uu, \coend^{\tens n+1}\bigr) \bigr\}_{n\in \N}$}, which is defined by setting for any $n\in \N$ and $f \in \Hom \bigl(\coend^{\tens n+1}, \uu\bigr) $,
\[
\Omega^n(f)=
\put(12,3){\small $\cdots$}
\put(44,37){\small $\cdots$}
\put(15,17){\small $f$}
\put(1,3){\small $0$}
\put(29,3){\small $n$}
\,\raisebox{-16.5mm}{\includegraphics[scale=0.85]{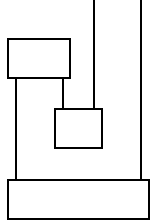}}\,.
\put(-38,-13){\small $\Omega$}
\put(-38,-41.5){\small $\Omega$}
\]
It remains to check that~$\omega^\bullet$ is a natural transformation between cocyclic~$\kk$-modules $\RT_{\mc B} \circ Y^\bullet$ and~$\coend^\bullet \circ \Phi$. The equations~\eqref{commcofac} and~\eqref{commcodegen} follow from axioms of a Hopf pairing~$\omega$ and isotopy invariance of graphical calculus. The equation~\eqref{eq-commute-cocyc} follows from equations~\eqref{invol} and \eqref{omegaS}, naturality of braiding, and isotopy invariance of graphical calculus. \hfill $\blacksquare$

\subsection[Sketch of computation of RT\_B circ Y\_\{bullet\}]{Sketch of computation of $\boldsymbol{\RT_{\mc B}\circ Y_\bullet}$} \label{RTtildeS}
Let us sketch the computation of $\RT_{\mc B}\circ Y_\bullet$, which is similar to computation of $\RT_{\mc B}\circ Y^\bullet$ given in detail in Sections~\ref{RTS} and~\ref{sec: finaltouch}.
The following lemma is analogously proved as Lemma~\ref{simpletangles}.
\begin{Lemma} \label{simpletangles2} If $T_4$ and $T_5$ are the following special ribbon graphs
	\[T_4 =\,
	\rsdraw{1.20}{0.7}{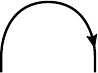}\;, \qquad T_5 =\,
	\rsdraw{0.45}{0.7}{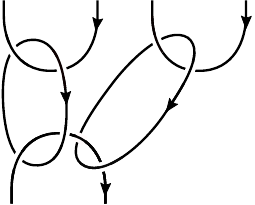}\;, \]
	then
	\begin{itemize}\itemsep=0pt
		\item[$(a)$] $\lvert T_4 \rvert= \varepsilon$,
		\item[$(b)$] $\lvert T_5 \rvert= \dim(\mc B)^{2}\Delta$.
	\end{itemize}
\end{Lemma}
Next, by using Lemma \ref{simpletangles2}, one obtains an analogue of Lemma \ref{calculinvisotopie}.

\begin{Lemma} \label{calculinvisotopie2} We have
	\begin{itemize}\itemsep=0pt
		\item[$(a)$] $
\lvert G_{Y_\bullet (d_i^n)} \rvert=\dim(\mc B)^{n}
\put(6,0){\small $\cdots$}
\put(56.5,0){\small $\cdots$}
\put(1,-23){\small $0$}
\put(34.5,-23){\small $i$}
\put(68,-23){\small $n$}
\,\raisebox{-9.5mm}{\includegraphics[scale=0.75]{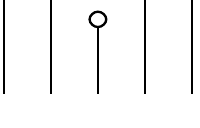}}
$,

\item[$(b)$] $
\lvert G_{Y_\bullet (s_j^n)} \rvert=\dim(\mc B)^{n+2}
\put(6,0){\small $\cdots$}
\put(74,0){\small $\cdots$}
\put(1,-23){\small $0$}
\put(43.5,-23){\small $j$}
\put(85.5,-23){\small $n$}
\,\raisebox{-9.5mm}{\includegraphics[scale=0.75]{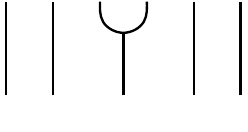}}
$,

		\item[$(c)$] $\lvert G_{Y_\bullet (t_n)} \rvert =
\begin{cases}
\dim(\mc B) \idrm_{\coend} & \mbox{if } n=0, \\
\dim(\mc B)^{n+1}
\put(24,-15){\small $\cdots$}
\put(0.8,-24){\small $0$}
\put(40,-24){\small $n-1$}
\put(70,-24){\small $n$}
\,\raisebox{-9.5mm}{\includegraphics[scale=0.75]{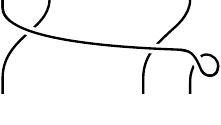}}	& \mbox{if } n\ge 1.	
\end{cases} $
	\end{itemize}
\end{Lemma}
By using Lemma~\ref{calculinvisotopie2}, one can compute that~$\RT_{\mc B}\circ Y_\bullet$ is equal to the cyclic~$\kk$-module given by the family \smash{$\bigl\{\Hom_{\mc B}\bigl(\uu,\coend^{\tens n+1}\bigr)\bigr\}_{n\in\N}$}, equipped with the faces $\{\mathrm{RT}_{\mc B}(Y_\bullet(d_i^n))\}_{n\in \N^*, 0\le i \le n}$, codege\-neracies \smash{$\bigl\{\mathrm{RT}_{\mc B}\bigl(Y_\bullet\bigl(s_j^n\bigr)\bigr)\bigr\}_{n\in \N, 0 \le j \le n}$}, and cocyclic operators $\{\mathrm{RT}_{\mc B}(Y_\bullet(t_n))\}_{n\in \N}$ given by formulas
\begin{gather*}
d_i^n(f)=
\put(15,-2.5){\small $\cdots$}
\put(60,-2.5){\small $\cdots$}
\put(41,-16){\small $f$}
\put(2,1){\small $0$}
\put(46,1){\small $i$}
\put(80,1){\small $n$}
\,\raisebox{-8.2mm}{\includegraphics[scale=1]{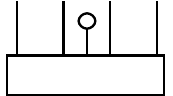}}\,,
\qquad
s_j^n=
\put(18,-2.5){\small $\cdots$}
\put(64,-2.5){\small $\cdots$}
\put(43,-16){\small $f$}
\put(2,1){\small $0$}
\put(50,1){\small $j$}
\put(87,1){\small $n$}
\,\raisebox{-8.2mm}{\includegraphics[scale=1]{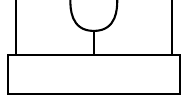}}\,,
\qquad
t_n=
\put(43,19){\small $\cdots$}
\put(41,-16){\small $f$}
\put(2,1){\small $0$}
\put(37,1){\small $n-1$}
\put(70,1){\small $n$}
\,\raisebox{-8.2mm}{\includegraphics[scale=1]{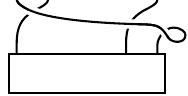}}\,.
\end{gather*}
The latter cyclic $\kk$-module is isomorphic to~$\coend_\bullet \circ \Phi^\opp$, where~$\Phi$ is the reindexing involution automorphism from Section~\ref{connes loday dual}. This again follows from the fact that the coend~$\coend$ is equipped with a non-degenerate Hopf pairing~$\omega$. Indeed, the isomorphism is given by the family~$\omega_\bullet= \omega^{\bullet}$, defined by formula~\eqref{omeganat}. This completes our proof of Theorem~\ref{CYCINCOBRT}. \hfill $\blacksquare$

\section{Related work} \label{related}
In this section we discuss some potentially related perspectives. From algebraic perspective, we review some constructions of paracyclic objects in a braided $\kk$-category and outline their connections to (co)cyclic $\kk$-modules from Section~\ref{algebraiccyclic}. These paracyclic objects can be further used to construct $r$-cyclic $\kk$-modules in the sense of Feigin and Tsygan \cite{feigin-tsygan}. From topolo\-gical perspective, we observe a relevance of the one-holed torus, which is a Hopf algebra in $\text{Cob}_3(1)^{\text{conn}}$, a cobordism category appearing in \cite{kerler} and different that the one from Section~\ref{seccobs}.

\subsection{Variations on the cyclic category}

For $r\in \N^*$, the $r$-cyclic category $\Delta C_r$ is defined in the same way as the cyclic category, except that the cyclic operators are replaced with the $r$\emph{-cyclic operators} $\{\tau_n\co n \to n\}_{n\in \N}$, which satisfy~\eqref{compcoccof}--\eqref{tnsigma0} and additionally, for any~$n\in \N$,
\begin{equation*}
 \tau_n^{r(n+1)}= \id_n.
\end{equation*}
Notice that $\Delta C_1=\Delta C$.
Another notion similar to that of the cyclic category is the paracyclic category $\Delta C_\infty$. In this setting, the cyclic operators are replaced with the \emph{paracyclic ope\-rators} $\{\tau_n\co n \to n\}_{n\in \N}$, which are isomorphisms satisfying \eqref{compcoccof}--\eqref{tnsigma0}.
Mutatis mutandi (see Section~\ref{simpcyc}) one defines $r$-(co)cyclic and para(co)cyclic objects in a category.

The cyclic duality from Section~\ref{connes loday dual} descends to the paracyclic category.
Indeed, the key ingredient of proof in \cite{loday98} is the existence of the so called extra dege\-neracies, which only relies on the invertibility of (co)cyclic operators. Since the same formulas apply, we abusively denote the paracyclic duality with $L \co \Delta C^\opp_\infty \to \Delta C_\infty$.

\subsection[Paracyclic and r-cyclic objects from braided Hopf algebras]{Paracyclic and $\boldsymbol{r}$-cyclic objects from braided Hopf algebras}
\label{parastuff}
Any Hopf algebra $H$ in a braided monoidal category $\mc B$ gives rise to a paracyclic object $\textbf{C}_\bullet(H)$ in $\mc B$, given by setting for any $n\in \N$, $\textbf{C}_n(H)=H ^{\tens n+1}$ and
\begin{equation} \label{counit-comult}
d_i^n=
\put(6,0){\small $\cdots$}
\put(56.5,0){\small $\cdots$}
\put(1,-23){\small $0$}
\put(34.5,-23){\small $i$}
\put(68,-23){\small $n$}
\,\raisebox{-9.5mm}{\includegraphics[scale=0.75]{page38-01}}\,,
\qquad
s_j^n=
\put(6,0){\small $\cdots$}
\put(74,0){\small $\cdots$}
\put(1,-22){\small $0$}
\put(43.5,-22){\small $j$}
\put(85.5,-22){\small $n$}
\,\raisebox{-9.5mm}{\includegraphics[scale=0.75]{page38-02}}\,,
\qquad
t_n=
\put(24,-15){\small $\cdots$}
\put(0.8,-24){\small $0$}
\put(40,-24){\small $n-1$}
\put(70,-24){\small $n$}
\,\raisebox{-9.5mm}{\includegraphics[scale=0.75]{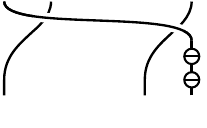}}\,.
\end{equation} This construction is functorial, that is, any morphism between coalgebras in $\mc B$ induces the morphism of corresponding paracyclic objects in $\mc B$.

Similarly, one can associate to $H$ the paracocyclic object~$\textbf{A}_\bullet(H)$ in $ \mc B$, given by setting for any~$n\in \N$, $\textbf{A}_n(H)=H^{\tens n+1}$ and
\[
\delta_i^n=
\put(7.2,0){\small $\cdots$}
\put(61.1,0){\small $\cdots$}
\put(1,-22.5){\small $0$}
\put(13,-22.5){\small $i-1$}
\put(54.5,-22.5){\small $i$}
\put(65,-22.5){\small $n-1$}
\,\raisebox{-9.3mm}{\includegraphics[scale=0.75]{page34-02}}\, ,
\qquad
\sigma_j^n =
\put(7.3,0){\small $\cdots$}
\put(72.5,0){\small $\cdots$}
\put(1,-23.5){\small $0$}
\put(33,-23.5){\small $j$}
\put(43,-23.5){\small $j+1$}
\put(76,-23.5){\small $n+1$}
\,\raisebox{-9.8mm}{\includegraphics[scale=0.75]{page34-04}}\,,
\qquad
\tau_n =
\put(40,-12){\small $\cdots$}
\put(2.8,-22){\small $0$}
\put(19.5,-22){\small $1$}
\put(69.5,-22){\small $n$}
\,\raisebox{-8.8mm}{\includegraphics[scale=0.75]{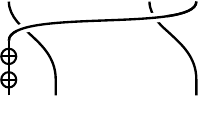}}\,.
\]
This construction is functorial, that is, any algebra morphism in $\mc B$ induces the morphism of corresponding paracocyclic objects in $\mc B$.

Notice that if~$\mc B$ is a ribbon category and if~$H$ is an involutive Hopf algebra in~$\mc B$ (that is,~$S^2=\theta_H$, where $\theta$ is the canonical twist of $\mc B$), then the paracocyclic operators~$\{t_n\}_{n\in \N}$ of~$\textbf{C}_\bullet(H)$ satisfy the ``twisted cyclicity condition'', that is, for all $n\in \N$,
\begin{equation}\label{twistedcyc}
t_n^{n+1}=(\theta_{H^{\tens n+1}})^{-1}.
\end{equation}
If $\mc B$ is additionally $\kk$-linear, then the cocyclic $\kk$-module $H^\bullet$ from Section~\ref{objetdedebut} can be obtained by compo\-sing~$\textbf{C}_\bullet(H)$ with the hom-functor~$\Hom_{\mc B}(-,\uu)$. Here, the cocyclicity condition~\eqref{cocyclicity} for~$H^\bullet$ follows by naturality of twists and the fact that~$\theta_\uu=\id_\uu$. In this vein, one derives from~$\textbf{C}_\bullet(H)$ an~$r$-cocyclic $\kk$-module as follows. Namely, if~$i$ is a simple object of~$\mc B$, then the twist~$\theta _i$ is a scalar multiple of the identity morphism. If this scalar is of finite order~$r$, then the composition of~$\textbf{C}_\bullet(H)$ with the hom-functor~$\Hom_{\mc B}(-,i)$ induces an~$r$-cocyclic (respectively, cyclic)~$\kk$-module. In a similar way, one obtains~$r$-cyclic modules from~$\textbf{A}_\bullet(H)$. For instance, according to \cite{BJT}, and built on theorems of Vafa~\cite{VAFA1988421} and M\"{u}ger~\cite{MUGER2003159}, all the twists on simple objects of a $\mathbb{C}$-linear ribbon fusion category are roots of unity.

Finally, we point out that the above construction of~$\textbf{C}_\bullet(H)$ (respectively,~$\textbf{A}_\bullet(H))$ may be restated for coalgebra~$H$ (respectively, algebra) in a balanced category $\mc B$, just like we did in Section~\ref{algebraiccyclic}. Namely, in the above formulas, the square of the antipode should be replaced with the twists. As already remarked in Section~\ref{passing}, by passing to $\textbf{C}_\bullet(H)\circ L$ (respectively,~$\textbf{A}_\bullet(H)\circ L^\opp$), these constructions fit into a more general framework of Akrami and Majid in \cite{cycliccocycles}, who considered ribbon algebras (or dually, coribbon coalgebras) in a braided category.

\subsection[Paracyclic and r-cyclic objects from Crane--Yetter Hopf algebra]{Paracyclic and $\boldsymbol{r}$-cyclic objects from Crane--Yetter Hopf algebra} \label{cy para}
Another relevant category is that of connected $3$-cobordisms, further denoted as $\conncob$, which first appeared in \cite{kerler}, which was studied over the last three decades, notably in \cite{asaeda, bobtcheva-piergallini, Ha, kerler2}, and most recently in \cite{beliakova-derenzi}.
For each $g\in \N$, fix a compact connected oriented surface $\Sigma_{g,1}$ of genus~$g$ and with one boundary component.
The objects of~$\conncob$ are surfaces~$\Sigma_{g,1}$ for each~$g\in \N$.
A morphism~$\Sigma_{g,1} \to \Sigma_{h,1}$ is given by the connected cobordism from~$\Sigma_{g,1}$ to~$\Sigma_{h,1}$.
In contrast to the category~$\textbf{3}\Cob_0$ which is recalled in Section~\ref{sect-cob3} and used throughout the paper, the category~$\conncob$ is a non-symmetric braided monoidal category, the monoidal product on objects is given by the connected sum of surfaces, on morphisms it is given by the connected sum of cobordisms and the unit object is $\Sigma_{0,1}$.

It is a result of Crane and Yetter \cite{Crane-Yetter}, that the one-holed torus $\Sigma_{1,1}$ has a structure of a Hopf algebra in ${\conncob}$.
By the general construction from Section~\ref{parastuff}, one can organize the family of surfaces
\[
\smash{\bigl\{\Sigma_{1,1}^{\tens g}\bigr\}_{g\in \N^*}=\{\Sigma_{g,1}\}_{g\in \N^*}}
\]
into a paracyclic object $\textbf{C}_\bullet(\Sigma_{1,1})$ (respectively, paracocyclic object $\textbf{A}_\bullet(\Sigma_{1,1}))$ in~${\conncob}$.
The exi\-stence of~$\textbf{C}_\bullet(\Sigma_{1,1})$ (respectively,~$\textbf{A}_\bullet(\Sigma_{1,1}))$ implies that any braided monoidal functor from ${\conncob}$ to a braided monoidal category~$\mc{B}$ induces a para(co)cyclic object in~$\mc B$. An example of such functors was recently studied by Beliakova and De Renzi in \cite{beliakova-derenzi}.
Let~$\mc B$ be a ribbon unimodular facto\-ri\-zable category and~$\kk$ an algebraically closed field.
Here the unimodularity means that~$\mc B$ is a finite~$\kk$-category in the sense of~\cite{etingof2016tensor} and that the projective cover of~$\uu$ is self-dual.
By~\cite[Theorem~1.2]{beliakova-derenzi}, there is a braided monoidal functor~$J_3 \co \conncob \to \mc B$, which, in particular, sends the one-holed torus to the end~$\aend$ of category~$\mc B$, that is, the end of functor~$(X,Y)\mapsto X\tens Y^*$.
Furthermore,
$J_3\circ \textbf{C}_\bullet(\Sigma_{1,1})=\textbf{C}_\bullet(\aend)$ and $J_3\circ \textbf{A}_\bullet(\Sigma_{1,1})=\textbf{A}_\bullet(\aend)$.
Note that factorizability in the sense of~\cite{beliakova-derenzi} is equivalent to non-degeneracy of the canonical pairing of the coend of~$\mc B$.
 The latter is also equivalent to invertibility of the Drinfeld map from~\cite[Proposition~4.11]{FGR}. In this setting, the end~$\aend$ and the coend~$\coend$ are isomorphic Hopf algebras. Therefore, by functoriality, paracyclic objects~$\textbf{C}_\bullet(\aend)$ and~$\textbf{C}_\bullet(\coend)$ (respectively, paracocyclic objects $\textbf{A}_\bullet(\aend)$ and $\textbf{A}_\bullet(\coend)$) in~$\mc B$ are isomorphic. Finally, since~$S^2_\aend= \theta_{\aend}$, paracyclic operators of~$\textbf{C}_\bullet(\aend)$ (respectively, paracocyclic operators of~$\textbf{A}_\bullet(\aend)$) satisfy the twisted cyclicity (respectively, cocyclicity) condition~\eqref{twistedcyc}. Hence, by composing~$J_3\circ \textbf{C}_\bullet(\Sigma_{1,1})$ with appropriate hom-functors (see Section~\ref{parastuff}), one obtains $r$-(co)cyclic $\kk$-modules.

\begin{Remark} Cobordisms in~${\conncob}$ admit surgery presentation by certain tangles (for a review, see~\cite[Section~4]{beliakova-derenzi}). On one hand, there is a resemblance between tangles which present structural morphisms of the Hopf algebra~$\Sigma_{1,1}$ and the special ribbon graphs which present generating morphisms of (co)cyclic objects in~$\textbf{3}\Cob_0$ (see Theorem~\ref{CYCINCOB} and a relevant Remark \ref{drinfeld-map-cobordisms}). Note also that in the construction of the latter objects, we do not use the monoidal structure of~$\textbf{3}\Cob_0$. One the other hand, in the contrast with the Kirby calculus on special ribbon graphs from Section~\ref{kirbycalculus}, \textbf{COUPON} and~\textbf{TWIST} move do not appear as moves between tangles presenting morphisms in~${\conncob}$.
It would be interesting to explore relationships between (co)cyclic objects in~$\textbf{3}\Cob_0$ which were constructed in Section~\ref{cobcyclic} and para(co)cyclic objects in~${\conncob}$ from Section~\ref{cy para}.
\end{Remark}

\begin{Remark} \label{drinfeld-map-cobordisms} A reviewer helpfully pointed out that there are two isomorphic algebra structure on the one-holed torus $\Sigma_{1,1}$. The structural morphisms of the Hopf algebra which ``models'' the end $\aend$ are defined by
\begingroup
\allowdisplaybreaks
\begin{gather*}
 m= \,
\rsdraw{0.55}{0.75}{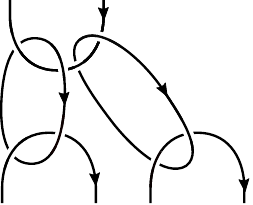} \;, \qquad \eta = \,
\rsdraw{0.55}{0.75}{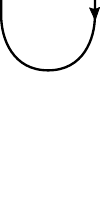} \;, \qquad \Delta = \,
\rsdraw{0.55}{0.75}{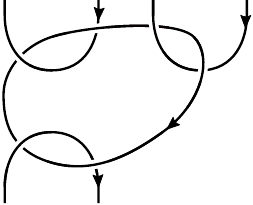} \;, \\
\varepsilon= \,
\rsdraw{0.55}{0.75}{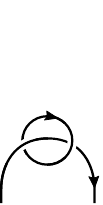} \;, \qquad S= \,
\rsdraw{0.55}{0.75}{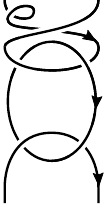} \;.
\end{gather*}
\endgroup
On the other hand, the Hopf algebra structure which ``models'' the coend $\coend$ is defined as follows.
The multiplication and the unit are defined respectively by tangles $T_3$ and $T_2$ from Lemma \ref{simpletangles}, while the comultiplication and the counit are defined by tangles $T_5$ and $T_4$ from Lemma \ref{simpletangles2}.
 Its antipode is given by the tangle \[ \,
\rsdraw{0.55}{0.75}{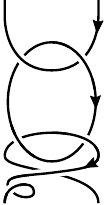} \;.\]
The two Hopf algebra structures are isomorphic via the map
 \[ \,
\rsdraw{0.55}{0.75}{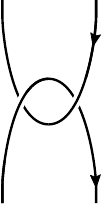} \qquad \text{of inverse} \qquad
\rsdraw{0.55}{0.75}{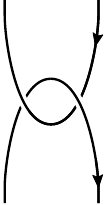} \;.\]
By functoriality of construction from Section~\ref{parastuff}, the associated para(co)cyclic objects are also isomorphic.
\end{Remark}

\begin{Remark} It is remarked in \cite{bartulovic2022cyclicsetsribbonstring} that (co)cyclic sets from ribbon string links, which are recalled in the introduction are isomorphic, via an appropriate bijection, to certain (co)cyclic sets from ribbon handles as in \cite{bv} (called bottom tangles in \cite{Ha}).
The latter comes from external Hopf algebra construction by Habiro \cite{Ha}. Namely, the object $\uparrow \downarrow$ in the category of the oriented ribbon tangles \cite[Section~\RomanNumeralCaps{1}.2.3]{turaevqinvariants} has a coalgebra structure via
\[\Delta = \,
\rsdraw{0.55}{0.75}{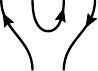} \;, \qquad \varepsilon= \,
\rsdraw{0.55}{0.75}{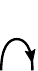}. \;\]
This fits into the context of special graphs as in Sections \ref{seccobs} and \ref{cobcyclic}, or tangles as in (the present) Section~\ref{related}.
Using above $\Delta$ and $\varepsilon$ and replacing $S^{-2}$ by $\theta_{\uparrow \downarrow}$ in formula \eqref{counit-comult}, we obtain a paracyclic object $\uparrow \downarrow_\bullet$ in the category of oriented ribbon tangles $\mc G$, where graphical presentation of operators of structure reminds of the informal discussion in the introduction of \cite{cycliccocycles}.
If we further restrict ourselves to the category $\mc G_+$ of special ribbon graphs modulo moves~$\textbf{K}1$,~$\textbf{K}2'$,~\textbf{COUPON}, and the~\textbf{TWIST} move, we obtain the defining special ribbon graphs from Section~\ref{Stildesketch} by postcomposing the structure operators $d_i^n$, $s_j^n$ and $t_n$ of $\uparrow \downarrow_\bullet$ with the diagrams representing the identity cobordisms (see Figure~\eqref{eq-id-graph}). More precisely,
\begingroup
\allowdisplaybreaks
\begin{gather*}
G_{Y_\bullet (d_i^n)} =
\put(52.5,16){\small $\cdots$}
\put(15.5,-32){\small $\cdots$}
\put(87,-32){\small $\cdots$}
\put(54.5,-1){\small $I_n$}
\put(5.5,-41){\small $0$}
\put(54.5,-41){\small $i$}
\put(103.5,-41){\small $n$}
\,\raisebox{-16.1mm}{\includegraphics[scale=1]{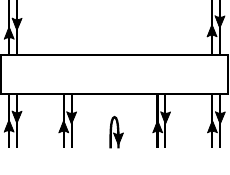}}\,,
\qquad
G_{Y_\bullet (s_j^n)}=
\put(52.5,16){\small $\cdots$}
\put(15.5,-32){\small $\cdots$}
\put(87,-32){\small $\cdots$}
\put(47,-1){\small $I_{n+2}$}
\put(5.5,-41){\small $0$}
\put(54.5,-41){\small $j$}
\put(103.5,-41){\small $n$}
\,\raisebox{-16.1mm}{\includegraphics[scale=1]{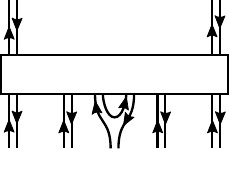}}\,,\\
G_{Y_\bullet (t_n)}=
\put(45,16){\small $\cdots$}
\put(37,-32){\small $\cdots$}
\put(42,-1){\small $I_{n+1}$}
\put(5.5,-41){\small $0$}
\put(60,-41){\small $n-1$}
\put(92,-41){\small $n$}
\,\raisebox{-15.7mm}{\includegraphics[scale=1]{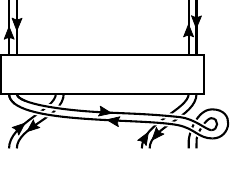}}\,.
\end{gather*}
\endgroup
In this setup, one can endow the object $\uparrow \downarrow$ with an \emph{external Hopf algebra} structure \cite[Sections~6.2 and 6.3]{Ha}.
The multiplication, the unit, and the antipode are here defined by setting for all admissible~$(i+j,k)$-graphs $\Gamma$,
\begingroup
\allowdisplaybreaks
\begin{gather*}
\check{m}_{(i,j,k)}(\Gamma) =
\put(52.5,16){\small $\cdots$}
\put(15.5,-32){\small $\cdots$}
\put(87,-32){\small $\cdots$}
\put(40,0){\small $\Gamma$}
\put(5.5,-41){\small $0$}
\put(62,-41){\small $i$}
\put(96.5,-41){\small $i+j$}
\put(13,22){\small $1$}
\put(112,22){\small $k$}
\,\raisebox{-16mm}{\includegraphics[scale=1]{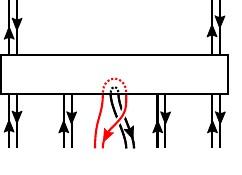}}\,,
\qquad
\check{u}_{(i,j,k)}(\Gamma) =
\put(44,16){\small $\cdots$}
\put(15.5,-32){\small $\cdots$}
\put(72,-32){\small $\cdots$}
\put(47,0){\small $\Gamma$}
\put(5.5,-41){\small $0$}
\put(47,-41){\small $i$}
\put(81.5,-41){\small $i+j$}
\put(13,22){\small $1$}
\put(97,22){\small $k$}
\,\raisebox{-16mm}{\includegraphics[scale=1]{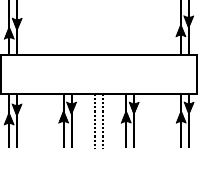}}\,,\\
\check{S}_{(i,j,k)}(\Gamma) =
\put(44,16){\small $\cdots$}
\put(15.5,-32){\small $\cdots$}
\put(72,-32){\small $\cdots$}
\put(47,0){\small $\Gamma$}
\put(5.5,-41){\small $0$}
\put(47,-41){\small $i$}
\put(81.5,-41){\small $i+j$}
\put(13,22){\small $1$}
\put(97,22){\small $k$}
\,\raisebox{-16mm}{\includegraphics[scale=1]{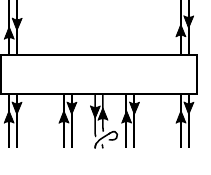}}\,.
\end{gather*}
\endgroup
 The red component in the formula for $\check m_{(i,j,k)}$ comes from duplicating (along the framing) the arc connecting input $2i+2$ with $2i+3$. The image $\check u_{(i,j,k)} (\Gamma)$ is obtained by deleting the arc connecting input $2i+2$ with $2i+3$ and the image \smash{$\check S_{(i,j,k)} (\Gamma)$} is obtained by changing the orientation of the arc connecting input $2i+2$ with $2i+3$.
 Note that \smash{$\check S^2_{(i,j,k)} (\Gamma) = \Gamma \circ \bigl( (\uparrow \downarrow)^{\tens i} \tens \theta_{\uparrow \downarrow} \tens (\uparrow \downarrow)^{\tens j} \bigr)$}. Using this data, we recover special ribbon graphs appearing in Section~\ref{construction}:
 \begin{gather*}
 G_{Y^\bullet (\delta_i^n)} = \check{u}_{(i,n-i, n+1)} (I_{n+1}), \qquad G_{Y^\bullet (\sigma_j^n)} = \check{m}_{(j,n-j, n+1)} (I_{n+1}),\\
G_{Y^\bullet (\tau_n)}=
\put(49,16){\small $\cdots$}
\put(59,-32){\small $\cdots$}
\put(45,-1){\small $I_{n+1}$}
\put(4.5,-41.5){\small $0$}
\put(31,-41.5){\small $1$}
\put(94,-41.5){\small $n$}
\,\raisebox{-16mm}{\includegraphics[scale=1]{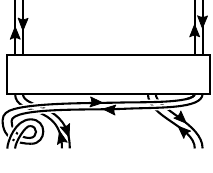}}\,.
 \end{gather*}
Let us illustrate the equality for codegeneracies in a special case:
\[
\check{m}_{(0,0,1)} =
\put(71,-42){\small $0$}
\,\raisebox{-16mm}{\includegraphics[scale=0.75]{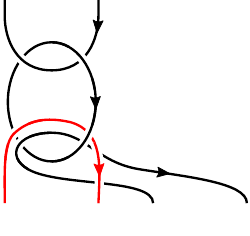}}\,=
\put(70,-42){\small $0$}
\,\raisebox{-16mm}{\includegraphics[scale=0.75]{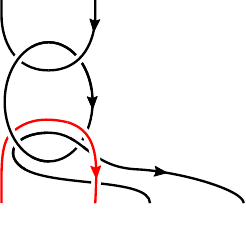}}\,=
\put(17,-42){\small $0$}
\put(71,-42){\small $1$}
\,\raisebox{-16mm}{\includegraphics[scale=0.75]{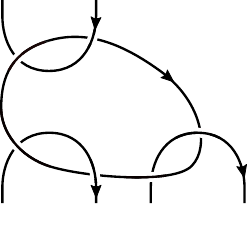}}\,,
\]
where the latter indeed equals~$G_{Y^\bullet (\sigma_0^0)}$. This gives a more formal construction of the (co)cyclic object in~$\mc G_+$. In both cases, if we omit \textbf{COUPON} and \textbf{TWIST} moves, we merely obtain para(co)cyclic objects.
\end{Remark}

\subsection*{Acknowledgements}

Most of the content of this paper is part of my PhD thesis. I would like to thank my supervisor Alexis Virelizier for his comments and advice. I would also like to thank Ulrich Kr\"{a}hmer, Christoph Schweigert, Kenichi Shimizu, and anonymous referees for discussions and comments. Finally, I~am grateful to the Laboratoire Paul Painlev\'e of the University of Lille and the Institut für Geometrie of the TU Dresden for their hospitality. This work was supported by the Labex CEMPI (ANR-11-LABX-0007-01), by the R\'egion Hauts-de-France, and by the FNS-ANR OChoTop grant (ANR-18-CE93-0002-01).

\pdfbookmark[1]{References}{ref}
\LastPageEnding

\end{document}